%% file: manuscript.tex
\documentclass[a4paper,12pt,twoside]{article}
\usepackage[latin1]{inputenc}
\usepackage{amsfonts,amssymb,amsmath,exscale}
\usepackage[dvips]{graphicx,color,psfrag}
\usepackage{mathtools}
\usepackage{fancyhdr}
\usepackage{breakcites}
\usepackage{rotating}
\usepackage{pifont}
\usepackage[hyperref]{xcolor}
\definecolor{myblue}{RGB}{0,0,139}
\usepackage[colorlinks=true,citecolor=blue,urlcolor=blue,linkcolor=blue]{hyperref}
\hypersetup{bookmarksnumbered}
\usepackage[square,numbers]{natbib}
\usepackage{tensor}
\usepackage{setspace}
\usepackage{bbm}
\usepackage[font=small,labelfont=bf]{caption}
\usepackage{sectsty}
\sectionfont{\fontsize{15}{0}\selectfont}
\subsectionfont{\fontsize{12}{0}\selectfont}
\fancypagestyle{plain}{%
  \fancyhead{}%
  \fancyfoot[c]{\sffamily\thepage}%
}
\makeatletter
\def\cleardoublepage{\clearpage\if@twoside \ifodd\c@page\else
  \hbox{}
  \vspace*{\fill}
  \thispagestyle{empty}
  \newpage
  \if@twocolumn\hbox{}\newpage\fi\fi\fi}
\usepackage{algorithm2e}

\usepackage{amsmath,amssymb}
\renewcommand{\citep}[1]{(\cite{#1}{})}

\renewcommand{\thefootnote}{$^*$}
\newcommand{\B}{{\cal B}}
\DeclareOldFontCommand{\sf}{\normalfont\sffamily}{\mathsf}
\DeclareOldFontCommand{\tt}{\normalfont\ttfamily}{\mathtt}
\DeclareOldFontCommand{\bf}{\normalfont\bfseries}{\mathbf}
\DeclareOldFontCommand{\it}{\normalfont\itshape}{\mathit}
\DeclareOldFontCommand{\sl}{\normalfont\slshape}{\relax}
\DeclareOldFontCommand{\sc}{\normalfont\scshape}{\relax}
\DeclareRobustCommand*{\cal}{\@fontswitch{\relax}{\mathcal}}
\DeclareRobustCommand*{\mit}{\@fontswitch{\relax}{\mathnormal}}
\DeclareMathAlphabet{\IKbb}{U}{bbm}{m}{sl}
\DeclareMathAlphabet{\Ikbb}{U}{bbmss}{m}{it}
\DeclareRobustCommand{\KIH}{{\IKbb H}}
\DeclareRobustCommand{\KIM}{{\IKbb M}}

\DeclareRobustCommand{\BC}{{\boldsymbol{\mathnormal C}}}

\DeclareRobustCommand{\BF}{{\boldsymbol{\mathnormal F}}}

\DeclareRobustCommand{\BN}{{\boldsymbol{\mathnormal N}}}

\DeclareRobustCommand{\BP}{{\boldsymbol{\mathnormal P}}}
\DeclareRobustCommand{\BQ}{{\boldsymbol{\mathnormal Q}}}

\DeclareRobustCommand{\BT}{{\boldsymbol{\mathnormal T}}}

\DeclareRobustCommand{\BX}{{\boldsymbol{\mathnormal X}}}

\DeclareRobustCommand{\Ba}{{\boldsymbol{\mathnormal a}}}

\DeclareRobustCommand{\Bk}{{\boldsymbol{\mathnormal k}}}

\DeclareRobustCommand{\Bn}{{\boldsymbol{\mathnormal n}}}

\DeclareRobustCommand{\Bx}{{\boldsymbol{\mathnormal x}}}

\DeclareRobustCommand{\calB}{{\mathcal B}}
\DeclareRobustCommand{\calD}{{\mathcal D}}
\DeclareRobustCommand{\calH}{{\mathcal H}}

\DeclareRobustCommand{\calP}{{\mathcal P}}
\DeclareRobustCommand{\calR}{{\mathcal R}}
\DeclareRobustCommand{\calT}{{\mathcal T}}
\DeclareRobustCommand{\calW}{{\mathcal W}}
\DeclareRobustCommand{\R}{{\mathcal R}}

\DeclareMathAlphabet{\Inbb}{U}{bbmss}{m}{n}

\DeclareMathAlphabet{\gothic}{U}{euf}{m}{n}

\DeclareMathAlphabet{\Bgothic}{U}{euf}{b}{n}

\newcommand{\tr}{\mathop{\operator@font tr}}
\newcommand{\dev}{\mathop{\operator@font dev}}

\DeclareMathAlphabet{\Ibb}{U}{msb}{m}{n}

\DeclareRobustCommand{\IN}{{\Ibb N}}
\DeclareRobustCommand{\IA}{{\Ibb A}}
\DeclareRobustCommand{\IQ}{{\Ibb Q}}
\DeclareRobustCommand{\IL}{{\Ibb L}}
\DeclareRobustCommand{\IP}{{\Ibb P}}

\newcommand{\een}{\]\@ignoretrue}
\DeclareRobustCommand{\Bone}{{\boldsymbol{\mathit 1}}}
\DeclareRobustCommand{\Bzero}{{\boldsymbol{\mathit 0}}}

\DeclareRobustCommand{\Bvarphi}{{\boldsymbol{\varphi}}}
\DeclareRobustCommand{\Blambda}{{\boldsymbol{\lambda}}}

\newcommand{\Arg}{\mathop{\operator@font Arg}}                
\newcommand{\Div}{\mathop{\operator@font Div}}                
\renewcommand{\div}{\mathop{\operator@font div}}              
\newcommand{\GRAD}{\mathop{\operator@font GRAD}}              
\newcommand{\Grad}{\mathop{\operator@font Grad}}              
\newcommand{\grad}{\mathop{\operator@font grad}}              
\newcommand{\WITH}{\quad\mbox{with}\quad}
\newcommand{\AND}{\quad\mbox{and}\quad}
\newcounter{SetOfEq}

\newcommand{\assemble}{\mathop{\textbf{\Large{\sf A}}}}

\newcommand{\lbar}{\overset{\rule[-.05 ex]{1.3ex}{.075ex}}}
\newcommand{\llbar}{\overset{\rule[-.05 ex]{2.5ex}{.075ex}}}
\newcommand{\ltilde}[1]{\widetilde{#1}}
\newcommand{\lhat}[1]{\widehat{#1}}

\newcommand\SFB[1]{\boldsymbol{\mathsf{#1}}}

\DeclareRobustCommand{\Bsfd}{\underline{\boldsymbol{\mathsf{d}}}}
\DeclareRobustCommand{\BsfK}{\underline{\boldsymbol{\mathsf{K}}}}
\DeclareRobustCommand{\BsfR}{\underline{\boldsymbol{\mathsf{R}}}} 
    \newcommand{\jump}{\jumpi}              
\def\jumpi#1{%
   \setbox8\hbox{$\left[\llap{$#1$}\right[$}%
   \left[\kern-0.3\wd8\left[#1\right]\kern-0.3\wd8\right]%
}

\begin{document}

\thispagestyle{empty}

\begin{center}
\textbf{\large Swelling-induced pattern transformations of periodic hydrogels -- from the wrinkling of internal surfaces to the buckling of thin films}

\medskip

E. Polukhov, L. Pytel \& M.-A. Keip%
\footnote{Corresponding author: marc-andre.keip@mechbau.uni-stuttgart.de, Phone: +49 711 685 66233, Fax: +49 711 685 66347}

\medskip

Institute of Applied Mechanics \\
Department of Civil and Environmental Engineering \\
University of Stuttgart, Stuttgart, Germany
\end{center}
\medskip

{%
\footnotesize
\noindent \textbf{Abstract.} 
We investigate pattern transformations of periodic hydrogel systems that are triggered by swelling-induced structural instabilities. The types of microstructures considered in the present work include single-phase and two-phase voided hydrogel structures as well as reinforced hydrogel thin films. While the observed transformations of the single-phase structures show good agreement with experimental findings, the two-phase materials provide novel patterns associated with wrinkling of internal surfaces. Furthermore, an extensive parametric study on the reinforced hydrogel thin films reveals new opportunities for the design of complex out-of-plane surface modes caused by swelling-induced instabilities. Next to the mentioned buckling-type instabilities, we encountered the development of micro-creases at the internal surfaces of periodic media before the loss of strong ellipticity of effective moduli.

\noindent \textbf{Keywords.} 
periodic hydrogel structures,
pattern transformations,
buckling,
Bloch-Floquet analysis,
swelling-induced instabilities,
strong ellipticity
}

\pagestyle{fancy}
\pagenumbering{arabic}

\section{Introduction}
\label{intro10}
\renewcommand{\thefootnote}{\arabic{footnote}}
\setcounter{footnote}{0} 

Hydrogels are polymeric, soft and biocompatible materials that undergo large deformations under swelling \cite{hoffman12}. At swollen state, their mechanical response depends on both the properties of the polymer and the solvent phase. While the elastic stiffness correlates with the density of cross-links of the polymer network, viscoelasticity is driven by the viscosity of the solvent \cite{bai+etal19}. Depending on how the polymer network is held together, hydrogels can be categorized as physical (reversible) or chemical (permanent), see \cite{rosiak+yoshii99} and \cite{hoffman12} for specific examples related to each category. In the former group, the network is interconnected by polymer-chain entanglements as well as possible additional forces resulting from ionic or hydrophobic interactions and hydrogen bonding \cite{zhang+etal18}. The latter group of hydrogels is based on covalently cross-linked networks. Both physical and chemical hydrogels can be inhomogeneous due to clustered regions of cross-links and swelling, as well as the existence of defects in the polymer network. Furthermore, hydrogels can be classified as conventional or (multi-)stimuli responsive \cite{hoffman12}. While both are hydrophilic and undergo swelling in a solvent, (multi-)stimuli responsive hydrogels may also respond to changes in pH, temperature or electric field \cite{hoffman95}. In the present work, we will focus on conventional hydrogels and describe their chemo-mechanical response related to the diffusion of a solvent through their polymer network.

\subsection{Engineering and biomedical applications of hydrogels}
%
Since hydrogels have favorable properties and show tunable, versatile response, they have a wide range of applications, such as fertilizer encapsulaters in agriculture \cite{rudzinski+etal02},
soft actuators \cite{ionov14} in engineering as well as
artificial muscles \cite{kaneko+etao02},
tissue-engineered auricles \cite{bichara+etal12},
contact lenses and wound dressings \cite{calo+khutoryanskiy15} in biomedicine, see also
\cite{rosiak+yoshii99,bichara+etal10,hoffman12} and the references therein. By virtue of intricate fabrication techniques, further applications are being realized not only at macroscopic scale, but also at micro- and nanoscopic scales \cite{ma+etal16,savina+etal16}. Here, we mention in particular the fabrication by means of photolithography, which has attracted increased attention for the manufacturing of patterned and composite hydrogels exhibiting controllable deformations \cite{ma+etal16,li+etal19}. We refer to \cite{wang+etal17,ma+etal19} for the use of photolithography to fabricate composite periodic hydrogel structures composed of high-swelling and non-swelling regions. By means of the latter, the authors were able to induce out-of-plane buckling in water yielding various periodic buckling patterns. Consequently, by designing associated hydrogel microstructures and by tuning material parameters, a rich set of deformation modes can be activated which increases their potential application as soft devices exploiting chemo-mechanical interactions.

Besides pattern transforming instabilities in periodic systems, surface wrinkling and creasing are observed in hydrogel structures. These instabilities can be exploited to further advance engineering and biomedical applications. In particular, wrinkling of hydrogels can be utilized for controlled formation of microgears \cite{yin+etal09}, generation of multicellular spheroids \cite{zhao+etal14}, controlled cell spreading \cite{jiang+etal02} among many other possibilities, see also \cite{lam+etal06,guvendiren+etal10c} and the detailed review \cite{dervaux+amar12} on instabilities in gels. We further refer to the recent experimental realization of nano-wrinkled architectures by means of laser direct assembly \cite{fan20223d}.

\subsection{Continuum-mechanical models for hydrogels}
%
Not least because of the above mentioned favorable properties and applications, hydrogels have gained increased attention in theoretical and computational research. Here, one main goal is to develop versatile and reliable models to simulate their response. Related models can give further insights into underlying features, provide access to the design of pertinent microstructures for specific applications and allow the investigation of the influence of various parameters on the respective behavior. In particular continuum-mechanical models are established and widely used to capture the response of hydrogels at macroscopic as well as microscopic scales. These models take into account chemo-mechanical interactions within the polymer network of a hydrogel, as well as between the polymer network and the diffusing solvent. We refer to \cite{zhao+hong+suo08,hong+zhao+zhou+suo08} for continuum-mechanical models of gels based on Flory--Huggins theory \cite{flory53,huggins41}, see also the Flory--Rehner type models \cite{flory+rehner43a,flory+rehner43b}. Finite-element implementations of the models \cite{zhao+hong+suo08,hong+zhao+zhou+suo08} have been applied to the simulation of transient diffusion processes in \cite{zhang+zhao+suo+jiang09} and drying-induced instabilities in \cite{hong+zhao+suo08}. We further refer to the continuum approaches addressed in \cite{chester+anand10,chester+dileo+anand15}.

Note that the finite-element formulation of the above mentioned models has been based on the chemical potential as a primary unknown. In combination with the mechanical deformation map, this corresponds to a saddle-point formulation. Because of the a priori indefiniteness of the underlying system of equations, related implementations could result in numerical challenges associated with \emph{inf-sup} \emph{instabilities} \cite{teichtmeister+etal19}. Such challenges can be overcome by employing a minimization-based variational formulation as proposed by \cite{boeger+etal17a} for hydrogels; see \cite{miehe+etal15} for a basis rooted in poromechanics.

A further similarity of the above mentioned approaches is the use of Biot's theory \cite{biot41,coussy04} relying on a Flory--Rehner-type energy function and Fick's law of diffusion. Alternative continuum-mechanical models have been established on the basis of the so-called Theory of Porous Media \cite{ehlers+acarturk09,ehlers+etal10,ehlers+wagner19}. We further refer to recent multiscale modeling approaches to describe the diffusion process in heterogeneous microstructures \cite{larsson+etal10,kaessmair+steinmann18,polukhov+keip20,pollmann+etal21,kaessmair+etal21}.

\subsection{Instabilities of hydrogels and their prediction}
%
\begin{figure}%
  \centering%
  \footnotesize%
  \psfrag{a}[c][c]{\shortstack{wrinkling}}
  \psfrag{b}[c][c]{\shortstack{creasing}}
  \psfrag{c}[c][c]{\shortstack{diamond-plate pattern}}
  \psfrag{d}[c][c]{\shortstack{initial \\ configuration}}
  \psfrag{e}[c][c]{\shortstack{deformed \\ configuration}}
  \psfrag{f}[cc][cc]{\shortstack{solvent diffuses \\ through external surface}}
  \psfrag{h}[cc][cc]{\shortstack{solvent diffuses \\ through internal surfaces}}
  \psfrag{unit}[c][c]{unit-cell RVE}
  \includegraphics*[scale=1.4]{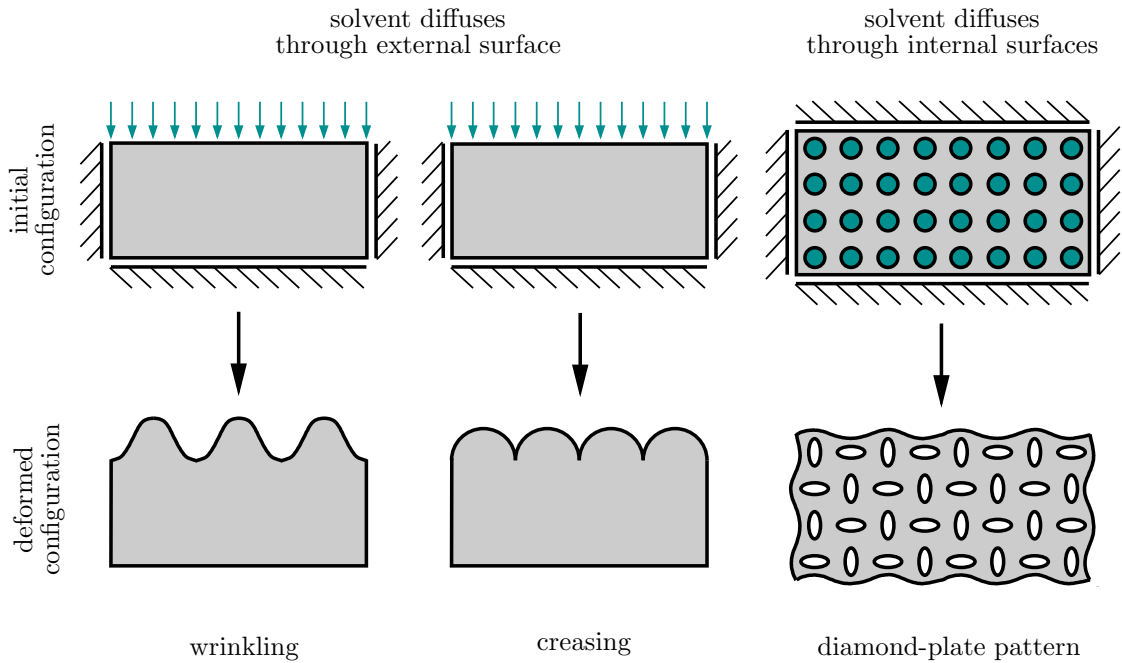}
\caption[Description of computational homogenization]{%
\textit{%
Examples of swelling-induced instabilities in constrained hydrogel structures}.}
\label{fig.intro1}%
\end{figure}%
Continuum-mechanical models of hydrogels have further been exploited to investigate complex instability phenomena such as surface wrinkling, creasing and pattern transformations as a result of deformation-diffusion processes \cite{dervaux+amar12} (see Fig.~\ref{fig.intro1} for an illustration).
Experimental investigations of surface instabilities have been documented, for example, in  \cite{tanaka1+etal87,trujillo+etal08,guvendiren+etal10b}. In \cite{guvendiren+etal10b}, wrinkling and creasing-type instabilities have been studied in pHEMA gels%
\footnote{%
The acronym pHEMA stands for poly(2-hydroxyethyl methacrylate) \cite{guvendiren+etal10b}.
},
where it has been shown that the resulting surface morphology strongly depends on cross-linker concentration and solvent-polymer interaction. In particular, the authors demonstrated that while a specific hydrogel layer develops wrinkling in water, the same layer shows a transition from wrinkling to creasing when exposed to alcohol and alcohol-water mixture. In \cite{guvendiren+etal10c}, surface wrinkling was exploited to control stem-cell morphology and differentiation. Further studies on wrinkling and creasing of hydrogel bilayers have been performed analytically in \cite{kang+huang10a,kang+huang10b,wu+etal13} and computationally in \cite{dortdivanlioglu+linder19,ilseng+etal19,sriram+polukhov+keip21}. In contrast to analytical investigations, computational models usually allow for the examination of instabilities within a transient framework. Furthermore, they are not limited with respect to the complexity of underlying geometries.

In addition to surface instabilities, pattern-transforming instabilities in periodic hydrogel structures have been investigated experimentally in \cite{zhang+etal08,zhu2012capillarity,wu+etal14}. In \cite{zhang+etal08}, the authors have shown that when PDMS gels%
\footnote{%
The acronym PDMS stands for poly(dimethyl siloxane) \cite{zhang+etal08}.
} with periodically embedded circular pores of 1~$\mu\text{m}$ diameter are swollen in toluene, an organic solvent, elastic buckling of circular pores into elliptic shapes is observed. In the buckled state, the major axes of the ellipses are mutually orthogonal to one another developing a so-called diamond-plate pattern (see Fig.~\ref{fig.intro1}).
The authors have utilized this instability for the one-step assembly of patterned functional surfaces using superparamagnetic nanoparticles.
In \cite{zhu2012capillarity}, similar buckling patterns are exploited to manipulate the optical properties of pHEMA-based periodic hydrogel membranes. In this case, the authors observed diamond-plate patterns during drying of the membrane below a critical temperature after immersing it in deionized water.
Further experimental studies on periodic hydrogel films are documented in \cite{wang+etal17,ma+etal19}, where microstructures with high-swelling and non-swelling regions developed out-of-plane buckling patterns. It can be shown that the observed patterns depend significantly on the morphological parameters of the microstructure.

Similar to the mentioned experimental studies, \cite{okumura+etal14,okumura+etal15} deal with the 2D numerical investigation of the effect of prestrain and imperfections on the instability of periodically perforated hydrogels under plane-strain conditions. Their simulations are performed on a larger domain of a periodic structure contained from $2\times 2$ and $10\times 10$ unit-cell representative volume elements. In line with experiments, they have found that the diamond-plate patterns are a dominant instability mode, which are obtained considering random imperfections of the voids without eigenvalue analysis. Following the theoretical and numerical modeling approaches of these works supplemented further with eigenvalue buckling analysis in Abaqus, in \cite{miyoshi+etal21, kikuchi+etal22} 3D pattern transformations on a gel film bonded to a soft substrate are investigated. The authors have studied sequential bifurcations on bifurcated paths considering that the first bifurcation is associated with a hexagonal dimple mode, which led to complex patterns such as herringbone and labyrinth and gave further insight into the nature of the formation of these patterns. We note that in these contributions, the chemical potential is assumed spatially constant within the hydrogel domain, so that transient effects (such as time-dependent diffusion) was not considered. Furthermore, the instability analyses are usually considered on larger domains, which can be reduced to unit-cell computations based on Bloch-Floquet representation theorem as discussed in this work.

All mentioned experimental and numerical studies highlight the wide-ranging opportunities to tune the response of periodic systems by exploiting instabilities and hence to obtain materials with enhanced functionalities that could be considered, for example, in the development of soft and smart devices.

%
\subsection{Novel features and objectives of the present contribution}
%
Motivated by the above mentioned works, we study numerically structural as well as material instabilities in periodic hydrogel systems. All considered numerical setups are based on corresponding experimental investigations and the prescription of experimentally motivated boundary conditions. In particular, we study instabilities of perforated composite hydrogel microstructures in two dimensions (cf.\ \cite{zhang+etal08,zhu2012capillarity,wu+etal14})  as well as in composite hydrogel thin  films in three spatial dimensions (cf.\ \cite{wang+etal17,ma+etal19}). Here, we do not only observe short-wavelength structural instabilities, which alter the periodicity of the microstructure, but also long-wavelength structural instabilities considering a transient diffusion process.

To model the transient diffusion process in hydrogels, we employ a rate-type variational minimization formulation following \cite{miehe+etal15,boeger+etal17a}, see also \cite{polukhov+keip20} for an embedding into a multiscale framework. For spatial discretization in two dimensions, we exploit Raviart--Thomas-type finite elements coupled with an implicit Euler scheme for time discretization. In three spatial dimensions, we discretize the domain with (non-conforming) hexahedral finite elements. Within both space-time discrete formulations, we implement Bloch--Floquet wave analysis to determine pattern-transforming instabilities together with the altered periodicity of the microstructure on a unit-cell based computations without considering any imperfections apriori. Furthermore, we derive a condition for the loss of strong ellipticity based on effective moduli of perforated microstructures. The latter is based on an extension of the surface-averaging approach of \cite{miehe+koch02} to the coupled problem of transient chemo-mechanics.

Our studies reveal various experimentally confirmed as well as novel buckling patterns of periodic hydrogel systems depending on their morphology and associated material properties. 
For example, it will be shown that micro-coated hydrogel microstructures develop complex combinations of the well-known diamond-plate pattern coupled with intricate micro-wrinkling along internal surfaces. 
To our best knowledge, such pattern transformations have not been recorded yet and could serve for additional applications of periodic hydrogels.
Furthermore, our numerical studies in three spatial dimensions indicate new opportunities for the design of complex instability-induced pattern transformations driven by the morphology of the underlying microstructure.

\subsection{Outline}
%
The present work is structured as follows. As a theoretical basis for the computational analysis of hydrogels, we discuss a variational formulation and its numerical implementation in Section~\ref{sec20}. Based on that, we address fundamental aspects of structural and material stability analysis as well as the associated numerical implementation in Section~\ref{sec30}. The developed framework will then be applied to the computational investigation of instabilities of two- and three-dimensional hydrogel microstructures as well as corresponding parametric studies in Section~\ref{sec40}. Here, we also document a set of newly observed pattern transformations. We close with a discussion of the main contributions and findings in Section~\ref{sec50}.

\section{Variational formulation of deformation-diffusion processes}
\label{sec20}

In this section, we discuss the continuum modeling and computational implementation of the coupled chemo-mechanical problem. We adopt a variational minimization formulation as it has proven convenient for the investigation of instabilities \cite{sriram+polukhov+keip21}. In the underlying variational principle, the deformation map and the solvent-volume flux are chosen as the primary unknowns \cite{miehe+etal15,boeger+etal17a}. The main ingredients of this formulation are an energy-storage functional, a dissipation-potential functional and an external-power functional, where the latter arises due to applied traction and chemical-potential boundary conditions. For the latter, we propose additional boundary conditions in terms of the jumps of primary fields, which proves useful for the investigation of instabilities in periodic hydrogel structures.

\subsection{Primary state variables and the second law of thermodynamics}
\label{sec21}

We choose the deformation map $\Bvarphi$ and the solvent-volume flux $\KIH$ as the primary state variables of the considered minimization principle. If we denote by $\calB_0$ and $\calB_t$ the reference and the current configuration of the body, respectively, we can describe the mapping of these fields as
\begin{equation}
 \label{phi-Hfields}
  \Bvarphi:
  \begin{cases}
   \calB_0\times\calT \to \B_t\subset \calR^3\\[0.5ex]
   (\BX,t) \mapsto \Bvarphi(\BX,t)
  \end{cases} 
  \AND
  \KIH:
  \begin{cases}
    \calB_0\times\calT \to \R^3 \\[0.5ex]
    (\BX,t) \mapsto \KIH(\BX,t)
  \end{cases}\, ,
\end{equation}
where the deformation map characterizes the deformed position $\Bx=\Bvarphi(\BX,t)$ of a material point $\BX\in\calB_0$ at a given time $t\in\calT$. The solvent-volume flux is a relative Lagrangian field characterizing local transport of the solvent volume relative to the motion of the polymer%
\footnote{%
We refer to \cite{coussy+etal98} for a micromechanically motivated description of poromechanics.
}.

An additional state variable associated with the solvent-volume flux is given by the solvent-volume content $s:\calB_0\times\calT \to \calR^+$. For an arbitrary subdomain $\calP_0\subseteq\calB_0$, the solvent-volume content is determined from the balance of solvent volume \cite{gurtin+etal10}
\begin{equation}
  \label{sfield}
   \frac{\mathrm{d}}{\mathrm{d}t}\int_{\calP_0}s\,\mathrm{d}V = -\int_{\partial\calP_0}\KIH\cdot\BN\,\mathrm{d}A
   \quad
   \Rightarrow
   \quad
   \dot{s} = - \Div\KIH
   \quad
   \mbox{in}
   \
   \calP_0 .
\end{equation}  
Note that the flux field is defined in the Lagrangian configuration, so that no convective terms related to the spatial change of the solvent-volume content arise.

Based on the second law of thermodynamics, we can write the following inequality for the evolution of the stored-energy density of an arbitrary subdomain $\calP_0\subseteq\calB_0$
\begin{equation}
  \label{2ndlaw}
  \frac{\mathrm{d}}{\mathrm{d}t}\int_{\calP_0}\psi\,\mathrm{d}V\le 
  \underbrace{\int_{\partial\calP_0}\BT\cdot\dot{\Bvarphi}\,\mathrm{d}A-\int_{\partial\calP_0}\mu h\,\mathrm{d}A}_{P_{ext}(\dot{\Bvarphi},\KIH)}\, ,
\end{equation}
where we have introduced the traction vector $\BT$ acting on the surface of the deformed subdomain, defined per unit undeformed area $\partial\calP_0$, as well as the surface flux $h = \KIH\cdot\BN$; $\mu$ is the chemical potential of the solvent and characterizes the supplied energy to the system as a result of the influx of solvent volume. We refer to \cite{coussy+etal98} for an analogous interpretation in poromechanics.
Localizing the above inequality and distinguishing local and diffusive parts of the dissipation, we obtain
\begin{equation}
  \label{localD}
  \calD_{loc}:=\BP:\dot{\BF}+\mu\dot{s}-\dot{\psi}\ge 0 \quad \mbox{and} \quad \calD_{dif}:=\KIH\cdot\KIM \ge 0 \, ,
\end{equation}
where $\BF:=\Grad\varphi$ is the deformation gradient with $J:=\det\BF>0$ and $\KIM:=-\Grad\mu$ denotes the negative gradient of the chemical potential. Since we assume diffusion to take place in an otherwise elastic solid, the local part of the dissipation must vanish, i.e., $\calD_{loc}=0$. Thus, \eqref{localD}\textsubscript{1} can be exploited to obtain the constitutive relations for the first Piola-Kirchhoff stress tensor $\BP=\partial_\BF\hat\psi(\BF,s)$ and for the chemical potential $\mu=\partial_s\hat\psi(\BF,s)$. The inequality \eqref{localD}\textsubscript{2} characterizes the dissipation due to the diffusion process, where $\KIM$ is the driving force for the solvent-volume transport. The latter can be determined from a dissipation-potential density $\hat\phi(\KIH;\BF,s)$ for a given constitutive state, i.e., $\KIM = \partial_\KIH\hat\phi$.

\subsection{Constitutive functions: The free-energy and the dissipation-potential density}
\label{sec22}

As we have seen above, the stress tensor, the chemical potential and the negative gradient of the chemical potential are local constitutive fields. They are determined from a free-energy density $\hat\psi(\BF,s)$ and a dissipation-potential density $\hat\phi(\KIH;\BF,s)$.

In what follows, we consider an additive free-energy function of the form \cite{hong+etal08,hong+liu+suo09,boeger+etal17a,chester+anand10}
\begin{equation}
  \label{psiadditive}
  \hat\psi(\BF,s) = \hat\psi_{mech}(\BF) + \hat\psi_{chem}(s) + \hat\psi(J,s)\, .
\end{equation}
While the mechanical contribution of the free-energy function is assumed to be of Neo-Hookean-type, the chemical contribution is specified as a Flory--Rehner-type function \cite{flory+rehner43b}. These contributions together with a coupling term yield the overall free-energy function
\begin{equation}
  \label{psi-final0}
\hat\psi(\BF,s) ={\displaystyle\frac{\gamma}{2}\big[\BF:\BF-3-2\ln(\det\BF)\big]} +
 \alpha\big[ {\displaystyle   s\ln\big(\frac{s}{1+s}\big)+\frac{\chi s}{1+s}}\big]
+ {\displaystyle \frac{   \epsilon}{2}\big(\det\BF-1-s\big)^2} \, ,
\end{equation}
where $\gamma$ is the shear modulus and $\alpha:=kT/\nu$ is a mixing modulus with Boltzmann constant $k$, absolute temperature $T$ and volume of a solvent molecule $\nu$; $\chi$ is the Flory--Huggins interaction parameter, $\epsilon$ is a penalty parameter to constrain $\det\BF=1+s$ for incompressible hydrogels. We refer to Table~\ref{table01} for further information on the given parameters.

The dissipation-potential function is chosen such that \eqref{localD}\textsubscript{2} is satisfied a priori. Therefore, $\hat\phi(\KIH;\BF,s)$ is assumed to be a convex and normalized function of the solvent-volume flux $\KIH$ according to
$\hat\phi(\Bzero) = 0$ and $\partial_\KIH\hat\phi(\Bzero) = \Bzero$
at the reference configuration. In the present work, we will consider a dissipation function related to Fickian diffusion given by
\begin{equation}
  \label{phi-final0}
    {\displaystyle \hat\phi(\KIH;\BC,s)=\frac{1}{2Ms}\BC:(\KIH\otimes\KIH)}\, ,
\end{equation}
where $M$ is referred to as a mobility parameter. The dissipation-potential function is formulated in terms of given values $\{\BF,s\}$ and therefore evaluated numerically at the previous time step $(\cdot)_n$ within an incremental variational formulation \cite{miehe+etal15,boeger+etal17a}, see Section~\ref{sec24} for details.

The above introduced constitutive functions are given with respect to a dry reference configuration. Since \eqref{psi-final0} is singular for the dry state $s(t=0)=0$, a pre-swollen state is usually considered as the reference configuration \cite{hong+etal08,hong+liu+suo09,boeger+etal17a}.
%
According to that, the constitutive functions \eqref{psi-final0} and \eqref{phi-final0} defined per unit volume of the pre-swollen reference configuration and appear as
\begin{equation}
  \label{eq:pre-psi}
  \hat \psi
= \frac{\gamma}{2J_0} [ J_0^{2/3} \BF : \BF - 3 - 2 \ln(JJ_0) ]
+ \frac{\alpha}{J_0} [ s \ln ( \frac{s}{1+s} )
                     + \frac{\chi s}{1+s} ]
+ \frac{\epsilon}{2 J_0} ( J J_0 - 1 - s )^2 
\end{equation}
and
\begin{equation}
  \label{eq:pre-phi}
  \hat \phi
= \frac{1}{2 J_0^{1/3} M s_n} \BC_n : ( \KIH \otimes \KIH ) \, ,
\end{equation}
where the pre-swollen state is characterized by a Jacobian $J_0$. The initial condition for the solvent-volume content associated with the pre-swollen reference configuration reads
\begin{equation}
  \label{eq:pre-s0}
  s(t=0)=s_0 = \frac{\gamma}{\epsilon}\big(J_0^{-1/3}-\frac{1}{J_0} \big) + J_0 -1 \, . 
\end{equation}
The chemical potential at the pre-swollen reference configuration is determined as \cite{boeger+etal17a}
\begin{equation}
  \label{mu0}
  \mu_0 = -\frac{\epsilon}{J_0}(J_0-1-s_0) + \frac{\alpha}{J_0}\big[\ln{\big(\frac{s_0}{1+s_0} \big)}
  + \frac{1}{1+s_0} + \frac{\chi}{(1+s_0)^2} \big] \, . 
\end{equation}
The latter relation is particularly required when applying chemical-potential boundary conditions in the minimization formulation, see Section~\ref{sec40}.

\subsection{Boundary conditions for periodic hydrogel structures}
\label{sec23}

\begin{figure}%
\centering%
\footnotesize%
\psfrag{D0}  [l][l]{$\calD_0$}
\psfrag{B0}  [r][r]{$\calB_t$}
\psfrag{H}   [c][c]{$\calH_0$}
\psfrag{bB3d}   [c][c]{$\partial\calB_{side}$}
\psfrag{bBtop}   [r][r]{$\partial\calB_{top}$}
\psfrag{dS}  [l][l]{$\partial\calB_t$}
\psfrag{X}   [c][c] {$\BX$}
\psfrag{x}   [c][c] {$\Bx$}
\psfrag{F}[l][l] {$\Bvarphi(\BX,t)$}
\psfrag{phi}[c][c] {$\Bvarphi(\BX),\,\Ba(\BX)$}
\psfrag{N+}[c][c] {$\BN^+$}
\psfrag{N-}[c][c] {$\BN^-$}
\psfrag{n+}[c][c] {$\Bn^+$}
\psfrag{n-}[c][c] {$\Bn^-$}
\psfrag{Xm}[c][c] {$\BX^-$}
\psfrag{Xp}[c][c] {$\BX^+$}
\psfrag{xm}[c][c] {$\Bvarphi^-$}
\psfrag{xp}[c][c] {$\Bvarphi^+$}
\psfrag{Xj}[c][c] {$\jump{\BX}$}  
\psfrag{xj}[c][c] {$\jump{\Bvarphi}$}
\psfrag{boundary}[c][c]{$\jump{\Bvarphi}=\lbar{\BF}\jump{\BX},\
  \jump{\KIH}=\Bzero\quad \text{on}\quad \partial\calB\quad \text{and}\quad \KIH=\KIH^p\quad \text{on}\quad \partial\calH$}
\psfrag{xjp}[c][c] {$\jump{\Bvarphi_p}$}
\psfrag{xmp}[c][c] {$\Bvarphi_p^-$}
\psfrag{xpp}[c][c] {$\Bvarphi_p^+$}
\psfrag{aa}[c][c] {\textbf{(a)}}
\psfrag{bb}[c][c] {\textbf{(b)}}
\psfrag{microstructure}[c][c]{\shortstack{microstructure of \\periodic hydrogel}}
\psfrag{unitcell}[c][c]{unit-cell RVE}  
\includegraphics*[scale=0.6]{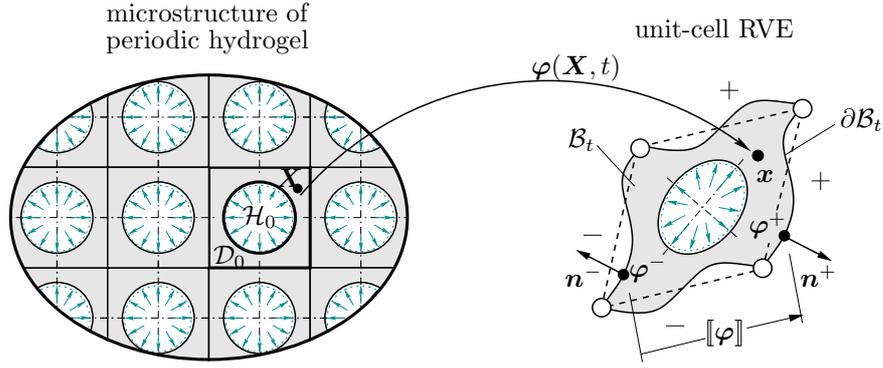}
\caption[Description of computational homogenization]{%
\textit{%
Boundary conditions for a two-dimensional periodic hydrogel microstructure}.
The response of periodic materials can be determined by means of unit-cell representative volume elements (unit-cell RVEs) of their microstructure. In the present approach, the unit-cell problem is solved in consideration of the jump conditions for the deformation map $\jump{\Bvarphi}=\lbar{\BF}\BX$ and the solvent-volume flux $\jump{\KIH}\cdot\BN=0$ across the external boundary according to \eqref{mechBC}. As can be seen in the figure, the solvent is assumed to enter the RVE through the boundary of the voids. In our computations, this is realized by applying the chemical potential at the void's boundary according to \eqref{chemBC}.}
\label{periodBC}%
\end{figure}%

Since our main aim is to investigate instabilities in periodic hydrogels, we will now provide a suitable setup for associated boundary conditions. As we will see, these are partly inspired from homogenization theory such that they allow for the definition of effective quantities, which could be used in the material stability analysis to be discussed later.

For illustration, let us consider the two-dimensional representation of a periodic hydrogel as depicted in Fig.~\ref{periodBC}. Here we assume that any out-of-plane deformation could be neglected and therefore typical plane-strain conditions apply. Since the structure is periodic, any stability analysis could be reduced to the smallest periodic domain, here given in terms of a square-shaped matrix with a circular void such that $\calB_0=\calD_0\cup\calH_0$ with $\calH_0$ denoting the domain of the void. In the following, we will refer to such a smallest possible periodic domain as a \textit{unit-cell representative volume element} (unit-cell RVE).

For the chemo-mechanical loading of the unit-cell RVE we prescribe boundary conditions of the generalized form
\begin{equation}
  \label{mechBC}
 \jump{\Bvarphi}\, =\lbar{\BF}\BX \quad \text{and} \quad \jump{\KIH}\cdot\BN=0 \quad \text{on} \quad \partial\calB_0 \, ,
\end{equation}
where we have introduced an \emph{effective} or \emph{macroscopic} deformation gradient $\lbar{\BF}$ \cite{miehe+schotte+schroeder99,polukhov+keip20}. The dual variable to the effective deformation gradient is the effective first Piola-Kirchhoff stress tensor $\lbar{\BP}$. Both quantities can be determined as the volume averages of their microscopic counterparts as%
\footnote{%
We note that when we prescribe the effective deformation $\lbar{\BF}$ via the RVE boundary, the effective stresses can be determined via \eqref{FbarPbar}\textsubscript{2}. Vice versa, it is also possible to apply effective stresses $\lbar{\BP}$ on the RVE. In that case, the effective deformation gradient can be determined via \eqref{FbarPbar}\textsubscript{1}, see \cite{miehe+koch02,miehe+etal16}.
}
\begin{equation}
  \label{FbarPbar}
  \lbar{\BF} = \frac{1}{|\calB_0|}\int_{\partial\calB_0}\Bvarphi\otimes\BN\,\mathrm{d}A
  \AND
  \lbar{\BP} = \frac{1}{|\calB_0|}\int_{\partial\calB_0}(\BP\BN)\otimes\BX\,\mathrm{d}A \, ,
\end{equation}
where $|\calB_0|$ denotes the volume of a unit-cell RVE. Note that according to \eqref{mechBC}\textsubscript{2}, we have defined vanishing jump conditions for the solvent-volume flux such that no effective flux is applied on the RVE. This setting is motivated by the assumption that the diffusion process is driven exclusively by suitable chemical-potential boundary conditions acting at the \textit{internal} boundary of unit-cell RVEs. Please refer to Section~\ref{sec40} for more details.

\subsection{Rate-type variational formulation of the deformation-diffusion problem}
\label{sec24}

Above, we have introduced the energy-storage function $\hat\psi(\BF,s)$, the dissipation-potential function $\hat\phi(\KIH)$ and the external-loading functional $P_{ext}(\dot{\Bvarphi},\KIH)$. Having these functions at hand, we define a rate-type potential to characterize the deformation-diffusion response \cite{miehe+etal15}
\begin{equation}
  \label{Pi}
  \begin{aligned}
  \Pi(\dot{\Bvarphi},\KIH)&:=\int_{\calB_0}\pi(\dot{\Bvarphi},\KIH)\,\mathrm{d}V - P_{ext}(\dot{\Bvarphi},\KIH)
  \WITH
  \pi(\dot{\Bvarphi},\KIH):=\frac{\mbox{d}}{\mbox{d}t}\hat\psi(\BF,s) + \hat\phi(\KIH;\BF,s) ,
  \end{aligned}
\end{equation}
where $\pi(\dot{\Bvarphi},\KIH)$ is the internal-power density. As a result, the variational principle for the considered problem in the continuous setting can be introduced as
\begin{equation}
  \label{var-form-continous}
  \{\dot{\Bvarphi}^\star,\KIH^\star\}=\Arg\big\{\inf_{\dot{\Bvarphi}\in\calW_{\dot{\Bvarphi}}}\inf_{\KIH\in\calW_{\KIH}}
  \Pi(\dot{\Bvarphi},\KIH) 
  \big\} \, ,
\end{equation}
where  $\calW_{\dot{\Bvarphi}}$ and $\calW_{\KIH}$ are the admissible spaces for the rate of the deformation map and the solvent-volume flux
\begin{equation}
  \label{Wcontinous}
  \begin{aligned}
  \calW_{\dot{\Bvarphi}}&:=\{\dot{\Bvarphi}\in H^1(\calB_0)\, \big|\,\dot{\Bvarphi}=\dot{\Bvarphi}_D \,\, \text{on} \,\, \partial\calB_0^\Bvarphi\} \AND \\
  \calW_{\KIH}&:=\{\KIH \in H(\Div,\calB_0)\, \big|\,\KIH=\KIH_D \,\, \text{on} \,\, \partial\calB_0^h\} \, . 
  \end{aligned}
\end{equation}
The Euler-Lagrange equations follow as the necessary conditions of the variational formulation as
\begin{equation}
  \label{Euler-equations}
  \boxed{
 \renewcommand*{\arraystretch}{1.25}
 \begin{array}{lrclcl}
  \mbox{Balance of linear momentum:}& \Div[\partial_{\BF}\hat\psi] &=& \Bzero & 
                                                        \text{in}  &\calB_0 \\ [1.5mm]
  \mbox{Inverse Fickian law:}& \nabla[\partial_s\hat\psi]+\partial_\KIH\hat\phi &=& \Bzero & \text{in}  &\calB_0 \\[1.5mm]
  \mbox{Mechanical traction:}& \partial_{\BF}\hat\psi\cdot\BN - {\BT} &=& \Bzero &
                                                         \text{on} &\partial\calB_0^T \\[1.5mm]
  \mbox{Chemical potential:}& -\partial_s\hat\psi + {\mu} &=& 0 & \text{on} &\partial\calB_0^{\mu} 
  \end{array} }
\end{equation}
%

\subsection{Space-time discretization of the variational formulation}
\label{sec241}

We solve \eqref{var-form-continous} using a conforming finite-element formulation. Before considering the discretization in space, we apply a modified version of the implicit Euler scheme in a time interval $[t_n,t]$ with $t$ being the current time step with $0<t\le T$ in a considered overall simulation time $T$. Note that for the variables evaluated at the current time step, we do not use subscripts. Consequently, the time-discrete version of \eqref{var-form-continous} reads%
\footnote{We remark that the variational formulation   \eqref{var-form-time-discrete} should also contain the energy-storage function as well as the external-power functional evaluated at the previous time step. However, since these terms are constant they vanish under variation and are therefore omitted from the equations to arrive at a compact representation.} \cite{miehe+etal15}
\begin{equation}
 \label{var-form-time-discrete}
\Pi^\tau(\Bvarphi,\KIH):=\int_{\calB_0}\hat\psi(\BF,s)+\tau\hat\phi(\KIH;\BF_n,s_n)\:\mathrm{d}V
-  P{}^\tau_{ext}({\Bvarphi},\KIH) \, ,
\end{equation}
where the solvent-volume content at the current time step is determined by $s=s_n-\tau\Div[\KIH]$. The dissipation-potential function $\hat\phi(\KIH;\BF_n,s_n)$ is evaluated using the deformation gradient and the solvent-volume content from the previous time step $t_n$ in order to stay consistent with the continuous variational setting \cite{miehe+etal15,boeger+etal17a}. For the explicit form of $\hat\phi$ we refer to \eqref{eq:pre-phi}. The time-discrete external-power functional reads
\begin{equation}
\label{Ptau}
P{}^\tau_{ext}({\Bvarphi},\KIH) = \int_{\partial\calB_0}\BT\cdot\Bvarphi\:\mathrm{d}A
-  \int_{\partial\calB_0}\tau\mu h\:\mathrm{d}A .
\end{equation}    
As a result, we obtain the time-discrete variational formulation with an initial condition for the solvent-volume content as
  \begin{equation}
    \label{Pitau}
  \{ {\Bvarphi^\star},\KIH^\star \}= \Arg \{\inf_{{\Bvarphi}\in\calW_{{\Bvarphi}}}
  \inf_{\KIH\in\calW_{\KIH}} \Pi^\tau ({\Bvarphi},\KIH) \} 
  \quad \text{with} \quad
  s(t_0) = s_0 . 
  \end{equation}
The above formulation is implemented into a conforming finite-element formulation by using a combination of Lagrange- and Raviart--Thomas-type shape functions \cite{boeger+etal17a,polukhov+keip20}. This then results in a minimization principle in the space-time discrete form as
 \begin{equation}
    \label{Pitauh}
  {\Bsfd^\star}= \Arg \{\min_{\Bsfd}
   \Pi^{\tau h} (\Bsfd) \} ,
 \end{equation}
where $\Bsfd$ contains all degrees of freedom comprising the node-based displacements and the edge-based normal-projected solvent-volume fluxes.

\section{Stability analysis of periodic hydrogel systems}
\label{sec30}

In this section, we describe methods to investigate the structural and the material stability of periodic hydrogel systems. For the structural stability analysis, we utilize the minimization structure of the underlying incremental system of equations
\cite{hill57,hill62,miehe+schotte+schroeder99}. In doing so, we investigate in particular \emph{short-wavelength} and \emph{long-wavelength} instabilities. The former are usually associated with a change in periodicity of an underlying microstructure, which could arise, for example, due to wrinkling of internal surfaces in perforated systems. We will detect these types of instabilities via Bloch--Floquet analysis. Long-wavelength instabilities are related to the loss of strong ellipticity of a corresponding homogenized continuum. These kinds of instabilities can be detected by checking the positive definiteness of an underlying acoustic tensor. In what follows, we will also present a corresponding numerical method for the analysis of related instabilities.

\subsection{Minimization-based structural stability analysis}
\label{sec31}

Within a minimization-based approach, hydrogel systems are considered to be structurally stable at an equilibrium state $\Bsfd^\star$, if all variations of the incremental potential from this equilibrium state under infinitesimally small perturbations are positive, that is \cite{miehe+schotte+schroeder99}
\begin{equation}
  \label{DeltaPi}
\Delta\Pi{}^{\tau h}(\Bsfd^\star,\delta\Bsfd):=
\Pi{}^{\tau h}(\Bsfd^\star+\epsilon\delta\Bsfd)-\Pi{}^{\tau h}(\Bsfd^\star)
>0 \, ,
\end{equation}
where $\epsilon$ is a small parameter. Considering a Taylor approximation of \eqref{DeltaPi} and truncating higher-order terms, we obtain
\begin{equation}
\Delta\Pi{}^{\tau h}(\Bsfd^\star,\delta\Bsfd) \approx
\frac{\mathrm{d}}{\mathrm{d}\epsilon}\bigg|_{\epsilon=0}\Pi{}^{\tau h}(\Bsfd^\star+\epsilon\delta\Bsfd) 
+\frac{1}{2}\frac{\mathrm{d}^2}{\mathrm{d}\epsilon^2}\bigg|_{\epsilon=0}\Pi{}^{\tau h}(\Bsfd^\star+\epsilon\delta\Bsfd)
>0 \, .
\end{equation}
Since the first variation vanishes at equilibrium, we arrive at
\begin{equation}
\Delta\Pi{}^{\tau h}(\Bsfd^\star,\delta\Bsfd)
\approx
\frac{1}{2}\frac{\mathrm{d}^2}{\mathrm{d}\epsilon^2}\bigg|_{\epsilon=0}\Pi{}^{\tau h}(\Bsfd^\star+\epsilon\delta\Bsfd) 
>0 \, .
\end{equation}
Consequently, structural stability is associated with the second variation of the incremental potential in the sense that
\begin{equation}
  \label{Pi2ndvar}
\Delta\Pi{}^{\tau h}(\Bsfd^\star,\delta\Bsfd)=
\frac{1}{2}\delta\Bsfd{}^T\big[\Pi{}^{\tau}(\Bsfd){}_{,\Bsfd\Bsfd}\big]_{\Bsfd^\star}\delta\Bsfd
=\Lambda\delta\Bsfd^T\delta\Bsfd
>0 .
\end{equation}
The latter results in an eigenvalue problem in terms of the tangent of the underlying system of equations according to
\begin{equation}
  {\Lambda} =
  \min_{\Bn\in\IN^3}\min_{\Bsfd}\big\{\Delta\Pi{}^{\tau h}(\Bsfd,\delta\Bsfd)
  \big/\|\delta\Bsfd\|^2\big\}
\begin{cases}
 {>   0} \quad \text{for stable state} \ \Bsfd^\star \\
 {\le 0} \quad \text{for unstable state} \ \Bsfd^\star 
\end{cases} \, ,
\end{equation}
where $\Bn\in\IN^3$ denotes the number of unit-cells RVEs with respect to three spatial directions \cite{mueller87}. The minimization over $\Bn$ allows for the detection of instabilities that would not be describable by means of a unit-cell RVE. However, the above approach is computationally inefficient since it requires the consideration of various realizations of RVEs that grow in size with growing $\Bn$. Therefore, we use Bloch--Floquet analysis to detect such critical structural instabilities on a unit-cell RVE \cite{geymonat+mueller+triantafyllidis93}.

Based on Bloch--Floquet analysis, the variations $\delta\Bsfd$ of the degrees of freedom can be represented in terms of unit-cell variations $\delta\Bsfd_{\calB_0}$ as
\begin{equation}
  \label{deldBFA}
  \delta\Bsfd(\BX)=\delta\Bsfd_{\calB_0}(\BX)\exp[\mathrm{i}\Bk\cdot\BX] 
  \quad \text{with} \quad 
  \delta\Bsfd^+_{\calB_0} = \delta\Bsfd^-_{\calB_0} 
  \quad \text{on} \quad
  \partial\calB_0=\partial\calB_0^+\cup\partial\calB_0^- ,
\end{equation}
where $\Bk$ is referred to as Bloch vector, which characterizes the wavelength of the variations. Consequently, this expression leads to the boundary conditions
\begin{equation}
  \label{BFA-BC}
  \delta\Bsfd(\BX^+) = \delta\Bsfd(\BX^-)\exp[\mathrm{i}\Bk\cdot(\BX^+-\BX^-)]
\end{equation}
for a unit-cell RVE. Here, $\BX^+$ and $\BX^-$ denote the reference placement of periodic nodes, which are located at opposite boundaries of the unit-cell RVE, see Fig.~\ref{periodBC}. It follows a stability criterion related to the tangent of a unit-cell RVE according to%
\footnote{%
For details related to the numerical implementation we refer to \cite{triantafyllidis+etal06,bertoldi+etal08,polukhov+vallicotti+keip18,polukhov+keip21}.
}
\begin{equation}
  \label{Lambdak}
  \boxed{
  {\Lambda} =
  \min_{k_i\in[0,\pi]}\min_{\Bsfd}\big\{\Delta\Pi{}^{\tau h}(\Bsfd,\delta\Bsfd)
  \big/\|\delta\Bsfd\|^2\big\}
\begin{cases}
 {>   0} \quad \text{for stable state} \ \Bsfd^\star \\
 {\le 0} \quad \text{for unstable state} \ \Bsfd^\star 
\end{cases}} \, .
\end{equation}
Depending on the values of the Bloch vector components $k_i\in\calR^3$, we distinguish three types of instabilities: (i) unit-cell periodic instabilities for $k_i=0$; (ii) short-wavelength instabilities for finite values of $k_i$ usually spanning RVEs made of multiple unit cells, where the explicit number of cells can be computed as $n_i=2\pi/k_i$; (iii) long-wavelength instabilities for $\Bk\to\Bzero$. We note that the latter could also be determined from the effective mechanical moduli of periodic hydrogels, see below.

%
\subsection{Loss of strong ellipticity: Material stability analysis at macroscale}
\label{sec32}

Having defined macroscopic mechanical boundary conditions in terms of $\lbar{\BF}$, we assume that at macroscale we need to fulfill the balance equation \cite{polukhov+keip20}
\begin{equation}
  \label{DivPbar}
  \llbar{\Div}[\lbar{\BP}{}(\lbar{\BF},\llbar{\Div}[\lbar{\KIH}])]=\lbar{\rho}_0\ddot{\lbar{\Bvarphi}} \, .
\end{equation}
In order to check if the latter equation satisfies the strong-ellipticity condition, we superimpose a plane wave of the form $\delta\lbar{\Bvarphi}=\lbar{\Bn}f(\lbar{\BX}\cdot\lbar{\BN}-\lbar{c}t)$, such that for infinitesimal $\delta\lbar{\Bvarphi}$ it holds that \cite{hill62,miehe+schotte+schroeder99}
\begin{equation}
  \label{DivdelPbar}
  \llbar{\Div}[\delta\lbar{\BP}{}]=\llbar{\Div}[\lbar{\IA}{}:\delta\lbar{\BF}]=\lbar{\rho}_0\delta\ddot{\lbar{\Bvarphi}}
  \quad \text{with} \quad \delta\lbar{\BF}=\llbar{\Grad}\delta\lbar{\Bvarphi} \, ,
\end{equation}
where $\lbar{\IA}$ are the effective mechanical moduli of the microstructure. Note that since we do not consider any macroscopic loading in terms of $\llbar{\Div}[\lbar{\KIH}]$, we do not have to take into account a corresponding variation in the above equation.
%
%
After some modifications, \eqref{DivdelPbar} appears as
\begin{equation}
  \label{Qbar}
  \lbar{\Bn}\cdot\lbar{\BQ}\cdot\lbar{\Bn} = \lbar{\rho}_0\lbar{c}{}^2
  \, ,
\end{equation}
where $\lbar{Q}_{ab}=(\lbar{\IA}{})\indices{_a^A_b^B}\lbar{N}_A\lbar{N}_B$ is the mechanical acoustic tensor defined in terms of the effective mechanical moduli. It follows the criterion for the material stability at the macroscopic scale
\begin{equation}
  \label{Lambdabar}
  \boxed{
  \lbar{\Lambda}= \min_{\|\lbar{\BN}\|=1}\big\{ \lbar{\Bn}\cdot\lbar{\BQ}\cdot\lbar{\Bn} \big\}
  \begin{cases}
   > 0 \quad \text{for materially stable state} \\
  \le 0 \quad \text{for materially unstable state} \\    
  \end{cases}  } \, .
\end{equation}
%

%
\subsection{Computation of the effective mechanical moduli of voided hydrogel structures}
\label{sec311}
We now describe the computation of the effective mechanical moduli of voided hydrogel structures following the approach proposed by \citet{miehe+koch02}. This approach allows to determine effective moduli based on surface integrals along the external boundary of RVEs. Starting point is the definition of the  effective mechanical stresses \eqref{FbarPbar}\textsubscript{2} and the effective chemical potential \cite{polukhov+keip20}
\begin{equation}
  \label{Pbarmubarsurface}
  \lbar{\BP} =  \frac{1}{|\calB_0|}\int_{\partial\calB_0}(\BP\BN)\otimes\BX\,\mathrm{d}A
  \AND \lbar{\mu} = \frac{1}{3|\calB_0|}\int_{\partial\calB_0}(\mu\BN)\cdot\BX\,\mathrm{d}A\, ,
\end{equation}
where $\partial \calB_0$ denotes the external boundary of an RVE. In the discrete setting, the above equations can be given as follows (see \cite{miehe+koch02} for the purely mechanical case)
\begin{equation}
    \label{Pbarlambda1}
    \lbar{\BP} =  \frac{1}{|\calB_0|}\sum_{I=1}^P\Blambda^{f}_I\otimes\jump{\BX^b_I}
    \AND
    \lbar{\mu} = \frac{1}{3|\calB_0|}\sum_{I=1}^P\lambda^{\mu}_I\BN\cdot\jump{\BX^b_I}\, ,
\end{equation}
where $\Blambda^f_I$, $\lambda^\mu_I$ and $\BX^b_I$ denote the nodal forces, the chemical potential and the coordinates associated with the node $I$ on the external boundary of an RVE, respectively. The number $P=M/2$ is defined from the total number of boundary nodes $M$. Introducing an array of nodal forces and chemical potentials $\Blambda:=[\Blambda^f_1,\lambda^\mu_1,...,\Blambda^f_P,\lambda^\mu_P]^T$ as well as an array containing the effective stresses (in vector notation) and the effective chemical potential  $\lbar{\SFB{P}}:=[\lbar{\BP},\lbar{\mu}]^T$, we can write \eqref{Pbarlambda1} in a compact form%
\footnote{$\lbar{\SFB{P}}:=[\lbar{P}_{11},\lbar{P}_{22},\lbar{P}_{12},\lbar{P}_{21},\lbar{\mu}]^T$ is considered for the corresponding array in two dimensions.}%
\begin{equation}
    \label{Pbarlambda2}  
    \lbar{\SFB{P}} =  \frac{1}{|\calB_0|}\sum_{I=1}^P\IQ_I\Blambda_I=\frac{1}{|\calB_0|}\IQ\Blambda\, ,
\end{equation}
where $\IQ:=[\IQ_1,...,\IQ_P]$ is the array of projection tensors with
\begin{equation}
  \IQ_I^T:=
  \begin{bmatrix}
    \jump{(\BX^b_I)_1} & 0 &     \jump{(\BX^b_I)_2} & 0 & 0 \\[2mm]
    0 &     \jump{(\BX^b_I)_2} & 0 &     \jump{(\BX^b_I)_1} & 0   \\[2mm]
    0 & 0 & 0& 0 & \displaystyle \frac{1}{3}\jump{\BX_I^b}\cdot\BN
  \end{bmatrix}
\end{equation}
in two dimensions.

In the following, we are only interested in the effective mechanical moduli. Since these relate the incremental effective stresses to the incremental effective deformation gradient according to $\Delta\lbar{\BP}=\lbar{\IA}:\Delta\lbar{\BF}$, we first determine the increment of the first Piola-Kirchhoff stress tensor in vector notation from \eqref{Pbarlambda2} as 
\begin{equation}
\Delta\lbar{\BP} = \IL\Delta\lbar{\SFB{P}}=\frac{1}{|\calB_0|}\IL\IQ\Delta\Blambda\, ,
\end{equation}
where $\IL$ is a projection tensor with entries $0$ or $1$.

 In order to obtain the final expression for the effective-moduli tensor, we need to have a relation between the increment of the generalized nodal forces $\Delta\Blambda$ and the effective deformation gradient $\lbar{\BF}$. We obtain such a relationship from the following variational principle in the space-time discrete setting 
\begin{equation}
  \label{Pitilde}
  \ltilde{\Pi}{}^{\tau h}(\Bsfd,\Blambda) = \Pi{}^{\tau h}(\Bsfd) + \Blambda\cdot(\jump{\Bsfd^b}-\IQ^T\lbar{\SFB{F}})\, ,
%
\end{equation}
where $\lbar{\SFB{F}}:=[\lbar{\BF},\llbar{\Div}\lbar{\KIH}]^T$ is an array given by $\lbar{\SFB{F}}=[\lbar{F}_{11},\lbar{F}_{22},\lbar{F}_{12},\lbar{F}_{21},\llbar{\Div}\lbar{\KIH}]^T$ in two dimensions.

The second term on the right-hand side of \eqref{Pitilde} considers the Lagrange multiplier method to enforce  the jump boundary conditions for the deformation map and the fluid-volume flux vector of the form
\begin{equation}
\jump{\Bvarphi}=\lbar{\BF}\jump{\BX} \AND \jump{\KIH}\cdot\BN=\frac{1}{3}\llbar{\Div}\lbar{\KIH}\jump{\BX}\cdot\BN \, ,
\end{equation}
refer to \eqref{mechBC} as well as to \cite{polukhov+keip20} for more details. Furthermore, $\BX^b:=[\BX^b_1,...,\BX^b_P]^T$ in \eqref{Pitilde} contains the coordinates of the nodes at the external boundary of an RVE. The Lagrange multiplier $\Blambda$  can be identified as the  array of nodal forces and chemical potentials. Its interpretation follows from the Euler-Lagrange equations of the variational principle. Furthermore, $\Bsfd^b:=[\Bsfd^b_1,..,\Bsfd^b_P]^T$ denotes the degrees of freedom at the boundary of the unit-cell RVE. In \eqref{Pitilde}, we can express the jump of the degrees of freedom in terms of the latter as
\begin{equation}
  \jump{\Bsfd^b}:=\assemble_{I=1}^{P}
  \begin{bmatrix}
    (\Bsfd^b_\Bvarphi)_I^+ -     (\Bsfd^b_\Bvarphi)_I^- \\[3mm]
    (\mathsf{d}^b_h)_I^+  +     (\mathsf{d}^b_h)_I^- 
  \end{bmatrix} = \IP\Bsfd^b \, ,
\end{equation}
where $\IP$ is a projection tensor with entries $\{-1,0,1\}$. Note that, since we are considering Raviart--Thomas-type elements, the edge-based solvent-volume flux degrees of freedom are anti-periodic at the external boundary \cite{polukhov+keip20}.

The microscopic problem is usually solved under constant macroscopic deformation gradient for each time step. Therefore, the first variation of \eqref{Pitilde} gives 
\begin{equation}
  \displaystyle
  \delta\ltilde{\Pi}{}^{\tau h} =
  \begin{bmatrix}
  \displaystyle    
    \delta\Bsfd^i \\
    \delta\Bsfd^b \\
    \delta\Blambda    
  \end{bmatrix}
  \cdot
  \begin{bmatrix}
  \displaystyle    
    \partial_{\Bsfd^i}\Pi^{\tau h} \\
    \partial_{\Bsfd^b}\Pi^{\tau h} + \IP\Blambda\\
    \IP\Bsfd^b - \IQ^T\lbar{\SFB{F}}
  \end{bmatrix}  =
  \begin{bmatrix}
  \displaystyle    
    \delta\Bsfd^i \\
    \delta\Bsfd^b \\
    \delta\Blambda    
  \end{bmatrix}
  \cdot  
  \begin{bmatrix}
  \displaystyle    
    \BsfR^i \\
    \BsfR^b + \IP\Blambda\\
    \IP\Bsfd^b - \IQ^T\lbar{\SFB{F}}
  \end{bmatrix}  = \Bzero  \, ,
\end{equation}  
where $\Bsfd^i$ denote the interior degrees of freedom of a unit-cell RVE. Once the above system of equations is solved at the microscale, the sensitivity of the system with respect to the macroscopic deformation gradient can be determined from
\begin{equation}
  \label{Delta2Pitilde}
    \displaystyle
    \Delta\delta\ltilde{\Pi}{}^{\tau h} =
  \begin{bmatrix}
  \displaystyle    
    \delta\Bsfd^i \\
    \delta\Bsfd^b \\
    \delta\Blambda    
  \end{bmatrix}
  \cdot      
  \begin{bmatrix}
  \displaystyle    
    \BsfK^{ii} & \BsfK^{ib} & \Bzero & \Bzero\\
    \BsfK^{bi} & \BsfK^{bb} & \IP^T & \Bzero \\
    \Bzero &     \IP & \Bzero &  - \IQ^T
  \end{bmatrix}
  \cdot
  \begin{bmatrix}
  \displaystyle    
    \Delta\Bsfd^i \\
    \Delta\Bsfd^b \\
    \Delta\Blambda \\
    \Delta\lbar{\SFB{F}}
  \end{bmatrix}
  = \Bzero  \, .
\end{equation}  
First, from the variations with respect to $\delta\Bsfd^i$ and $\delta\Bsfd^b$ in
\eqref{Delta2Pitilde}\textsubscript{1} and \eqref{Delta2Pitilde}\textsubscript{2}, we obtain elimination equations for $\Delta\Bsfd^i$ and $\Delta\Bsfd^b$ as
\begin{equation}
\Delta\Bsfd^i=-(\BsfK^{ii})^{-1}\BsfK^{ib}\Delta\Bsfd^b \AND \Delta\Bsfd^b=-(\underline{\boldsymbol{\mathsf{\ltilde{K}}}} {}^{bb})^{-1}\IP^T\Delta\Blambda\, ,
\end{equation}
where $\underline{\boldsymbol{\mathsf{\ltilde{K}}}} {}^{bb} = \BsfK{}^{bb} - \BsfK{}^{bi}(\BsfK{}^{ii})^{-1}\BsfK{}^{ib}$. Next, we consider the variation with respect to $\delta\Blambda$ in \eqref{Delta2Pitilde} and obtain the searched expression of $\Delta\Blambda$ given by
\begin{equation}
\Delta\Blambda = \big[\IP(\underline{\boldsymbol{\mathsf{\ltilde{K}}}} {}^{bb})^{-1}\IP^T\big]{}^{-1}\IQ^T\Delta\lbar{\SFB{F}} \, .
\end{equation}  
Finally, we can compute the effective mechanical moduli as
\begin{equation}
  \label{Pbarlambda}
  \boxed{
    \Delta\lbar{\BP} =  \frac{1}{|\calB_0|}\IL\IQ\big[\IP(\underline{\boldsymbol{\mathsf{\ltilde{K}}}} {}^{bb})^{-1}\IP^T\big]{}^{-1}\IQ^T\IL^T\Delta\lbar{\BF} \quad \Rightarrow \quad
    \lbar{\IA}: = \frac{1}{|\calB_0|}\IL\IQ\big[\IP(\underline{\boldsymbol{\mathsf{\ltilde{K}}}} {}^{bb})^{-1}\IP^T\big]{}^{-1}\IQ^T\IL^T  } \, .
\end{equation}

\section{Numerical investigation of instabilities in hydrogel structures}
\label{sec40}

We now investigate computationally microscopic and macroscopic instabilities in two- and three-dimensional periodic hydrogels, where all  considered boundary value problems are motivated from experimental studies (e.g., \cite{zhang+etal08,zhu2012capillarity,wu+etal14}). The two-dimensional studies of single- and two-phase voided hydrogel microstructures under plane-strain conditions are summarized in Section~\ref{sec41}. Corresponding investigations of structrual instabilities of three-dimensional hydrogel thin films are documented in Section~\ref{sec45}. In all cases, we examine the influence of physical properties and microscopic morphology on the onset and the type of instabilities. 


\subsection{Two-dimensional hydrogel structures}
\label{sec41}

In this section, we study instabilities in single-phase and two-phase voided hydrogel microstructures. While the single-phase hydrogels consist of a hydrogel matrix and voids, the two-phase hydrogels further contain stiff, permeable coatings surrounding the voids as described in Fig.~\ref{2D_BVPs}. 
The latter structures are motivated from experimental studies on constrained plane bilayer films, which reveal a rich set of wrinkling patterns \cite{guvendiren+etal09,guvendiren+etal10a,guvendiren+etal10b}; see also the analytical and numerical studies \cite{kang+huang10a,kang+huang10b,cao+hutchinson12,wu+etal13,dortdivanlioglu+linder19,sriram+polukhov+keip21}. Associated wrinkling phenomena have not only been reported for plane bilayers, but also for bilayer tubes \cite{jin+etal18,sriram+polukhov+keip21} and could be exploited in various engineering and biomedical applications such as controlled formation of micro-gears \cite{yin+etal09}, generation of multi-cellular spheroids \cite{zhao+etal14}, controlled cell spreading \cite{jiang+etal02}, see also \cite{lam+etal06,guvendiren+etal10c}. To the best of our knowledge, \textit{swelling-induced micro-wrinkling at internal surfaces of periodic hydrogels} has not been reported in the literature yet. The sole experimental study coming close to the results reported here is the one recently documented in \cite{fan20223d}.

In the following, we first discuss suitable boundary conditions and numerical setups of boundary value problems to analyze the given hydrogel microstructures. Next, we provide some numerical studies in order to give further insight into the response of the voided hydrogel microstructures. Eventually, we perform a set of computational investigations demonstrating instabilities in two-dimensional hydrogel microstructures.

\subsubsection{Boundary conditions in the two-dimensional setting}
\label{sec411}

%
\begin{figure}%
\centering%
\footnotesize%
\psfrag{D0}  [l][l]{$\calD_0$}
\psfrag{A}  [c][c]{\color{red} $A$}
\psfrag{B0}  [r][r]{$\calB_t$}
\psfrag{H}   [c][c]{$\calH_0$}
  \psfrag{dS}  [l][l]{$\partial\calB_t$}
  \psfrag{X}   [c][c] {$\BX$}
  \psfrag{x}   [c][c] {$\Bx$}
  \psfrag{F}[l][l] {$\Bvarphi(\BX,t)$}
  \psfrag{phi}[l][l] {$\Bvarphi(\BX),\,\Ba(\BX)$}
  \psfrag{N+}[c][c] {$\BN^+$}
  \psfrag{N-}[c][c] {$\BN^-$}
  \psfrag{n+}[c][c] {$\Bn^+$}
  \psfrag{n-}[c][c] {$\Bn^-$}
  \psfrag{Xm}[c][c] {$\BX^-$}
  \psfrag{Xp}[c][c] {$\BX^+$}
  \psfrag{xm}[c][c] {$\Bvarphi^-$}
  \psfrag{xp}[c][c] {$\Bvarphi^+$}
  \psfrag{Xj}[c][c] {$\jump{\BX}$}  
  \psfrag{xj}[c][c] {$\jump{\Bvarphi}$}
  \psfrag{boundary}[c][c]{$\jump{\Bvarphi}=\lbar{\BF}\jump{\BX},\
    \jump{\KIH}=\Bzero\quad \text{on}\quad \partial\calB\quad \text{and}\quad \KIH=\KIH^p\quad \text{on}\quad \partial\calH$}
  \psfrag{xjp}[c][c] {$\jump{\Bvarphi_p}$}
  \psfrag{xmp}[c][c] {$\Bvarphi_p^-$}
  \psfrag{xpp}[c][c] {$\Bvarphi_p^+$}
  \psfrag{aa}[c][c] {\textbf{(a)}}
  \psfrag{bb}[c][c] {\textbf{(b)}}
  \psfrag{micro}[c][c]{\shortstack{microstructure of \\periodic hydrogels}}
  \psfrag{unit}[c][c]{unit-cell RVE}  
\includegraphics*[scale=0.6]{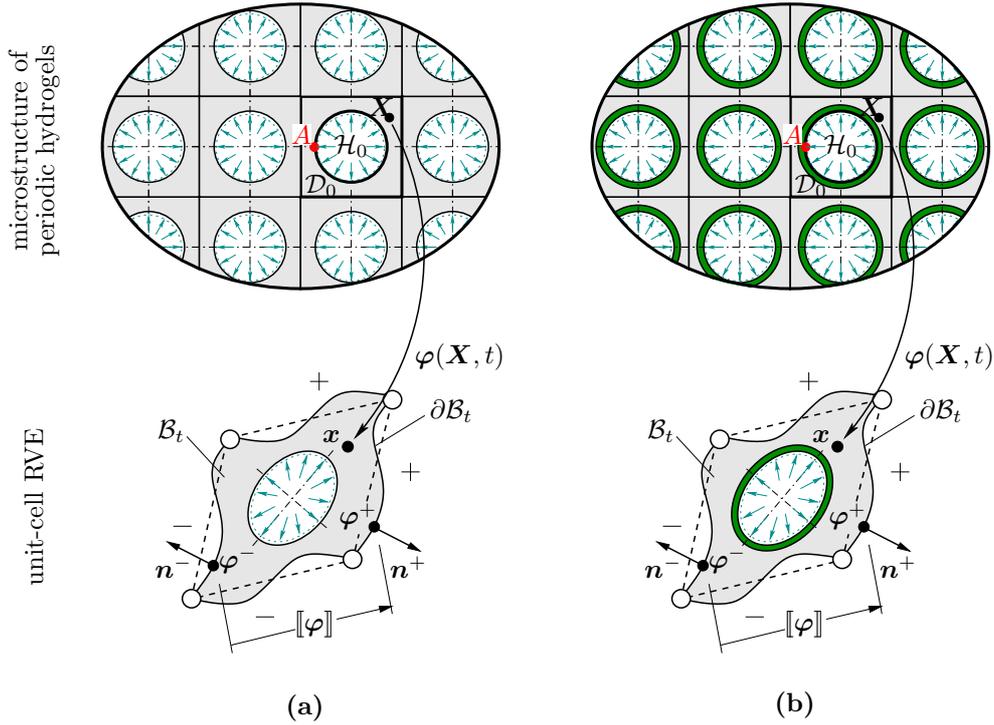}
\caption[Description of computational homogenization]{%
\textit{%
Unit-cell RVEs of two periodic hydrogel microstructures%
}.
We consider \textbf{(a)} single-phase hydrogels 
(consisting of hydrogel matrix and voids) 
and \textbf{(b)} two-phase hydrogels (consisting of hydrogel matrix, hydrogel coatings, and voids). 
The coating is realized as a stiff film that is perfectly bonded to the matrix (highlighted in green). It is further assumed that the solvent enters the hydrogel through the internal boundary $\partial\calH_0$, where we prescribe suitable boundary conditions in terms of the chemical potential. At the external boundary $\partial\calB_0$, the jump conditions for the deformation as well as for the solvent-volume flux are applied as described in \ref{sec23}.
}%
\label{2D_BVPs}%
\end{figure}%
Since we consider periodic hydrogel microstructures, all computations are reduced to a periodic unit-cell RVE as discussed in Section~\ref{sec23}. Thus, the jump conditions for the deformation map and the solvent-volume flux according to \eqref{mechBC} are considered at the external boundary of the RVE.
%
%
In what follows, we set $\lbar{\BF} = \Bone$ to mimic an effectively constrained hydrogel composite. We further assume that a solvent diffuses into the hydrogel solely through the boundary of the voids $\partial\calH_0$ and thus consider the mechanical traction and chemical potential on $\partial\calH_0$ as
\begin{equation}
  \label{chemBC}
  \BP\cdot\BN=\Bzero \AND
 \mu = \lambda(t)\mu_0 \quad \text{on} \quad \partial\calH_0  
  \quad \text{with} \quad
  \lambda(t) =
  \begin{dcases}
    1-t \, &\text{for} \quad    t\le 1\, \text{sec}\\
    0   \, &\text{for} \quad    t >  1\, \text{sec}  
  \end{dcases}  \, ,
\end{equation}
where $\mu_0$ is the chemical potential of the reference state, see \eqref{mu0}. As can be seen in the latter equation, the chemical potential is linearly increased from $\mu_0$ to zero within one second \cite{boeger+etal17a,dortdivanlioglu+linder19}.

In the present work, we consider a time step of $\tau = 4 \cdot 10^{-3} \, \text{sec}$ in order to determine instability points accurately. We note that a suitable time-step size is related to the material parameters of the boundary value problem as well as its numerical discretization \cite{teichtmeister+etal19}.

In Fig.~\ref{2D_BVPs}, we depict illustrations of microstructures considered in the first part of our studies. The microstructures are discretized using Raviart--Thomas-type finite elements of the type Q$_2$RT$_0$, where the displacement field is interpolated using bi-quadratic Lagrangian shape functions associated with node-based degrees of freedom and the solvent-volume flux field is interpolated using linear Raviart--Thomas shape functions associated with edge-based degrees of freedom \cite{boeger+etal17a,polukhov+keip20}. Depending on the geometry of the problem, approximately $5,000$ to $9,000$ finite elements (corresponding to approximately $20,000$ to $30,000$ nodes) are used for the discretization of a unit-cell RVE. The used material parameters are summarized in Table~\ref{table01}.

%
\subsubsection{Response of voided hydrogel microstructures}%
%
\label{sec42}

\begin{table}[t]
\footnotesize
\centering
\renewcommand{\arraystretch}{1.3}
\caption{Material parameters of the considered hydrogel microstructures.}

\begin{tabular}{llllll}
\hline
Parameter & Name, unit & Value & Eq. \\
\hline
$\gamma$   & Shear modulus,  N/mm\textsuperscript{2}        & $0.01$     & \eqref{eq:pre-psi} \\
$\alpha$   & Mixing modulus, N/mm\textsuperscript{2}        & $20.0$     & \eqref{eq:pre-psi}    \\
$\chi$     & Mixing parameter, --                           & $0.1$      & \eqref{eq:pre-psi}   \\
$M$        & Mobility parameter, mm\textsuperscript{4}/(Ns) & $10^{-4}$  & \eqref{eq:pre-phi}\\
$\epsilon$ & Penalty parameter, N/mm\textsuperscript{2}     & $10\gamma$ &  \eqref{eq:pre-psi}  \\
$J_0$      & Pre-swollen Jacobian, --                       & $1.01$     &  \eqref{eq:pre-psi} $\&$ \eqref{eq:pre-phi}  \\
\hline
\end{tabular}
\label{table01}
\end{table}

\begin{figure}[h]
  \hspace{-5mm}
  \scalebox{0.73}{\input{./time-r}}  
  \hspace{-10mm}
  \scalebox{0.73}{\input{./time-delR}}
  \vspace{-9mm}
  \hspace{-5mm}
  
  { \text{(a)} \hspace{-5mm} \hspace{0.5\textwidth} \text{(b)}}
  
  \vspace{8mm}
  \hspace{-5mm}
  \scalebox{0.73}{\input{./time-h_v}}
  \hspace{-10mm}
  \scalebox{0.73}{\input{./time-pbar11}}
  \vspace{-9mm}
  \hspace{-5mm}
  
  { \text{(c)} \hspace{-5mm} \hspace{0.5\textwidth} \text{(d)}}  
    
  \caption{\emph{Response of voided hydrogel microstructures versus simulation time.}
    \textbf{a)} Radius $r_A$ of voids and
    \textbf{b)} change in the radius $\Delta r_A$ of voids measured at the point $A$ as depicted in Fig.~\ref{2D_BVPs} with the initial radius $r_{initial} = \sqrt{f_0/\pi}$;
    \textbf{c)} Solvent volume entering the hydrogel matrix through the void boundary $h_{void}=-\int_{\partial\calH_0}\KIH\cdot\BN\, \mathrm{d}A$;
    \textbf{d)} Effective first Piola-Kirchhoff stress $\lbar{P}_{11}$ normalized with respect to the shear modulus of the matrix $\gamma$.}
  \label{fig:res-vs-time_voided}  
\end{figure}
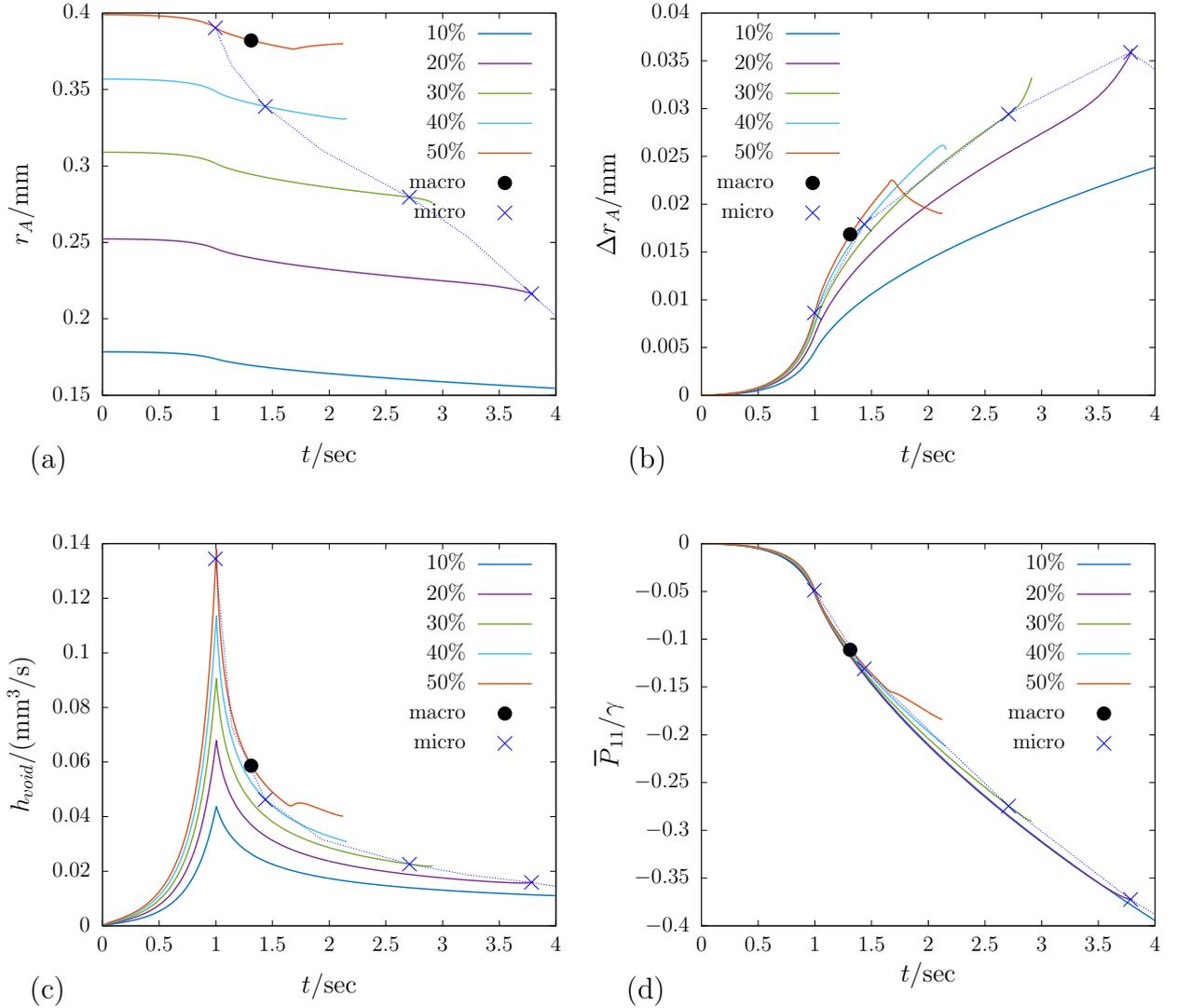

\begin{figure}%
\centering%
\normalsize%
\psfrag{t1} [c][c]{$t=0.5\, \text{sec}$}
\psfrag{t2} [c][c]{$t=1.0\, \text{sec}$}
\psfrag{t3} [c][c]{$t=1.5\, \text{sec}$}
\psfrag{t4} [c][c]{$t=2.12\, \text{sec}$}
\includegraphics*[scale=1.25]{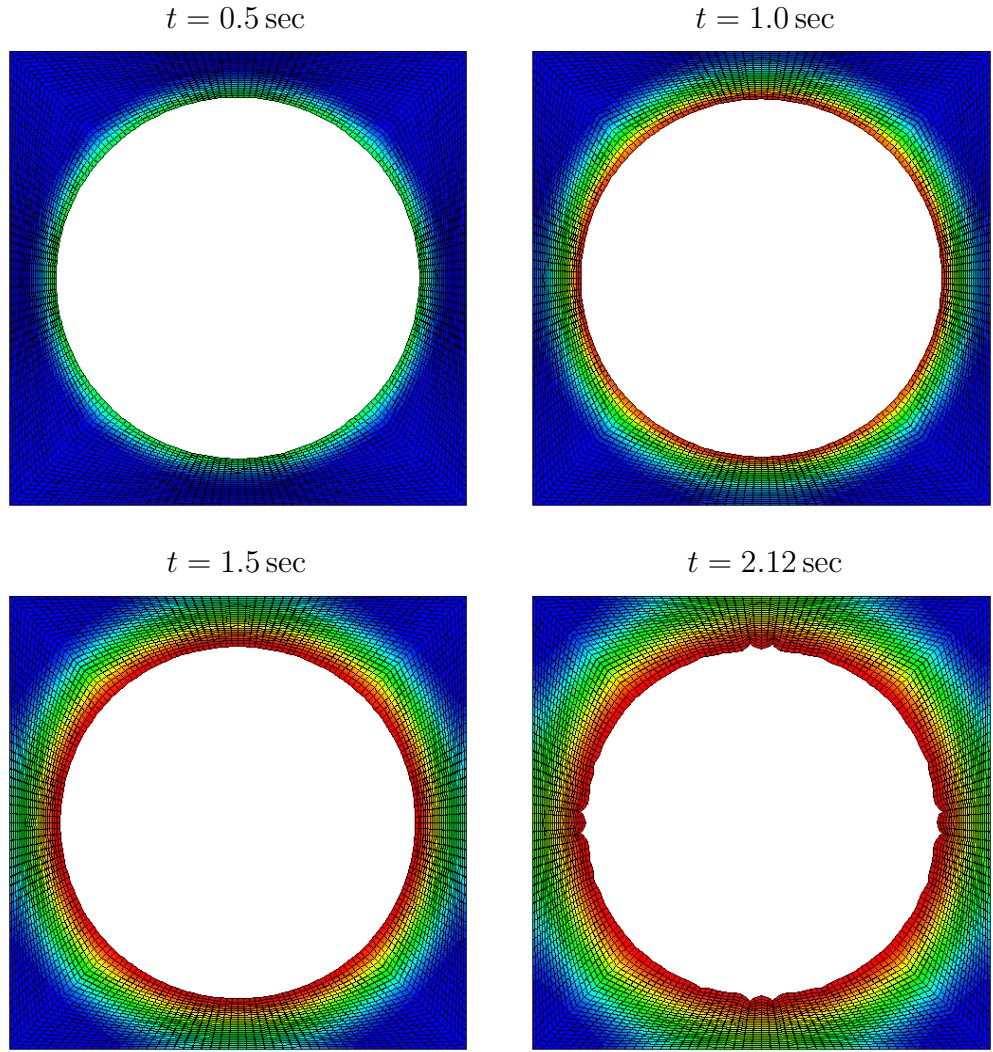}
\caption[]{
  \emph{Snapshots of the deformed hydrogel microstructure with the void-volume fraction $f_0=50\%$.} The diffusion of solvent volume takes place through the boundary of the void  $\partial\calH_0$, where we prescribe the chemical potential according to \eqref{chemBC}. The contour indicates the chemical potential, where blue corresponds to the minimum and red to the maximum value, respectively. We observe that at the time instant $t=2.12\, \text{sec}$, the void boundary has developed creases. It is observed that for this microstructure the strong ellipticity of the effective response is lost at $t=1.312\,\text{sec}$.}
\label{f50contour_mu}%
\end{figure}%

\begin{figure}[h]
  \hspace{-5mm}
  \scalebox{0.73}{\input{./mu_v-rn}}
  \hspace{-10mm}
  \scalebox{0.73}{\input{./mu_v-delRn}}
  \vspace{-9mm}
  \hspace{-5mm}
  
  { \text{(a)} \hspace{-5mm} \hspace{0.5\textwidth} \text{(b)}}
  
  \vspace{8mm}
  \hspace{-5mm}
  \scalebox{0.73}{\input{./mu_v-h_vn}}
  \hspace{-10mm}
  \scalebox{0.73}{\input{./mu_v-pbar11n}}
  \vspace{-9mm}
  \hspace{-3mm}
  
  { \text{(c)} \hspace{-5mm} \hspace{0.5\textwidth} \text{(d)}}  
    
  \caption{\emph{Response of hydrogel microstructures versus chemical potential prescribed at the internal boundaries of the voids.}
    \textbf{a)} Radius $r_A$;
    \textbf{b)} Change in the radius $\Delta r_A$ of the voids against the normalized chemical potential $\mu^N/\alpha$ at point $A$;
    \textbf{c)} Solvent volume entering the hydrogel matrix through the void boundary;
    \textbf{d)} Effective first Piola-Kirchhoff stress $\lbar{P}_{11}$ normalized with respect to the shear modulus of the matrix $\gamma$.}
  \label{fig:res-vs-mu_voided}  
\end{figure}
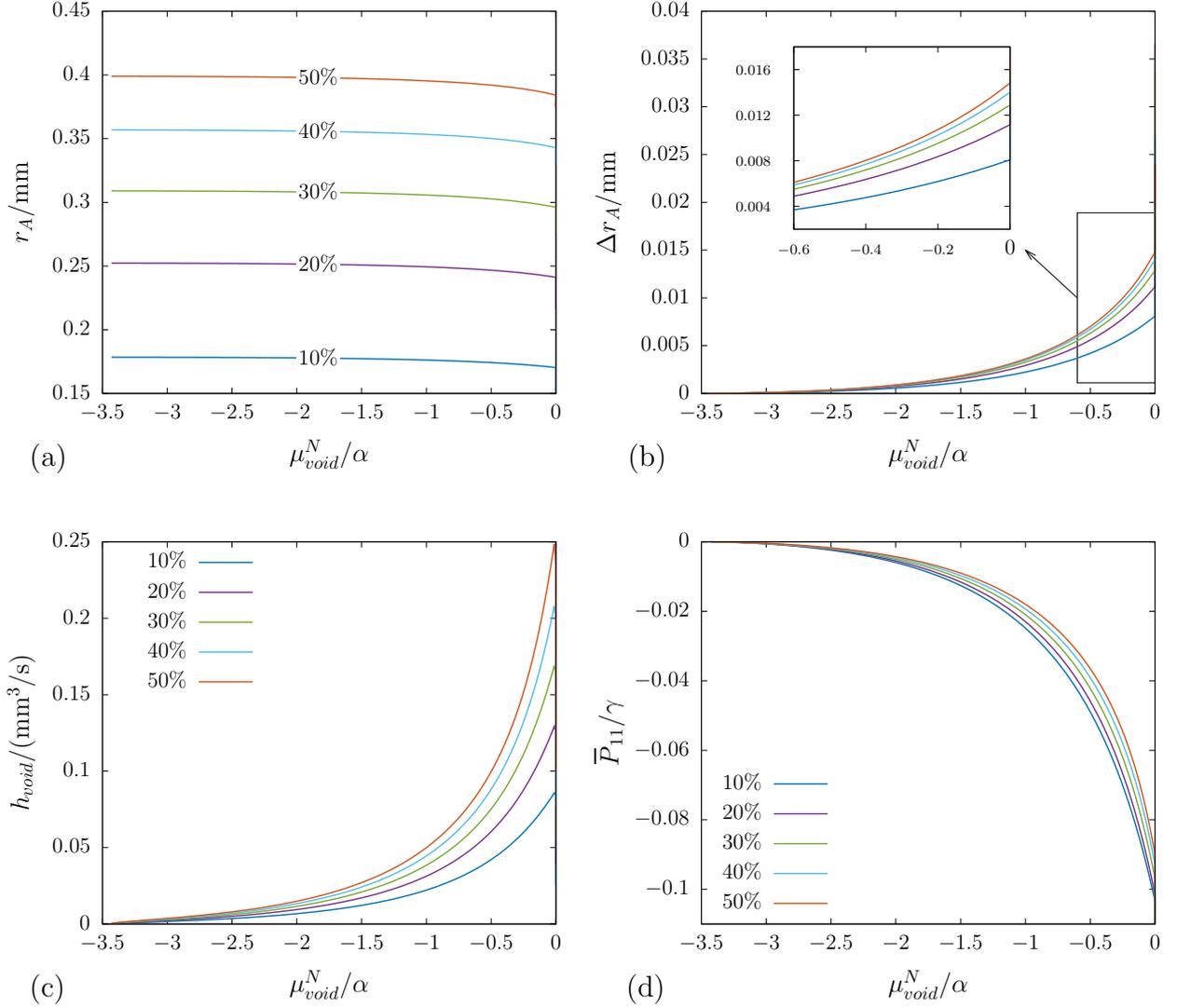

To obtain a first insight into the behavior of voided hydrogel microstructures, we analyze a sequence of RVEs with the void-volume fractions $f_0\in [10\%-50\%]$ and illustrate their response over the simulation time in Fig.~\ref{fig:res-vs-time_voided}. In Fig.~\ref{fig:res-vs-time_voided}\textbf{a} and \textbf{b}, we observe that the radius related to a given void-volume fraction $f_0$ changes significantly due to the influx of solvent, even after loading $t\ge 1\,\text{sec}$.
In Fig.~\ref{fig:res-vs-time_voided}\textbf{c}, we illustrate the solvent volume entering through the boundary of the voids over time according to $h_{void}=-\int_{\partial\calH_0}\KIH\cdot\BN\, \mathrm{d}A$ with $\BN$ being a unit normal vector pointing outward the hydrogel matrix. We observe that the influx of the solvent volume increases as we decrease the chemical potential at the void boundary until it finally reaches its maximum at $t=1\,\text{sec}$ and $\mu^N_{void}=0.0\,\text{N/mm\textsuperscript{2}}$.
In Fig.~\ref{fig:res-vs-time_voided}\textbf{d}, we have plotted the effective first Piola-Kirchhoff stress $\lbar{P}_{11}=\lbar{P}_{22}$ normalized with respect to the matrix shear modulus $\gamma$, see \eqref{FbarPbar}\textsubscript{2}. Notably, the normalized effective stress response is almost identical for all considered microstructures.

In Fig.~\ref{fig:res-vs-mu_voided}, we depict associated response functions over the applied chemical potential $\mu^N_{void}$ at the void boundary $\partial\calH_0$. While such an illustration gives insight into the chemically driven response, is not capable of describing the behavior beyond $t>1\,\text{sec}$. Therefore, the instability studies documented in the following sections will be discussed with respect to time rather than applied chemical potential.

\textbf{Creasing instabilities.} Before we come to the analysis of structural instabilities, we take a look at several snapshots of the unit-cell RVE with 50\% void-volume fraction during its deformation, see Fig.~\ref{f50contour_mu}. There we observe the development of \textit{creasing} during the swelling process. A subsequent analysis of effective moduli indicates a loss of strong ellipticity at $t = 1.312 \, \text{sec}$ (black bullet in Fig.~\ref{fig:res-vs-time_voided}), i.e., before creasing. In this context we note that according to \cite{mielke+sprenger98,cawte20} creasing is linked to \textit{quasi-convexity at the boundary} \cite{ball+marsden86}, which is a stronger condition than quasi-convexity---and consequently rank-one convexity---in the bulk. Thus, there may be a certain relation between the loss of strong ellipticity of effective moduli and the development of creasing at internal boundaries of the microstructure; see also \cite{silling91,yang+etal21,pandurangi+etal22}. However, an in-depth study of such a relation is beyond the scope of the present work.
%

%
%

%

\subsubsection{Instabilities in single-phase hydrogels with voids}
%
\label{sec43}

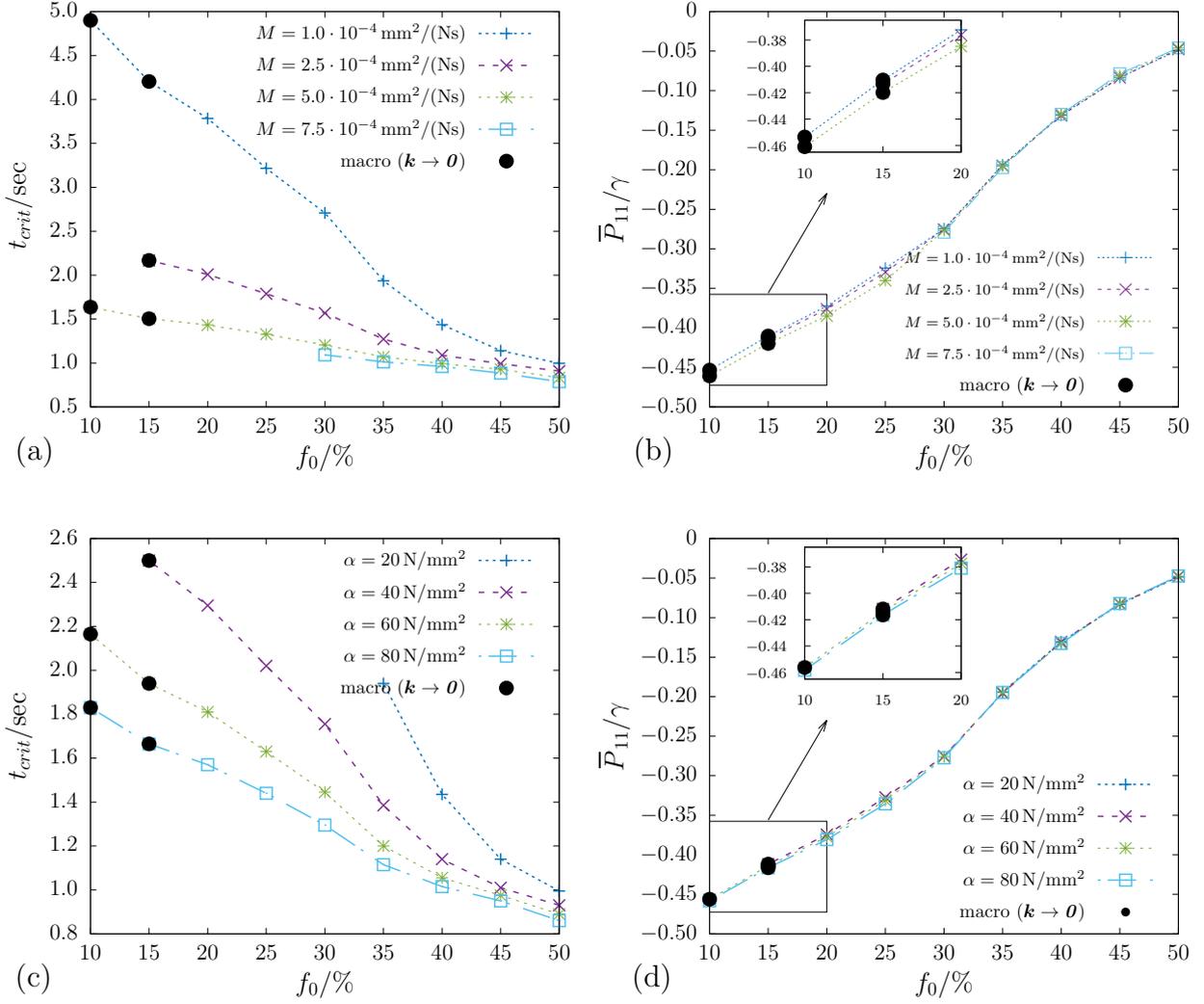
\begin{figure}[h]
  \hspace{-5mm}
  \scalebox{0.73}{\input{./f0-p11bar-vs-time}}
  \hspace{-10mm}
  \scalebox{0.73}{\input{./f0-p11bar-vs-M}}
  \vspace{-12mm}
  \hspace{-5mm}
  
  { \text{(a)} \hspace{-5mm} \hspace{0.5\textwidth} \text{(b)}}
  
  \vspace{8mm}
  \hspace{-5mm}
  \scalebox{0.73}{\input{./f0-vs-time-alpha}}
  \hspace{-10mm}
  \scalebox{0.73}{\input{./f0-vs-p11bar-alpha}}
  \vspace{-12mm}
  \hspace{-5mm}
  
  { \text{(c)} \hspace{-5mm} \hspace{0.5\textwidth} \text{(d)}}  
    
  \caption{\emph{Onset of instabilities depending on the mobility $M$ and the chemical parameter $\alpha$.} Results are plotted with respect to the critical time $t$ (\textbf{a} and \textbf{c}) and the normalized effective first Piola-Kirchhoff stress $\lbar{P}_{11}/\gamma_m$ (\textbf{b} and \textbf{d}) with $\gamma$ being the shear modulus of the matrix. We observe that while the mobility and the chemical parameters have significant influence on the critical time, the normalized critical stresses are all in the same range. Note that next to short-wavelength (pattern-transforming) instabilities, we observe long-wavelength instabilities for the void-volume fractions $f_0=10\%$ and $f_0=15\%$.}
  \label{2DBVPa}
\end{figure}

We now study the onset of instabilities for the voided hydrogel microstructures as described in Fig.~\ref{2D_BVPs}\textbf{a}. In particular, we analyze the influence of the void-volume fraction $f_0$, the mobility parameter $M$ and the chemical parameter $\alpha$ of the matrix.
As in the preceding studies, the remaining material parameters are taken from Table~\ref{table01}.

The onset of instabilities for the considered initial boundary value problem is depicted in Fig.~\ref{2DBVPa}\textbf{a--d}. In Fig.~\ref{2DBVPa}\textbf{a} and \textbf{c}, we plot the critical time versus the void-volume fraction depending on the mobility parameter $M$ and the chemical parameter $\alpha$, respectively. We observe that both parameters have a significant impact on the time at which a critical instability occurs. As expected, instabilities are initiated earlier for larger $M$ (that is, for faster diffusion through the matrix) and larger $\alpha$ (that is, for an enhanced chemo-mechanical interaction). In Fig.~\ref{2DBVPa}\textbf{b} and \textbf{d}, we plot the normalized critical first Piola-Kirchhoff stress $\lbar{P}_{11}=\lbar{P}_{22}$ at the instability points for different values of the mobility and the chemical parameter. We observe that the critical stresses associated with the various hydrogel microstructures are all in the same range. This indicates that the instability point is strongly correlated with the effective normalized compressive stress state. We further observe that the short-wavelength (microscopic) instabilities are primary for the microstructures with the void-volume fractions $f_0\in\{20\%-50\%\}$ and that the long-wavelength (macroscopic) instabilities are primary for the microstructures with the void-volume fractions $f_0\in\{10\%-15\%\}$. As a result of the microscopic instabilities, the hydrogel develops the so-called diamond-plate pattern, resulting in a change of periodicity of the microstructure as described in \cite{zhang+etal08,zhu2012capillarity,wu+etal14} and depicted in Fig.~\ref{buckling_modes_h002}\textbf{f}.

We note that the macroscopic instabilities were detected by means of Bloch--Floquet wave analysis with $\Bk\to\Bzero$ and coincide with the loss of strong ellipticity of the effective mechanical moduli as discussed in Section~\ref{sec32}.

\subsubsection{Instabilities of two-phase hydrogels with voids}
\label{sec44}

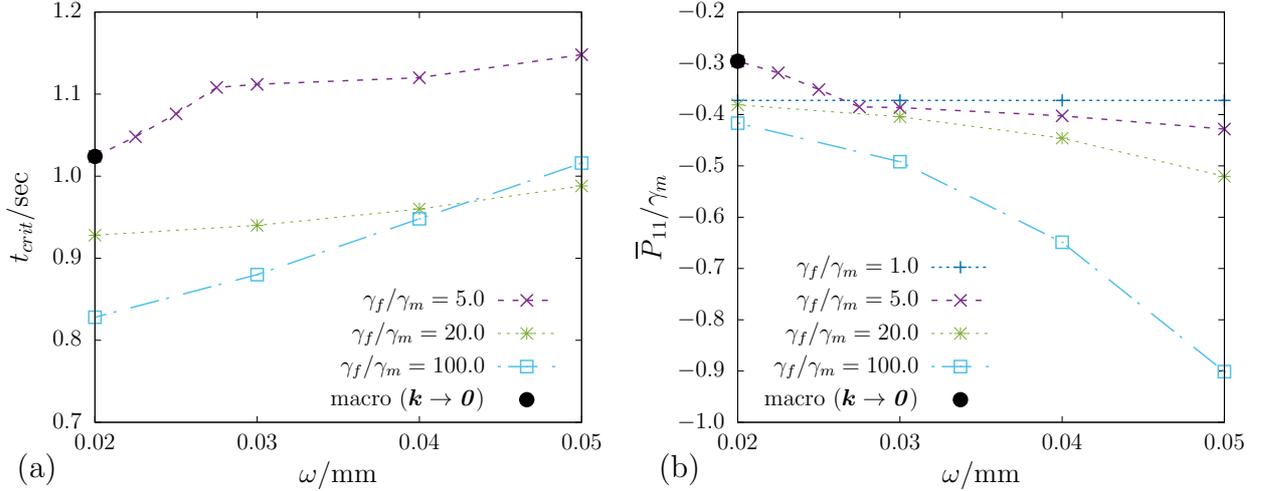
\begin{figure}[h]
  \hspace{-5mm}
  \scalebox{0.73}{\input{./2lay_h-vs-time_ga}}  
  \hspace{-10mm}
  \scalebox{0.73}{\input{./2lay_h-vs-pbar11_ga}}
  \vspace{-12mm}
  \hspace{-5mm}
    
  { \text{(a)} \hspace{-5mm} \hspace{0.5\textwidth} \text{(b)}}
  \caption{\emph{Onset of instabilities for the composite hydrogels with void-volume fraction $f_0=20\%$ for varying shear-modulus ratios $\gamma_c/\gamma_m$ of the coating ($c$) and the matrix ($m$) as well as coating thickness $\omega$.} The results are plotted with respect to \textbf{a)} the critical time $t$ and \textbf{b)} the normalized effective first Piola-Kirchhoff stress $\lbar{P}_{11}/\gamma_m$. Again, most instabilities induce pattern transformations from one unit-cell RVE to $2\times 2$ unit-cell RVEs. Furthermore, we observe macroscopic instabilities ($\Bk\to\Bzero$).}
\label{2DBVPb1}
\end{figure}
%
\begin{figure}[htb]%
\centering%
\footnotesize%
\psfrag{a} [c][c]{(a)}
\psfrag{b} [c][c]{(b)}
\psfrag{c} [c][c]{(c)}
\psfrag{d} [c][c]{(d)}
\psfrag{e} [c][c]{(e)}
\psfrag{f} [c][c]{(f)}
\psfrag{ga1} [c][c]{$\gamma_c/\gamma_m = 6.0$}
\psfrag{ga2} [c][c]{$\gamma_c/\gamma_m = 8.0$}
\psfrag{ga3} [c][c]{$\gamma_c/\gamma_m = 10.0$}
\psfrag{ga4} [c][c]{$\gamma_c/\gamma_m = 14.0$}
\psfrag{ga5} [c][c]{$\gamma_c/\gamma_m = 20.0$}
\psfrag{ga6} [c][c]{$\gamma_c/\gamma_m = \{30.0-100.0\}$}
\includegraphics*[scale=1.6]{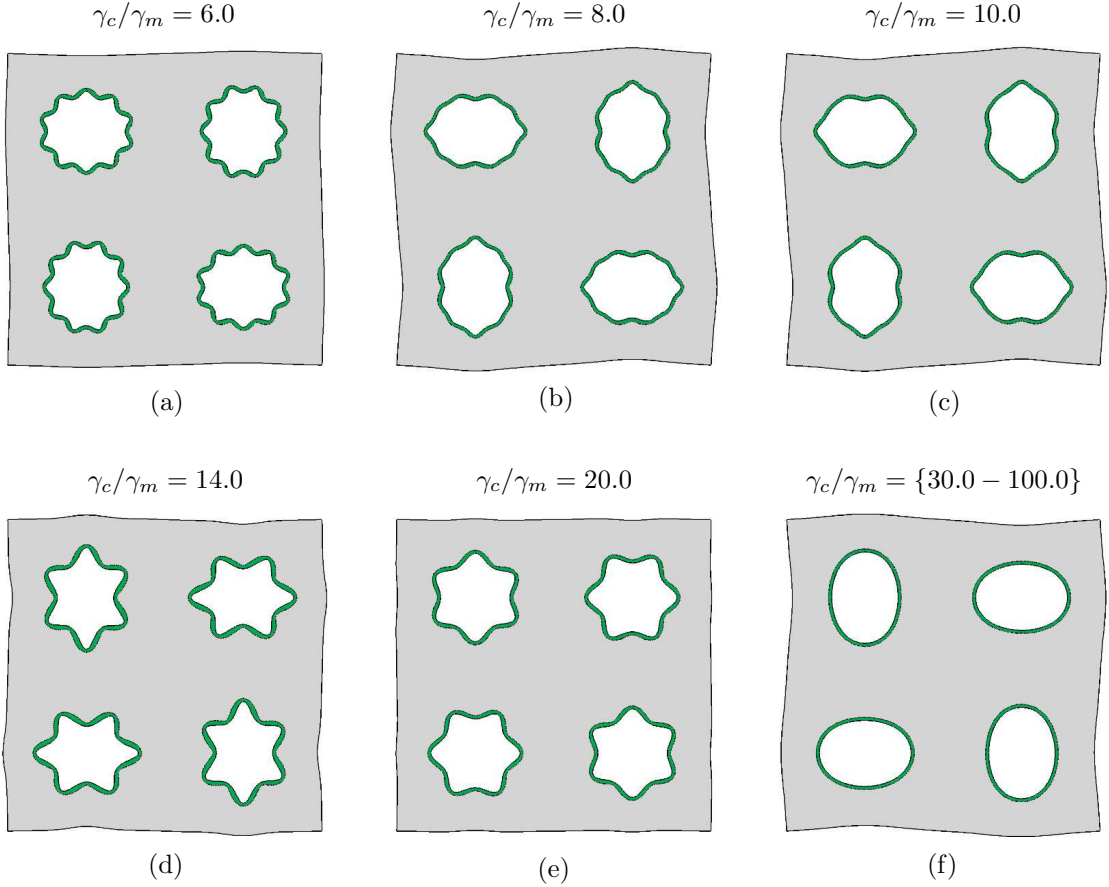}
\caption[]{
  \emph{Buckling modes of two-phase hydrogels depending on the shear modulus ratio of the coating ($c$) and the matrix ($m$).} The periodic hydrogels are composed of hydrogel matrix, hydrogel coatings with thickness $\omega=0.02\, \text{mm}$, and voids with volume fraction $f_0=20\%$. Depending on the shear-modulus ratio of coating and matrix $\gamma_c/\gamma_m$, pattern-transforming instabilities with and without wrinkling of internal surfaces can be observed.}
\label{buckling_modes_h002}%
\end{figure}%
%
%
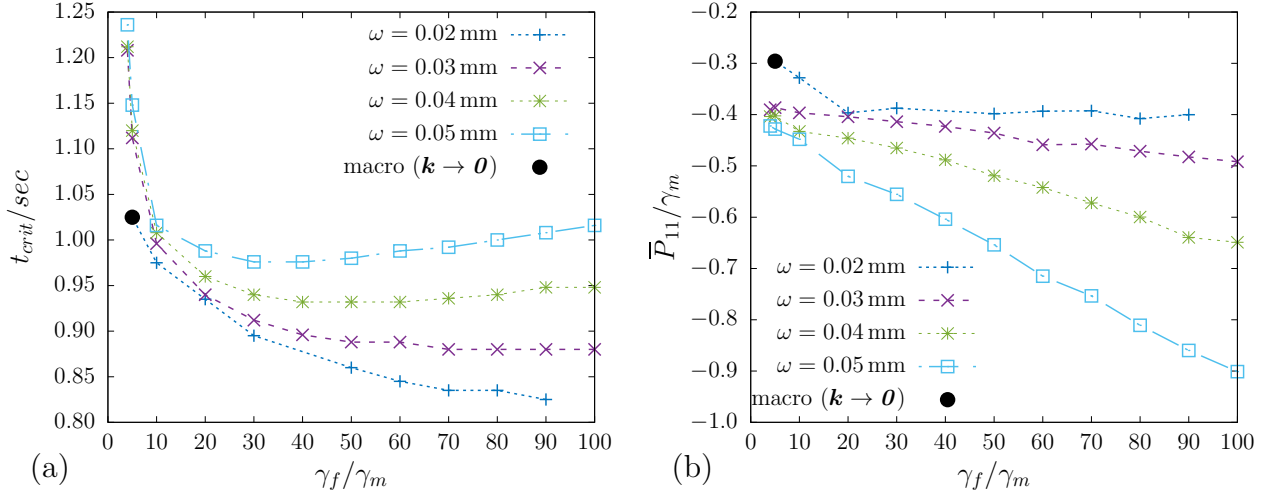
\begin{figure}[h]
  \hspace{-5mm}
  \scalebox{0.73}{\input{./2lay_gamma-vs-time_h}}  
  \hspace{-10mm}
  \scalebox{0.73}{\input{./2lay_gamma-vs-pbar11_h}}  
  \vspace{-12mm}
  \hspace{-5mm}
    
  { \text{(a)} \hspace{-5mm} \hspace{0.5\textwidth} \text{(b)}}
  \caption{\emph{Onset of instabilities for the composite hydrogels with the void-volume fraction $f_0=20\%$ depending on the shear modulus ratio $\gamma_c/\gamma_m$ of the coating ($c$)  and the matrix ($m$) as well as the coating thickness $\omega$.} The results are plotted with respect to \textbf{a)} the critical time $t$ and \textbf{b)} the normalized effective first Piola-Kirchhoff stress $\lbar{P}_{11}/\gamma_m$. The short-wavelength buckling patterns span $2\times 2$ unit-cell RVEs.}
  \label{2DBVPb2}
\end{figure}
%
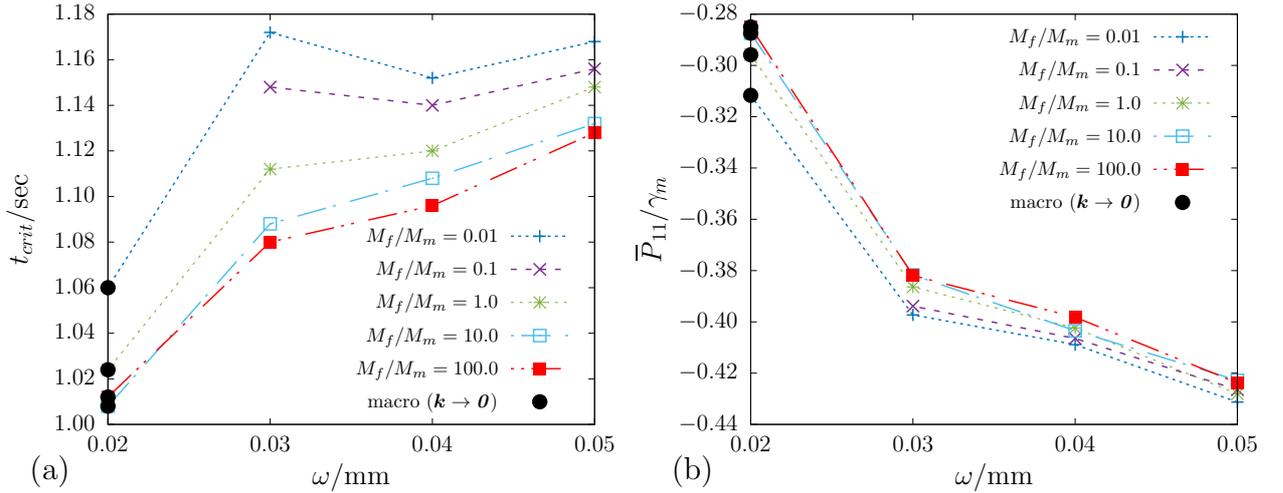
\begin{figure}[h]
  \hspace{-5mm}
  \scalebox{0.73}{\input{./2lay_time-vs-h_M}}  
  \hspace{-10mm}
  \scalebox{0.73}{\input{./2lay_p11bar-vs-h_M}}  
  \vspace{-12mm}
  \hspace{-5mm}
    
  { \text{(a)} \hspace{-5mm} \hspace{0.5\textwidth} \text{(b)}}
  \caption{\emph{Onset of instabilities for the composite hydrogels with the shear-modulus ratio $\gamma_c/\gamma_m=5$ and the void-volume fraction $f_0=20\%$ depending on the mobility-parameter ratio $M_m/M_c$ of the coating ($c$) and the matrix ($m$) as well as the coating thickness $\omega$.} The results are plotted with respect to \textbf{a)} the critical time $t$ and \textbf{b)} the normalized effective first Piola-Kirchhoff stress $\lbar{P}_{11}/\gamma_m$. The pattern-transforming instabilities span $2\times 2$ unit-cell RVEs.}
\label{2DBVPb3}  
\end{figure}

We proceed with analyzing instabilities of two-phase voided hydrogel microstructures as described in Fig.~\ref{2D_BVPs}\textbf{b}. As can be seen in this figure, we take into account microstructures composed of hydrogel matrix $(m)$, voids ($v$) and hydrogel coatings $(c)$ surrounding the voids.

In the following, we investigate the influence of the thickness of the coating as well as of the ratios between the shear moduli and the mobility parameters of the hydrogel coating and the hydrogel matrix on the formation of instabilities. Here, we take into account a fixed void-volume fraction $f_0=20\%$ and boundary conditions as described in  Section~\ref{sec41}. Material parameters for the matrix and the coating are listed in Table~\ref{table01}.

In Fig.~\ref{2DBVPb1}, we illustrate the critical time $t_{crit}$ and the effective first Piola-Kirchhoff stress $\lbar{P}_{11}=\lbar{P}_{22}$ versus the coating thickness $\omega$ for various values of shear-modulus ratio $\gamma_c/\gamma_m$ of the coating and the matrix. In all studies we assume that the mechanical stiffness of the coating is higher than that of the matrix $\gamma_m=0.01\,\text{N/mm\textsuperscript{2}}$. In the present example the mobility parameter of both matrix and coating is set to $M_m=M_c=10^{-4}\,\text{mm\textsuperscript{4}/(Ns)}$.

From Fig.~\ref{2DBVPb1}\textbf{a}, we observe that the shear-modulus ratio $\gamma_c/\gamma_m$ has a significant influence on the critical instability point, such that microscopic instabilities occur much earlier than in case of the corresponding single-phase voided microstructures. Similarly, the critical effective stresses strongly depend on the shear-modulus ratio, see Fig.~\ref{2DBVPb1}\textbf{b}. In that context we observe in particular that for $\gamma_c/\gamma_m=5.0$ and the coating thickness $\omega<0.025$, instabilities occur much earlier at much lower critical stresses.

We note that while most of the triggered instabilities show a critical periodicity of $2\times 2$ unit cells, the observed microcopic patterns are highly dependent on the thickness of the coating and the ratio of shear moduli, see Fig.~\ref{buckling_modes_h002}. In particular, two distinct kinds of pattern transformations given by {(i)} micro-wrinkling and {(ii)} diamond-plate patterns (mutual deformation of voids into ellipses) are recorded. In the first case, pattern transformations due to wrinkling of the coating are detected for $\gamma_c/\gamma_m< 30$. Here, the number of micro-wrinkles can be tuned by means of the shear moduli of the matrix and the coating. When the shear-modulus ratio exceeds the value of $\gamma_c/\gamma_m = 30$, diamond-plate patterns with smooth internal boundaries as in case of the single-phase hydrogel structures are observed.%
\footnote{%
We note that while the Bloch--Floquet analysis indicated a long-wavelength instability for the shear-modulus ratio $\gamma_c/\gamma_m= 5.0$ (see black bullet in Fig.~\ref{2DBVPb1}), the analysis of effective moduli did not show a loss of strong ellipticity. A further analysis of strong ellipticity beyond the critical point was not possible due to missing convergence. Refining the time-step size did not help to overcome this issue.
We refer to \cite{dortdivanlioglu+linder19}, who observed a competition between short- and long-wavelength instabilities for small shear-modulus ratios in plane bilayer systems. Their results indicate that for a refined mesh, the long-wavelength instabilities become the primary ones.%
}

In Fig.~\ref{2DBVPb2}, we further study instabilities for a wide range of shear-modulus ratios and coating thicknesses $\omega\in\{0.02,\,0.03,\,0.04,\,0.05\}\,\text{mm}$. The graphs describe the critical time $t_{crit}$ and the effective first Piola-Kirchhoff stress $\lbar{P}_{11}$ depending on the shear-modulus ratio $\gamma_c/\gamma_m$. Similar to the previous cases, we observe that for small values $\gamma_c/\gamma_m<10$, the critical instability time increases exponentially. As the stiffness of the coating increases, the corresponding compressive stresses rise and induce earlier buckling. A slightly different response is observed for the thicker and stiffer coatings due to a stabilizing effect, see the intersection of corresponding curves in Fig.~\ref{2DBVPb1}\textbf{a}.

In Fig.~\ref{2DBVPb3}, we study the influence of various mobility-parameter ratios $M_c/M_m$ of the coating and the matrix on the onset and the type of instabilities in consideration of a shear-modulus ratio of $\gamma_c/\gamma_m=5.0$ with $\gamma_m=0.01\,\text{N/mm\textsuperscript{2}}$ and a mobility parameter of the coating $M_c=10^{-4}\,\text{mm\textsuperscript{4}/(Ns)}$. We illustrate the critical time $t_{crit}$ and the effective first Piola-Kirchhoff stress $\lbar{P}_{11}=\lbar{P}_{22}$ depending on the thickness $\omega$ of the coating. Since the mobility parameter is associated with the diffusivity of the hydrogel, the solvent tends to accumulate in the coating when $M_c/M_m>1$. As a result, the coating reaches the critical stress faster and the instabilities occur at lower critical time, see Fig.~\ref{2DBVPb3}\textbf{a}. In Fig.~\ref{2DBVPb3}\textbf{b}, we depict the corresponding effective critical stresses over the coating thickness $\omega$ depending on the mobility-parameter ratio. In all considered problems, we observe the buckling instabilities resulting in $2\times 2$ unit-cell RVEs. For a coating thickness of $\omega=0.02\,\text{mm}$, the long-wavelength instabilities are detected again.

%
\subsection{Three-dimensional hydrogel films}
\label{sec45}

The previous studies were limited to a two-dimensional plane-strain scenario. We now extend the analysis to the three-dimensional case and investigate structural instabilities in hydrogel thin films. The considered films comprise cylindrical hydrogel inclusions that are embedded into a hydrogel matrix in a periodic manner, see Fig.~\ref{periodBC3}. We assume that the diffusivity of the inclusions is low compared with the matrix, so that swelling of the inclusions is not pronounced. The considered setup is motivated from the experimental works \cite{wang+etal17,ma+etal19}.

Before we provide the results of a parametric study in consideration of various geometrical properties of microstructures, we discuss suitable boundary conditions for the numerical realization of three-dimensional hydrogel films considering plane periodicity. In that context, we also comment on relevant ingredients of the numerical discretization and finite-element implementation. 

%
\subsubsection{Boundary conditions in the three-dimensional setting}
%
\label{sec451}
%
\begin{figure}[htb]%
\centering%
\footnotesize%
\psfrag{D0}  [l][l]{$\calD$}
\psfrag{B0}  [l][l]{$\calB$}
\psfrag{H}   [c][c]{$\calH$}
\psfrag{bB3d}   [c][c]{$\partial\calB_0^{side}$}
\psfrag{bBtop}   [r][r]{$\partial\calB_0^{top}$}
  \psfrag{dS}  [l][l]{$\partial\calB_t$}
  \psfrag{X}   [c][c] {$\BX$}
  \psfrag{x}   [c][c] {$\Bx$}
  \psfrag{F}[c][c] {$\BF$}
  \psfrag{phi}[c][c] {$\Bvarphi(\BX),\,\Ba(\BX)$}
  \psfrag{N+}[c][c] {$\BN^+$}
  \psfrag{N-}[c][c] {$\BN^-$}
  \psfrag{n+}[c][c] {$\Bn^+$}
  \psfrag{n-}[c][c] {$\Bn^-$}
  \psfrag{Xm}[c][c] {$\BX^-$}
  \psfrag{Xp}[c][c] {$\BX^+$}
  \psfrag{xm}[c][c] {$\Bvarphi^-$}
  \psfrag{xp}[c][c] {$\Bvarphi^+$}
  \psfrag{Xj}[c][c] {$\jump{\BX}$}  
  \psfrag{xj}[c][c] {$\jump{\Bvarphi}$}
  \psfrag{xjp}[c][c] {$\jump{\Bvarphi}$}
  \psfrag{xmp}[c][c] {$\Bvarphi^-$}
  \psfrag{xpp}[c][c] {$\Bvarphi^+$}
  \psfrag{aa}[c][c] {\text{(a)}}
  \psfrag{bb}[c][c] {\text{(b)}}
  \psfrag{a}[c][c] {$a$}
  \psfrag{b}[c][c] {$b$}
  \psfrag{alpha}[l][l] {$\beta$}      
  \psfrag{micro}[c][c]{\shortstack{microstructure of \\periodic hydrogel film}}
  \psfrag{unit}[c][c]{unit-cell RVE}  
\includegraphics*[scale=0.6]{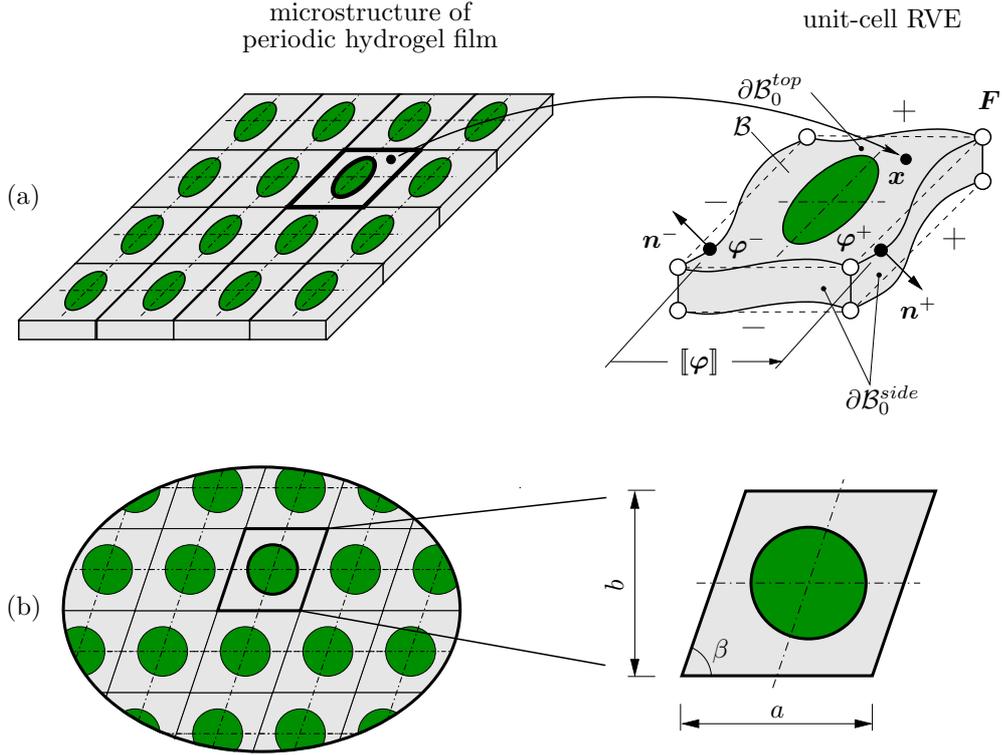}
\caption[Description of computational homogenization]{%
\textit{%
Boundary value problem for the three-dimensional composite hydrogel thin films}.
\textbf{a)} The periodic microstructure of a thin film is represented by a unit-cell RVE, where the jump conditions for the deformation map and the solvent-volume flux refer only to the side boundaries, see \eqref{3Djumpbc}. In the three-dimensional setting it is assumed that the solvent diffuses through the top and the bottom surfaces, which is realized by applying associated boundary conditions for the chemical potential. \textbf{b)} A top view of the thin film indicated that the periodic microstructure can be described by means of unit-cell RVEs given by right parallelogrammic prism with cylindrical inclusions. In all considered examples, we assume that the inclusions have lower diffusivity and higher mechanical stiffness than the matrix.}
\label{periodBC3}%
\end{figure}%

Because of the periodic nature of the problem, the following computations can again be reduced to a selected unit-cell RVE. In this connection, we consider the jump boundary conditions across the lateral boundaries of the RVE
\begin{equation}
  \label{3Djumpbc}
  \jump{\Bvarphi} = \lhat{\lbar{\BF}} \jump{\BX}
  \AND
  \jump{\KIH}\cdot\BN=0 \quad \text{on} \quad \partial\calB_0^{\textrm{side}}\, ,
\end{equation}  
where $\lhat{\lbar{\BF}}$ is the effective deformation gradient in the plane of periodicity, see also Fig.~\ref{periodBC3}\textbf{a}. We assume that the hydrogel is mechanically constrained in plane and that no stresses arise in out-of-plane direction, that is
\begin{equation}
  \label{3Djumpbc-2}
    \lhat{\lbar{\BF}}
    =
    \left[
    \begin{array}{ccc}
        1 & 0 & 0     \\
        0 & 1 & 0     \\
        0 & 0 & \cdot \\
    \end{array}
    \right]
\AND
\BP \cdot \BN = \Bzero
\ \mbox{on} \
\partial\calB_0^{top,bottom} \, .
\end{equation}
On the chemical side, the solvent is assumed to enter the hydrogel only through the top and the bottom surfaces of the film, which is reflected by the chemical-potential boundary condition
\begin{equation}
\label{chemBC3D}
\mu = \lambda(t)\mu_0
\quad \text{on} \quad
\partial\calB_0^{top,bottom}
\quad \text{with} \quad
\lambda(t) =
  \begin{dcases}
    1-t\, &\text{for} \quad  t\le 1\, \text{sec} \\
    0  \, &\text{for} \quad  t >  1\, \text{sec} 
  \end{dcases}  \, ,
\end{equation}
where $\mu_0$ is the chemical potential of the reference state, refer to \eqref{mu0}. In what follows, we consider a time-step size of $\tau=4\cdot 10^{-3}\,\text{sec}$.

The geometry of the boundary value problem and the corresponding unit-cell RVE are described in Fig.~\ref{periodBC3}\textbf{b}. Here, the film has the thickness $\omega=0.02\, \text{mm}$ and plane dimensions of $a\times b$ with a fixed length of $a=1.0\,\text{mm}$. The inclusions with radius $r=0.25\,\text{mm}$ have lower diffusivity and higher elastic stiffness than the matrix ($M_i/M_m=0.01$, $\gamma_i/\gamma_m=10$). The remaining material parameters of both phases are taken from Table~\ref{table01}.

The microstructures described in Fig.~\ref{periodBC3}\textbf{b} are discretized using hexahedral finite elements of the type H$_1$H$_1$. Here, the displacement field and the solvent-volume flux are assumed as nodal degrees of freedom resulting in six degrees of freedom per node. Depending on the morphology, $5,500$ to $9,000$ finite elements ($8,500$ to $14,000$ nodes) are needed for the discretization of a unit-cell RVE.

%
\subsubsection{Instabilities of three-dimensional two-phase hydrogel thin films}
\label{sec452}
\begin{figure}[t]
  \hspace{-5mm}
  \scalebox{0.73}{\input{./b-vs-time-alpha}}  
  \hspace{-10mm}
  \scalebox{0.73}{\input{./b-vs-pbar-alpha}}  
  \vspace{-12mm}
  \hspace{-5mm}
    
  { \text{(a)} \hspace{-5mm} \hspace{0.5\textwidth} \text{(b)}}
  \caption{\emph{Onset of instabilities for the three-dimensional composite hydrogel films depending on their morphology.} A shear-modulus ratio of $\gamma_i/\gamma_m=10$ between inclusions and matrix as well as a morphology with $a=1.0\,\text{mm}$ and $r=0.25\,\text{mm}$ are considered. Instabilities are depicted with respect to \textbf{a)} the critical time $t$ and \textbf{b)} the normalized effective first Piola-Kirchhoff stress $\lbar{P}_{11}/\gamma_m$.}
\label{3DBVPb1}  
\end{figure}
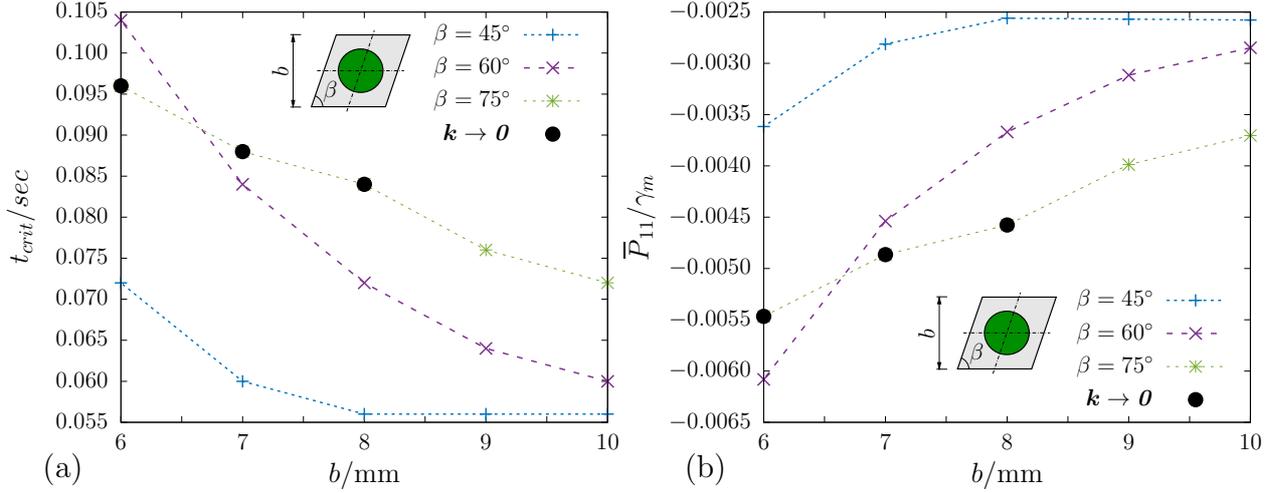
%

%
\begin{figure}[htb]%
\centering%
\footnotesize%
\psfrag{a} [c][c]{(a)}
\psfrag{b} [c][c]{(b)}
\psfrag{c} [c][c]{(c)}
\psfrag{d} [c][c]{(d)}
\psfrag{e} [c][c]{(e)}
\psfrag{f} [c][c]{(f)}
\psfrag{min} [l][l]{min}
\psfrag{max} [l][l]{max}
\psfrag{x1} [l][l]{$X_1$}
\psfrag{x2} [r][r]{$X_2$}
\psfrag{x3} [r][r]{$X_3$}
\psfrag{zz} [c][c]{$\delta{v}_3$}
\psfrag{ga1} [c][c]{$b = 0.6\,\text{mm}$}
\psfrag{ga2} [c][c]{$b = 0.8\,\text{mm}$}
\psfrag{ga3} [c][c]{$b = 1.0\,\text{mm}$}
\psfrag{ga4} [c][c]{$b = 0.6\,\text{mm}$}
\psfrag{ga5} [c][c]{$b = 0.8\,\text{mm}$}
\psfrag{ga6} [c][c]{$b = 1.0\,\text{mm}$}
\includegraphics*[scale=1.45]{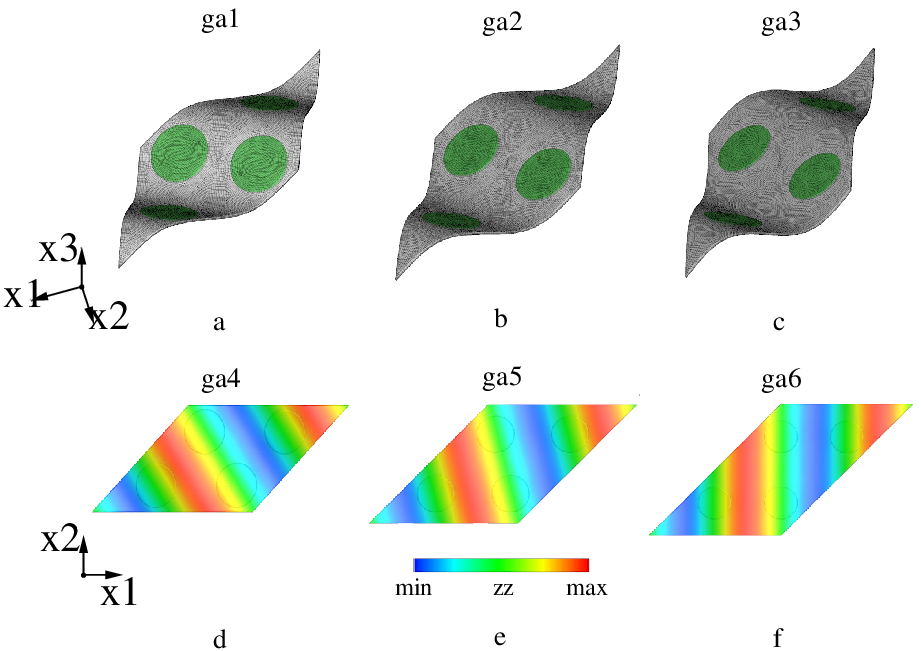}
\caption[]{
  \emph{Buckling modes of composite hydrogels with $\beta=45^\circ$ depending on the width $b$ of the unit-cell RVE.}
  \textbf{a)--c)} Different buckling patterns spanning $2\times 2$ unit-cell RVEs;
  \textbf{d)--e)} Contour plots of the critical eigenvectors ($\delta v_3$; top view) reveal the change in the direction of the buckling depending on the width of the unit-cell RVE.%
}
\label{buckling_modes_3D}%
\end{figure}%
%
%
%
\begin{figure}[htb]%
\centering%
\footnotesize%
\psfrag{a} [c][c]{(a)}
\psfrag{b} [c][c]{(b)}
\psfrag{c} [c][c]{(c)}
\psfrag{d} [c][c]{(d)}
\psfrag{e} [c][c]{(e)}
\psfrag{f} [c][c]{(f)}
\psfrag{min} [l][l]{min}
\psfrag{max} [l][l]{max}
\psfrag{x1} [l][l]{$X_1$}
\psfrag{x2} [r][r]{$X_2$}
\psfrag{x3} [r][r]{$X_3$}
\psfrag{ga1} [c][c]{$b = 0.6\,\text{mm}$}
\psfrag{ga2} [c][c]{$b = 0.8\,\text{mm}$}
\psfrag{ga3} [c][c]{$b = 1.0\,\text{mm}$}
\psfrag{ga4} [c][c]{$b = 0.6\,\text{mm}$}
\psfrag{ga5} [c][c]{$b = 0.8\,\text{mm}$}
\psfrag{ga6} [c][c]{$b = 1.0\,\text{mm}$}
\psfrag{v3} [c][c]{$\delta{v}_3$}
\includegraphics*[scale=1.45]{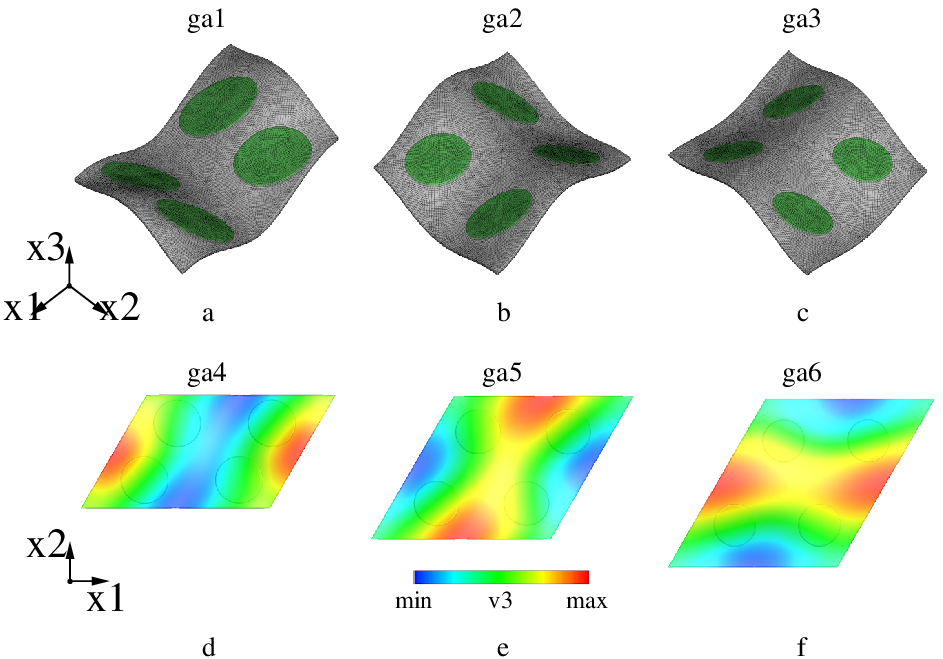}
\caption[]{
  \emph{Buckling modes of composite hydrogels with $\beta=60^\circ$ depending on the width $b$ of the unit-cell RVE.}
  \textbf{a)--c)} Different buckling patterns spanning $2\times 2$ unit-cell RVEs;
  \textbf{d)--e)} Contour plots of the critical eigenvectors ($\delta v_3$; top view) reveal the change in the direction of the buckling depending on the width of the unit-cell RVE.%
}
\label{buckling_modes_3D-2}%
\end{figure}%
%
\begin{figure}[htb]%
\centering%
\footnotesize%
\psfrag{a} [c][c]{(a)}
\psfrag{b} [c][c]{(b)}
\includegraphics*[scale=0.65]{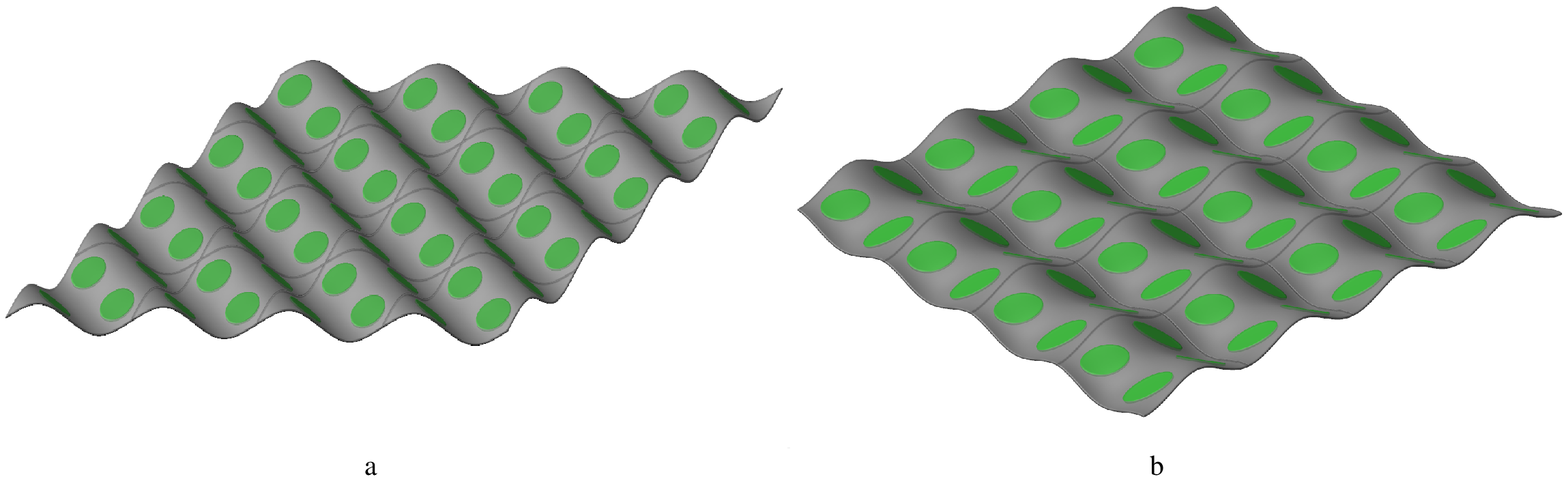}
\caption[]{
  \emph{Buckling modes of composite hydrogel microstructures composed of unit-cell RVEs with width $b=0.8\, \text{mm}$.}
  \textbf{a)} We observe wrinkling of the microstructure for $\beta=45^\circ$ and \textbf{b)} cooperative buckling for $\beta=60^\circ$.%
}
\label{3Dmicrostructure}%
\end{figure}%

In Fig.~\ref{3DBVPb1}, we illustrate the onset of instabilities depending on the morphological parameters of the RVE. We observe that both the width of the unit cell $b$ as well as the angle $\beta$ have a significant influence on the critical time and the critical stress. Since the composite effectively becomes softer with increasing $b$, instabilities occur earlier (cf. Euler buckling of beams).
A similar behavior can be observed for decreasing angle $\beta$, which results in an effectively larger distance between inclusions in the $X_2$-direction for a fixed value of $b$.
While in most cases we observe buckling into a configuration described by $2 \times 2$ unit-cell RVEs, some computations indicate a transition from short-wavelength to long-wavelength instabilities. These points are indicated by black bullets in Fig.~\ref{3DBVPb1}.

In Fig.~\ref{buckling_modes_3D}, we illustrate representative buckling patterns for different widths $b$ of the film for an angle $\beta=45^\circ$ (see Fig.~\ref{buckling_modes_3D-2} for $\beta=60^\circ$). In case of $\beta=45^\circ$, we observe \textit{wrinkling} of the microstructure, see Fig.~~\ref{buckling_modes_3D}\textbf{a--c}, where the direction of the wrinkling pattern can be influenced by adjusting the width $b$ of the unit-cell RVE (see Fig.~~\ref{buckling_modes_3D}\textbf{d--f}). For $\beta=60^\circ$, we observe substantially different buckling patterns in the form of cooperative saddle-point shapes, see Fig.~\ref{buckling_modes_3D-2}. For better insight into the overall buckling patterns, we refer to Fig.~\ref{3Dmicrostructure}, where we plot larger domains of the buckled films.

\section{Conclusion}
\label{sec50}

In the present work, we studied instabilities in two-dimensional and three-dimensional composite hydrogels with periodic microstructures.
In two dimensions, we investigated the onset of instabilities for single-phase voided microstructures as well as two-phase coated microstructures. We observed that in case of single-phase hydrogels, the periodicity of the microstructure may change from a unit-cell RVE to a larger RVE containing $2\times 2$ unit cells,
forming the well-known diamond-plate pattern frequently observed in experiments. For two-phase composites with thin coating surrounding the voids, we observed that the buckling modes can be substantially different from those of the single-phase hydrogel microstructures. Here, the change in periodicity and the development of patterns can be attributed either to wrinkling of the coating or to cooperative deformation of voids into ellipses. A detailed parametric study revealed that the observed patterns depend on the thickness of the coating as well as on the material parameters of coating and matrix.
Our simulations of three-dimensional composite hydrogel films revealed various buckling patterns depending on the arrangement of small reinforcing particles within the hydrogel matrix of the microstructure. These findings could provide access to the design of complex out-of-plane buckling patterns to be tuned by adjustment of microstructures.

\section*{Acknowledgements}

Funded by Deutsche Forschungsgemeinschaft (DFG, German Research Foundation) under Germany's Excellence Strategy -- EXC 2075 -- 390740016. This funding is gratefully acknowledged.

\end{document}

%% file: time-r.tex
\begingroup
  \makeatletter
  \providecommand\color[2][]{%
    \GenericError{(gnuplot) \space\space\space\@spaces}{%
      Package color not loaded in conjunction with
      terminal option `colourtext'%
    }{See the gnuplot documentation for explanation.%
    }{Either use 'blacktext' in gnuplot or load the package
      color.sty in LaTeX.}%
    \renewcommand\color[2][]{}%
  }%
  \providecommand\includegraphics[2][]{%
    \GenericError{(gnuplot) \space\space\space\@spaces}{%
      Package graphicx or graphics not loaded%
    }{See the gnuplot documentation for explanation.%
    }{The gnuplot epslatex terminal needs graphicx.sty or graphics.sty.}%
    \renewcommand\includegraphics[2][]{}%
  }%
  \providecommand\rotatebox[2]{#2}%
  \@ifundefined{ifGPcolor}{%
    \newif\ifGPcolor
    \GPcolortrue
  }{}%
  \@ifundefined{ifGPblacktext}{%
    \newif\ifGPblacktext
    \GPblacktexttrue
  }{}%
  \let\gplgaddtomacro\g@addto@macro
  \gdef\gplbacktext{}%
  \gdef\gplfronttext{}%
  \makeatother
  \ifGPblacktext
    \def\colorrgb#1{}%
    \def\colorgray#1{}%
  \else
    \ifGPcolor
      \def\colorrgb#1{\color[rgb]{#1}}%
      \def\colorgray#1{\color[gray]{#1}}%
      \expandafter\def\csname LTw\endcsname{\color{white}}%
      \expandafter\def\csname LTb\endcsname{\color{black}}%
      \expandafter\def\csname LTa\endcsname{\color{black}}%
      \expandafter\def\csname LT0\endcsname{\color[rgb]{1,0,0}}%
      \expandafter\def\csname LT1\endcsname{\color[rgb]{0,1,0}}%
      \expandafter\def\csname LT2\endcsname{\color[rgb]{0,0,1}}%
      \expandafter\def\csname LT3\endcsname{\color[rgb]{1,0,1}}%
      \expandafter\def\csname LT4\endcsname{\color[rgb]{0,1,1}}%
      \expandafter\def\csname LT5\endcsname{\color[rgb]{1,1,0}}%
      \expandafter\def\csname LT6\endcsname{\color[rgb]{0,0,0}}%
      \expandafter\def\csname LT7\endcsname{\color[rgb]{1,0.3,0}}%
      \expandafter\def\csname LT8\endcsname{\color[rgb]{0.5,0.5,0.5}}%
    \else
      \def\colorrgb#1{\color{black}}%
      \def\colorgray#1{\color[gray]{#1}}%
      \expandafter\def\csname LTw\endcsname{\color{white}}%
      \expandafter\def\csname LTb\endcsname{\color{black}}%
      \expandafter\def\csname LTa\endcsname{\color{black}}%
      \expandafter\def\csname LT0\endcsname{\color{black}}%
      \expandafter\def\csname LT1\endcsname{\color{black}}%
      \expandafter\def\csname LT2\endcsname{\color{black}}%
      \expandafter\def\csname LT3\endcsname{\color{black}}%
      \expandafter\def\csname LT4\endcsname{\color{black}}%
      \expandafter\def\csname LT5\endcsname{\color{black}}%
      \expandafter\def\csname LT6\endcsname{\color{black}}%
      \expandafter\def\csname LT7\endcsname{\color{black}}%
      \expandafter\def\csname LT8\endcsname{\color{black}}%
    \fi
  \fi
    \setlength{\unitlength}{0.0500bp}%
    \ifx\gptboxheight\undefined%
      \newlength{\gptboxheight}%
      \newlength{\gptboxwidth}%
      \newsavebox{\gptboxtext}%
    \fi%
    \setlength{\fboxrule}{0.5pt}%
    \setlength{\fboxsep}{1pt}%
\begin{picture}(7200.00,5040.00)%
    \gplgaddtomacro\gplbacktext{%
      \csname LTb\endcsname
      \put(982,704){\makebox(0,0)[r]{\strut{}$0.15$}}%
      \put(982,1549){\makebox(0,0)[r]{\strut{}$0.2$}}%
      \put(982,2394){\makebox(0,0)[r]{\strut{}$0.25$}}%
      \put(982,3239){\makebox(0,0)[r]{\strut{}$0.3$}}%
      \put(982,4084){\makebox(0,0)[r]{\strut{}$0.35$}}%
      \put(982,4929){\makebox(0,0)[r]{\strut{}$0.4$}}%
      \put(1114,484){\makebox(0,0){\strut{}$0$}}%
      \put(1735,484){\makebox(0,0){\strut{}$0.5$}}%
      \put(2357,484){\makebox(0,0){\strut{}$1$}}%
      \put(2978,484){\makebox(0,0){\strut{}$1.5$}}%
      \put(3600,484){\makebox(0,0){\strut{}$2$}}%
      \put(4221,484){\makebox(0,0){\strut{}$2.5$}}%
      \put(4842,484){\makebox(0,0){\strut{}$3$}}%
      \put(5464,484){\makebox(0,0){\strut{}$3.5$}}%
      \put(6085,484){\makebox(0,0){\strut{}$4$}}%
    }%
    \gplgaddtomacro\gplfronttext{%
      \csname LTb\endcsname
      \put(234,2816){\rotatebox{-270}{\makebox(0,0){\strut{}\large $r_A/\text{mm}$}}}%
      \put(3599,30){\makebox(0,0){\strut{}\large $t/\text{sec}$}}%
      \csname LTb\endcsname
      \put(5098,4701){\makebox(0,0)[r]{\strut{}$10\%$}}%
      \csname LTb\endcsname
      \put(5098,4371){\makebox(0,0)[r]{\strut{}$20\%$}}%
      \csname LTb\endcsname
      \put(5098,4041){\makebox(0,0)[r]{\strut{}$30\%$}}%
      \csname LTb\endcsname
      \put(5098,3711){\makebox(0,0)[r]{\strut{}$40\%$}}%
      \csname LTb\endcsname
      \put(5098,3381){\makebox(0,0)[r]{\strut{}$50\%$}}%
      \csname LTb\endcsname
      \put(5098,3051){\makebox(0,0)[r]{\strut{}macro}}%
      \csname LTb\endcsname
      \put(5098,2721){\makebox(0,0)[r]{\strut{}micro}}%
    }%
    \gplbacktext
    \put(0,0){\includegraphics{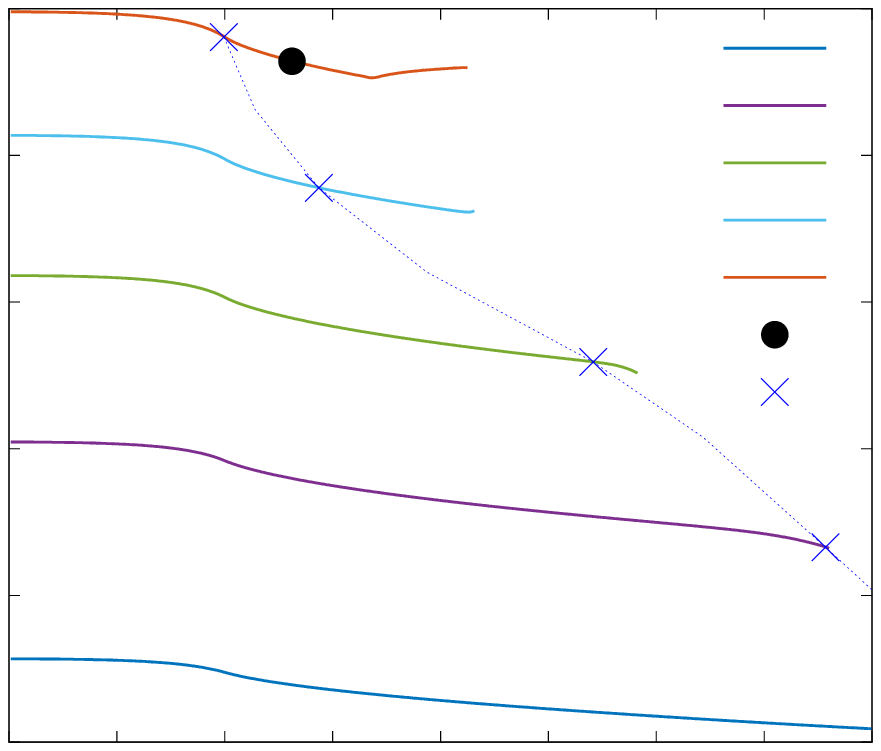}}%
    \gplfronttext
  \end{picture}%
\endgroup

%% file: time-delR.tex
\begingroup
  \makeatletter
  \providecommand\color[2][]{%
    \GenericError{(gnuplot) \space\space\space\@spaces}{%
      Package color not loaded in conjunction with
      terminal option `colourtext'%
    }{See the gnuplot documentation for explanation.%
    }{Either use 'blacktext' in gnuplot or load the package
      color.sty in LaTeX.}%
    \renewcommand\color[2][]{}%
  }%
  \providecommand\includegraphics[2][]{%
    \GenericError{(gnuplot) \space\space\space\@spaces}{%
      Package graphicx or graphics not loaded%
    }{See the gnuplot documentation for explanation.%
    }{The gnuplot epslatex terminal needs graphicx.sty or graphics.sty.}%
    \renewcommand\includegraphics[2][]{}%
  }%
  \providecommand\rotatebox[2]{#2}%
  \@ifundefined{ifGPcolor}{%
    \newif\ifGPcolor
    \GPcolortrue
  }{}%
  \@ifundefined{ifGPblacktext}{%
    \newif\ifGPblacktext
    \GPblacktexttrue
  }{}%
  \let\gplgaddtomacro\g@addto@macro
  \gdef\gplbacktext{}%
  \gdef\gplfronttext{}%
  \makeatother
  \ifGPblacktext
    \def\colorrgb#1{}%
    \def\colorgray#1{}%
  \else
    \ifGPcolor
      \def\colorrgb#1{\color[rgb]{#1}}%
      \def\colorgray#1{\color[gray]{#1}}%
      \expandafter\def\csname LTw\endcsname{\color{white}}%
      \expandafter\def\csname LTb\endcsname{\color{black}}%
      \expandafter\def\csname LTa\endcsname{\color{black}}%
      \expandafter\def\csname LT0\endcsname{\color[rgb]{1,0,0}}%
      \expandafter\def\csname LT1\endcsname{\color[rgb]{0,1,0}}%
      \expandafter\def\csname LT2\endcsname{\color[rgb]{0,0,1}}%
      \expandafter\def\csname LT3\endcsname{\color[rgb]{1,0,1}}%
      \expandafter\def\csname LT4\endcsname{\color[rgb]{0,1,1}}%
      \expandafter\def\csname LT5\endcsname{\color[rgb]{1,1,0}}%
      \expandafter\def\csname LT6\endcsname{\color[rgb]{0,0,0}}%
      \expandafter\def\csname LT7\endcsname{\color[rgb]{1,0.3,0}}%
      \expandafter\def\csname LT8\endcsname{\color[rgb]{0.5,0.5,0.5}}%
    \else
      \def\colorrgb#1{\color{black}}%
      \def\colorgray#1{\color[gray]{#1}}%
      \expandafter\def\csname LTw\endcsname{\color{white}}%
      \expandafter\def\csname LTb\endcsname{\color{black}}%
      \expandafter\def\csname LTa\endcsname{\color{black}}%
      \expandafter\def\csname LT0\endcsname{\color{black}}%
      \expandafter\def\csname LT1\endcsname{\color{black}}%
      \expandafter\def\csname LT2\endcsname{\color{black}}%
      \expandafter\def\csname LT3\endcsname{\color{black}}%
      \expandafter\def\csname LT4\endcsname{\color{black}}%
      \expandafter\def\csname LT5\endcsname{\color{black}}%
      \expandafter\def\csname LT6\endcsname{\color{black}}%
      \expandafter\def\csname LT7\endcsname{\color{black}}%
      \expandafter\def\csname LT8\endcsname{\color{black}}%
    \fi
  \fi
    \setlength{\unitlength}{0.0500bp}%
    \ifx\gptboxheight\undefined%
      \newlength{\gptboxheight}%
      \newlength{\gptboxwidth}%
      \newsavebox{\gptboxtext}%
    \fi%
    \setlength{\fboxrule}{0.5pt}%
    \setlength{\fboxsep}{1pt}%
\begin{picture}(7200.00,5040.00)%
    \gplgaddtomacro\gplbacktext{%
      \csname LTb\endcsname
      \put(982,704){\makebox(0,0)[r]{\strut{}$0$}}%
      \put(982,1232){\makebox(0,0)[r]{\strut{}$0.005$}}%
      \put(982,1760){\makebox(0,0)[r]{\strut{}$0.01$}}%
      \put(982,2288){\makebox(0,0)[r]{\strut{}$0.015$}}%
      \put(982,2817){\makebox(0,0)[r]{\strut{}$0.02$}}%
      \put(982,3345){\makebox(0,0)[r]{\strut{}$0.025$}}%
      \put(982,3873){\makebox(0,0)[r]{\strut{}$0.03$}}%
      \put(982,4401){\makebox(0,0)[r]{\strut{}$0.035$}}%
      \put(982,4929){\makebox(0,0)[r]{\strut{}$0.04$}}%
      \put(1114,484){\makebox(0,0){\strut{}$0$}}%
      \put(1735,484){\makebox(0,0){\strut{}$0.5$}}%
      \put(2357,484){\makebox(0,0){\strut{}$1$}}%
      \put(2978,484){\makebox(0,0){\strut{}$1.5$}}%
      \put(3600,484){\makebox(0,0){\strut{}$2$}}%
      \put(4221,484){\makebox(0,0){\strut{}$2.5$}}%
      \put(4842,484){\makebox(0,0){\strut{}$3$}}%
      \put(5464,484){\makebox(0,0){\strut{}$3.5$}}%
      \put(6085,484){\makebox(0,0){\strut{}$4$}}%
    }%
    \gplgaddtomacro\gplfronttext{%
      \csname LTb\endcsname
      \put(102,2816){\rotatebox{-270}{\makebox(0,0){\strut{}\large $\Delta r_A/\text{mm}$}}}%
      \put(3599,30){\makebox(0,0){\strut{}\large $t/\text{sec}$}}%
      \csname LTb\endcsname
      \put(1906,4701){\makebox(0,0)[r]{\strut{}$10\%$}}%
      \csname LTb\endcsname
      \put(1906,4371){\makebox(0,0)[r]{\strut{}$20\%$}}%
      \csname LTb\endcsname
      \put(1906,4041){\makebox(0,0)[r]{\strut{}$30\%$}}%
      \csname LTb\endcsname
      \put(1906,3711){\makebox(0,0)[r]{\strut{}$40\%$}}%
      \csname LTb\endcsname
      \put(1906,3381){\makebox(0,0)[r]{\strut{}$50\%$}}%
      \csname LTb\endcsname
      \put(1906,3051){\makebox(0,0)[r]{\strut{}macro}}%
      \csname LTb\endcsname
      \put(1906,2721){\makebox(0,0)[r]{\strut{}micro}}%
    }%
    \gplbacktext
    \put(0,0){\includegraphics{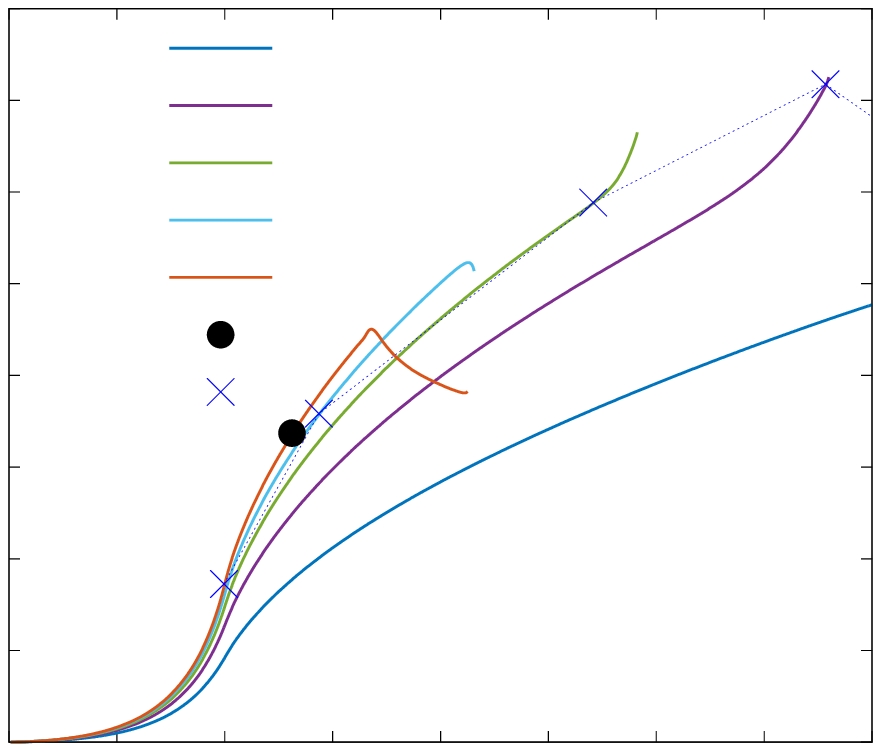}}%
    \gplfronttext
  \end{picture}%
\endgroup

%% file: time-h_v.tex
\begingroup
  \makeatletter
  \providecommand\color[2][]{%
    \GenericError{(gnuplot) \space\space\space\@spaces}{%
      Package color not loaded in conjunction with
      terminal option `colourtext'%
    }{See the gnuplot documentation for explanation.%
    }{Either use 'blacktext' in gnuplot or load the package
      color.sty in LaTeX.}%
    \renewcommand\color[2][]{}%
  }%
  \providecommand\includegraphics[2][]{%
    \GenericError{(gnuplot) \space\space\space\@spaces}{%
      Package graphicx or graphics not loaded%
    }{See the gnuplot documentation for explanation.%
    }{The gnuplot epslatex terminal needs graphicx.sty or graphics.sty.}%
    \renewcommand\includegraphics[2][]{}%
  }%
  \providecommand\rotatebox[2]{#2}%
  \@ifundefined{ifGPcolor}{%
    \newif\ifGPcolor
    \GPcolortrue
  }{}%
  \@ifundefined{ifGPblacktext}{%
    \newif\ifGPblacktext
    \GPblacktexttrue
  }{}%
  \let\gplgaddtomacro\g@addto@macro
  \gdef\gplbacktext{}%
  \gdef\gplfronttext{}%
  \makeatother
  \ifGPblacktext
    \def\colorrgb#1{}%
    \def\colorgray#1{}%
  \else
    \ifGPcolor
      \def\colorrgb#1{\color[rgb]{#1}}%
      \def\colorgray#1{\color[gray]{#1}}%
      \expandafter\def\csname LTw\endcsname{\color{white}}%
      \expandafter\def\csname LTb\endcsname{\color{black}}%
      \expandafter\def\csname LTa\endcsname{\color{black}}%
      \expandafter\def\csname LT0\endcsname{\color[rgb]{1,0,0}}%
      \expandafter\def\csname LT1\endcsname{\color[rgb]{0,1,0}}%
      \expandafter\def\csname LT2\endcsname{\color[rgb]{0,0,1}}%
      \expandafter\def\csname LT3\endcsname{\color[rgb]{1,0,1}}%
      \expandafter\def\csname LT4\endcsname{\color[rgb]{0,1,1}}%
      \expandafter\def\csname LT5\endcsname{\color[rgb]{1,1,0}}%
      \expandafter\def\csname LT6\endcsname{\color[rgb]{0,0,0}}%
      \expandafter\def\csname LT7\endcsname{\color[rgb]{1,0.3,0}}%
      \expandafter\def\csname LT8\endcsname{\color[rgb]{0.5,0.5,0.5}}%
    \else
      \def\colorrgb#1{\color{black}}%
      \def\colorgray#1{\color[gray]{#1}}%
      \expandafter\def\csname LTw\endcsname{\color{white}}%
      \expandafter\def\csname LTb\endcsname{\color{black}}%
      \expandafter\def\csname LTa\endcsname{\color{black}}%
      \expandafter\def\csname LT0\endcsname{\color{black}}%
      \expandafter\def\csname LT1\endcsname{\color{black}}%
      \expandafter\def\csname LT2\endcsname{\color{black}}%
      \expandafter\def\csname LT3\endcsname{\color{black}}%
      \expandafter\def\csname LT4\endcsname{\color{black}}%
      \expandafter\def\csname LT5\endcsname{\color{black}}%
      \expandafter\def\csname LT6\endcsname{\color{black}}%
      \expandafter\def\csname LT7\endcsname{\color{black}}%
      \expandafter\def\csname LT8\endcsname{\color{black}}%
    \fi
  \fi
    \setlength{\unitlength}{0.0500bp}%
    \ifx\gptboxheight\undefined%
      \newlength{\gptboxheight}%
      \newlength{\gptboxwidth}%
      \newsavebox{\gptboxtext}%
    \fi%
    \setlength{\fboxrule}{0.5pt}%
    \setlength{\fboxsep}{1pt}%
\begin{picture}(7200.00,5040.00)%
    \gplgaddtomacro\gplbacktext{%
      \csname LTb\endcsname
      \put(982,704){\makebox(0,0)[r]{\strut{}$0$}}%
      \put(982,1308){\makebox(0,0)[r]{\strut{}$0.02$}}%
      \put(982,1911){\makebox(0,0)[r]{\strut{}$0.04$}}%
      \put(982,2515){\makebox(0,0)[r]{\strut{}$0.06$}}%
      \put(982,3118){\makebox(0,0)[r]{\strut{}$0.08$}}%
      \put(982,3722){\makebox(0,0)[r]{\strut{}$0.1$}}%
      \put(982,4325){\makebox(0,0)[r]{\strut{}$0.12$}}%
      \put(982,4929){\makebox(0,0)[r]{\strut{}$0.14$}}%
      \put(1114,484){\makebox(0,0){\strut{}$0$}}%
      \put(1735,484){\makebox(0,0){\strut{}$0.5$}}%
      \put(2357,484){\makebox(0,0){\strut{}$1$}}%
      \put(2978,484){\makebox(0,0){\strut{}$1.5$}}%
      \put(3600,484){\makebox(0,0){\strut{}$2$}}%
      \put(4221,484){\makebox(0,0){\strut{}$2.5$}}%
      \put(4842,484){\makebox(0,0){\strut{}$3$}}%
      \put(5464,484){\makebox(0,0){\strut{}$3.5$}}%
      \put(6085,484){\makebox(0,0){\strut{}$4$}}%
    }%
    \gplgaddtomacro\gplfronttext{%
      \csname LTb\endcsname
      \put(234,2816){\rotatebox{-270}{\makebox(0,0){\strut{}\large $h_{void}/\text{(mm\textsuperscript{3}/s)}$}}}%
      \put(3599,30){\makebox(0,0){\strut{}\large $t/\text{sec}$}}%
      \csname LTb\endcsname
      \put(5098,4701){\makebox(0,0)[r]{\strut{}  $10\%$}}%
      \csname LTb\endcsname
      \put(5098,4371){\makebox(0,0)[r]{\strut{}  $20\%$}}%
      \csname LTb\endcsname
      \put(5098,4041){\makebox(0,0)[r]{\strut{}  $30\%$}}%
      \csname LTb\endcsname
      \put(5098,3711){\makebox(0,0)[r]{\strut{}  $40\%$}}%
      \csname LTb\endcsname
      \put(5098,3381){\makebox(0,0)[r]{\strut{}  $50\%$}}%
      \csname LTb\endcsname
      \put(5098,3051){\makebox(0,0)[r]{\strut{}macro}}%
      \csname LTb\endcsname
      \put(5098,2721){\makebox(0,0)[r]{\strut{}micro}}%
    }%
    \gplbacktext
    \put(0,0){\includegraphics{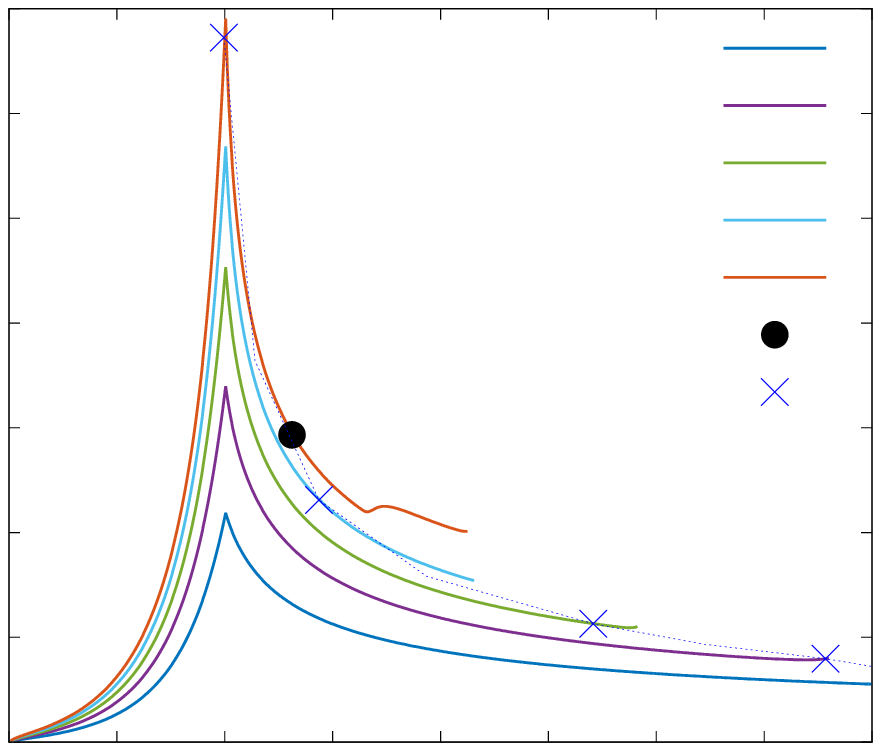}}%
    \gplfronttext
  \end{picture}%
\endgroup

%% file: time-pbar11.tex
\begingroup
  \makeatletter
  \providecommand\color[2][]{%
    \GenericError{(gnuplot) \space\space\space\@spaces}{%
      Package color not loaded in conjunction with
      terminal option `colourtext'%
    }{See the gnuplot documentation for explanation.%
    }{Either use 'blacktext' in gnuplot or load the package
      color.sty in LaTeX.}%
    \renewcommand\color[2][]{}%
  }%
  \providecommand\includegraphics[2][]{%
    \GenericError{(gnuplot) \space\space\space\@spaces}{%
      Package graphicx or graphics not loaded%
    }{See the gnuplot documentation for explanation.%
    }{The gnuplot epslatex terminal needs graphicx.sty or graphics.sty.}%
    \renewcommand\includegraphics[2][]{}%
  }%
  \providecommand\rotatebox[2]{#2}%
  \@ifundefined{ifGPcolor}{%
    \newif\ifGPcolor
    \GPcolortrue
  }{}%
  \@ifundefined{ifGPblacktext}{%
    \newif\ifGPblacktext
    \GPblacktexttrue
  }{}%
  \let\gplgaddtomacro\g@addto@macro
  \gdef\gplbacktext{}%
  \gdef\gplfronttext{}%
  \makeatother
  \ifGPblacktext
    \def\colorrgb#1{}%
    \def\colorgray#1{}%
  \else
    \ifGPcolor
      \def\colorrgb#1{\color[rgb]{#1}}%
      \def\colorgray#1{\color[gray]{#1}}%
      \expandafter\def\csname LTw\endcsname{\color{white}}%
      \expandafter\def\csname LTb\endcsname{\color{black}}%
      \expandafter\def\csname LTa\endcsname{\color{black}}%
      \expandafter\def\csname LT0\endcsname{\color[rgb]{1,0,0}}%
      \expandafter\def\csname LT1\endcsname{\color[rgb]{0,1,0}}%
      \expandafter\def\csname LT2\endcsname{\color[rgb]{0,0,1}}%
      \expandafter\def\csname LT3\endcsname{\color[rgb]{1,0,1}}%
      \expandafter\def\csname LT4\endcsname{\color[rgb]{0,1,1}}%
      \expandafter\def\csname LT5\endcsname{\color[rgb]{1,1,0}}%
      \expandafter\def\csname LT6\endcsname{\color[rgb]{0,0,0}}%
      \expandafter\def\csname LT7\endcsname{\color[rgb]{1,0.3,0}}%
      \expandafter\def\csname LT8\endcsname{\color[rgb]{0.5,0.5,0.5}}%
    \else
      \def\colorrgb#1{\color{black}}%
      \def\colorgray#1{\color[gray]{#1}}%
      \expandafter\def\csname LTw\endcsname{\color{white}}%
      \expandafter\def\csname LTb\endcsname{\color{black}}%
      \expandafter\def\csname LTa\endcsname{\color{black}}%
      \expandafter\def\csname LT0\endcsname{\color{black}}%
      \expandafter\def\csname LT1\endcsname{\color{black}}%
      \expandafter\def\csname LT2\endcsname{\color{black}}%
      \expandafter\def\csname LT3\endcsname{\color{black}}%
      \expandafter\def\csname LT4\endcsname{\color{black}}%
      \expandafter\def\csname LT5\endcsname{\color{black}}%
      \expandafter\def\csname LT6\endcsname{\color{black}}%
      \expandafter\def\csname LT7\endcsname{\color{black}}%
      \expandafter\def\csname LT8\endcsname{\color{black}}%
    \fi
  \fi
    \setlength{\unitlength}{0.0500bp}%
    \ifx\gptboxheight\undefined%
      \newlength{\gptboxheight}%
      \newlength{\gptboxwidth}%
      \newsavebox{\gptboxtext}%
    \fi%
    \setlength{\fboxrule}{0.5pt}%
    \setlength{\fboxsep}{1pt}%
\begin{picture}(7200.00,5040.00)%
    \gplgaddtomacro\gplbacktext{%
      \csname LTb\endcsname
      \put(982,704){\makebox(0,0)[r]{\strut{}$-0.4$}}%
      \put(982,1232){\makebox(0,0)[r]{\strut{}$-0.35$}}%
      \put(982,1760){\makebox(0,0)[r]{\strut{}$-0.3$}}%
      \put(982,2288){\makebox(0,0)[r]{\strut{}$-0.25$}}%
      \put(982,2816){\makebox(0,0)[r]{\strut{}$-0.2$}}%
      \put(982,3345){\makebox(0,0)[r]{\strut{}$-0.15$}}%
      \put(982,3873){\makebox(0,0)[r]{\strut{}$-0.1$}}%
      \put(982,4401){\makebox(0,0)[r]{\strut{}$-0.05$}}%
      \put(982,4929){\makebox(0,0)[r]{\strut{}$0$}}%
      \put(1114,484){\makebox(0,0){\strut{}$0$}}%
      \put(1735,484){\makebox(0,0){\strut{}$0.5$}}%
      \put(2357,484){\makebox(0,0){\strut{}$1$}}%
      \put(2978,484){\makebox(0,0){\strut{}$1.5$}}%
      \put(3600,484){\makebox(0,0){\strut{}$2$}}%
      \put(4221,484){\makebox(0,0){\strut{}$2.5$}}%
      \put(4842,484){\makebox(0,0){\strut{}$3$}}%
      \put(5464,484){\makebox(0,0){\strut{}$3.5$}}%
      \put(6085,484){\makebox(0,0){\strut{}$4$}}%
    }%
    \gplgaddtomacro\gplfronttext{%
      \csname LTb\endcsname
      \put(102,2816){\rotatebox{-270}{\makebox(0,0){\strut{}\large $\lbar{P}{}_{11}/\gamma$}}}%
      \put(3599,154){\makebox(0,0){\strut{}\large $t/\text{sec}$}}%
      \csname LTb\endcsname
      \put(5098,4701){\makebox(0,0)[r]{\strut{}  $10\%$}}%
      \csname LTb\endcsname
      \put(5098,4371){\makebox(0,0)[r]{\strut{}  $20\%$}}%
      \csname LTb\endcsname
      \put(5098,4041){\makebox(0,0)[r]{\strut{}  $30\%$}}%
      \csname LTb\endcsname
      \put(5098,3711){\makebox(0,0)[r]{\strut{}  $40\%$}}%
      \csname LTb\endcsname
      \put(5098,3381){\makebox(0,0)[r]{\strut{}  $50\%$}}%
      \csname LTb\endcsname
      \put(5098,3051){\makebox(0,0)[r]{\strut{}macro}}%
      \csname LTb\endcsname
      \put(5098,2721){\makebox(0,0)[r]{\strut{}micro}}%
    }%
    \gplbacktext
    \put(0,0){\includegraphics{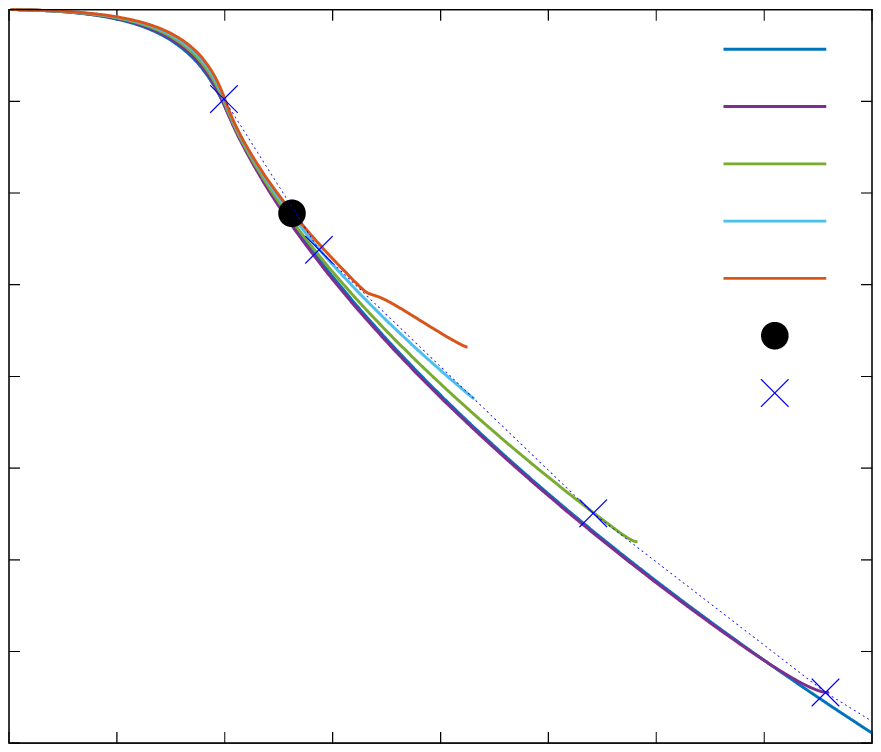}}%
    \gplfronttext
  \end{picture}%
\endgroup

%% file: mu_v-rn.tex
\begingroup
  \makeatletter
  \providecommand\color[2][]{%
    \GenericError{(gnuplot) \space\space\space\@spaces}{%
      Package color not loaded in conjunction with
      terminal option `colourtext'%
    }{See the gnuplot documentation for explanation.%
    }{Either use 'blacktext' in gnuplot or load the package
      color.sty in LaTeX.}%
    \renewcommand\color[2][]{}%
  }%
  \providecommand\includegraphics[2][]{%
    \GenericError{(gnuplot) \space\space\space\@spaces}{%
      Package graphicx or graphics not loaded%
    }{See the gnuplot documentation for explanation.%
    }{The gnuplot epslatex terminal needs graphicx.sty or graphics.sty.}%
    \renewcommand\includegraphics[2][]{}%
  }%
  \providecommand\rotatebox[2]{#2}%
  \@ifundefined{ifGPcolor}{%
    \newif\ifGPcolor
    \GPcolortrue
  }{}%
  \@ifundefined{ifGPblacktext}{%
    \newif\ifGPblacktext
    \GPblacktexttrue
  }{}%
  \let\gplgaddtomacro\g@addto@macro
  \gdef\gplbacktext{}%
  \gdef\gplfronttext{}%
  \makeatother
  \ifGPblacktext
    \def\colorrgb#1{}%
    \def\colorgray#1{}%
  \else
    \ifGPcolor
      \def\colorrgb#1{\color[rgb]{#1}}%
      \def\colorgray#1{\color[gray]{#1}}%
      \expandafter\def\csname LTw\endcsname{\color{white}}%
      \expandafter\def\csname LTb\endcsname{\color{black}}%
      \expandafter\def\csname LTa\endcsname{\color{black}}%
      \expandafter\def\csname LT0\endcsname{\color[rgb]{1,0,0}}%
      \expandafter\def\csname LT1\endcsname{\color[rgb]{0,1,0}}%
      \expandafter\def\csname LT2\endcsname{\color[rgb]{0,0,1}}%
      \expandafter\def\csname LT3\endcsname{\color[rgb]{1,0,1}}%
      \expandafter\def\csname LT4\endcsname{\color[rgb]{0,1,1}}%
      \expandafter\def\csname LT5\endcsname{\color[rgb]{1,1,0}}%
      \expandafter\def\csname LT6\endcsname{\color[rgb]{0,0,0}}%
      \expandafter\def\csname LT7\endcsname{\color[rgb]{1,0.3,0}}%
      \expandafter\def\csname LT8\endcsname{\color[rgb]{0.5,0.5,0.5}}%
    \else
      \def\colorrgb#1{\color{black}}%
      \def\colorgray#1{\color[gray]{#1}}%
      \expandafter\def\csname LTw\endcsname{\color{white}}%
      \expandafter\def\csname LTb\endcsname{\color{black}}%
      \expandafter\def\csname LTa\endcsname{\color{black}}%
      \expandafter\def\csname LT0\endcsname{\color{black}}%
      \expandafter\def\csname LT1\endcsname{\color{black}}%
      \expandafter\def\csname LT2\endcsname{\color{black}}%
      \expandafter\def\csname LT3\endcsname{\color{black}}%
      \expandafter\def\csname LT4\endcsname{\color{black}}%
      \expandafter\def\csname LT5\endcsname{\color{black}}%
      \expandafter\def\csname LT6\endcsname{\color{black}}%
      \expandafter\def\csname LT7\endcsname{\color{black}}%
      \expandafter\def\csname LT8\endcsname{\color{black}}%
    \fi
  \fi
    \setlength{\unitlength}{0.0500bp}%
    \ifx\gptboxheight\undefined%
      \newlength{\gptboxheight}%
      \newlength{\gptboxwidth}%
      \newsavebox{\gptboxtext}%
    \fi%
    \setlength{\fboxrule}{0.5pt}%
    \setlength{\fboxsep}{1pt}%
\begin{picture}(7200.00,5040.00)%
    \gplgaddtomacro\gplbacktext{%
      \csname LTb\endcsname
      \put(982,704){\makebox(0,0)[r]{\strut{} $0.15$}}%
      \put(982,1408){\makebox(0,0)[r]{\strut{}$0.2$}}%
      \put(982,2112){\makebox(0,0)[r]{\strut{}$0.25$}}%
      \put(982,2816){\makebox(0,0)[r]{\strut{}$0.3$}}%
      \put(982,3521){\makebox(0,0)[r]{\strut{}$0.35$}}%
      \put(982,4225){\makebox(0,0)[r]{\strut{}$0.4$}}%
      \put(982,4929){\makebox(0,0)[r]{\strut{}$0.45$}}%
      \put(1114,484){\makebox(0,0){\strut{}$-3.5$}}%
      \put(1824,484){\makebox(0,0){\strut{}$-3$}}%
      \put(2534,484){\makebox(0,0){\strut{}$-2.5$}}%
      \put(3244,484){\makebox(0,0){\strut{}$-2$}}%
      \put(3955,484){\makebox(0,0){\strut{}$-1.5$}}%
      \put(4665,484){\makebox(0,0){\strut{}$-1$}}%
      \put(5375,484){\makebox(0,0){\strut{}$-0.5$}}%
      \put(6085,484){\makebox(0,0){\strut{}$0$}}%
    }%
    \gplgaddtomacro\gplfronttext{%
      \csname LTb\endcsname
      \put(234,2816){\rotatebox{-270}{\makebox(0,0){\strut{}\large $r_A/\text{mm}$}}}%
      \put(3599,30){\makebox(0,0){\strut{}\large $\mu^N_{void}/\alpha$}}%
      \csname LTb\endcsname
      \put(3244,1077){\makebox(0,0)[l]{\strut{}\colorbox{white}{$10\%$}}}%
      \put(3244,2118){\makebox(0,0)[l]{\strut{}\colorbox{white}{$20\%$}}}%
      \put(3244,2917){\makebox(0,0)[l]{\strut{}\colorbox{white}{$30\%$}}}%
      \put(3244,3591){\makebox(0,0)[l]{\strut{}\colorbox{white}{$40\%$}}}%
      \put(3244,4185){\makebox(0,0)[l]{\strut{}\colorbox{white}{$50\%$}}}%
    }%
    \gplbacktext
    \put(0,0){\includegraphics{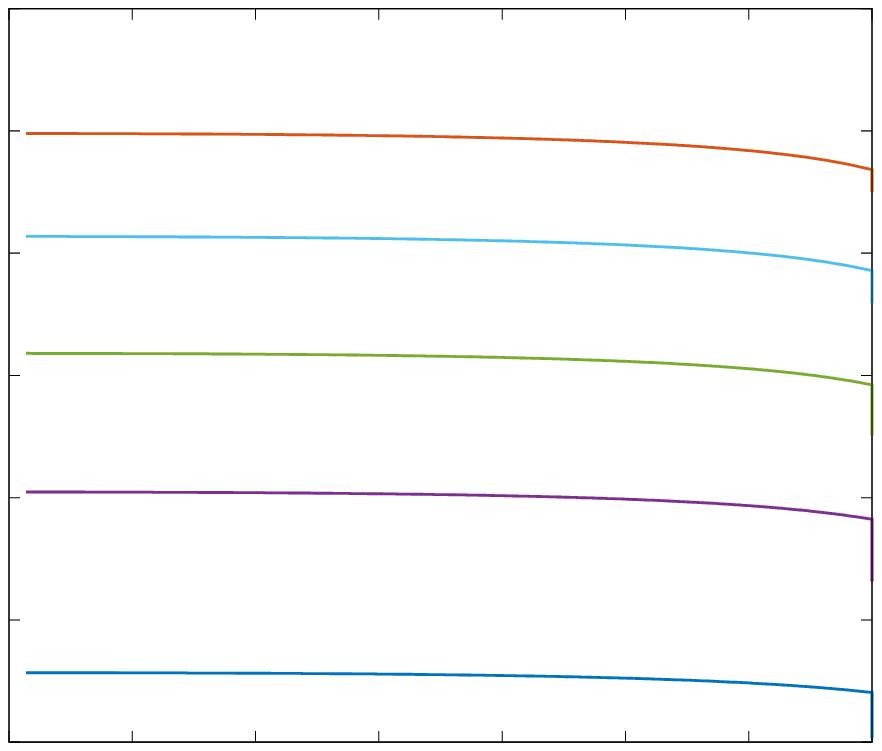}}%
    \gplfronttext
  \end{picture}%
\endgroup

%% file: mu_v-delRn.tex
\begingroup
  \makeatletter
  \providecommand\color[2][]{%
    \GenericError{(gnuplot) \space\space\space\@spaces}{%
      Package color not loaded in conjunction with
      terminal option `colourtext'%
    }{See the gnuplot documentation for explanation.%
    }{Either use 'blacktext' in gnuplot or load the package
      color.sty in LaTeX.}%
    \renewcommand\color[2][]{}%
  }%
  \providecommand\includegraphics[2][]{%
    \GenericError{(gnuplot) \space\space\space\@spaces}{%
      Package graphicx or graphics not loaded%
    }{See the gnuplot documentation for explanation.%
    }{The gnuplot epslatex terminal needs graphicx.sty or graphics.sty.}%
    \renewcommand\includegraphics[2][]{}%
  }%
  \providecommand\rotatebox[2]{#2}%
  \@ifundefined{ifGPcolor}{%
    \newif\ifGPcolor
    \GPcolortrue
  }{}%
  \@ifundefined{ifGPblacktext}{%
    \newif\ifGPblacktext
    \GPblacktexttrue
  }{}%
  \let\gplgaddtomacro\g@addto@macro
  \gdef\gplbacktext{}%
  \gdef\gplfronttext{}%
  \makeatother
  \ifGPblacktext
    \def\colorrgb#1{}%
    \def\colorgray#1{}%
  \else
    \ifGPcolor
      \def\colorrgb#1{\color[rgb]{#1}}%
      \def\colorgray#1{\color[gray]{#1}}%
      \expandafter\def\csname LTw\endcsname{\color{white}}%
      \expandafter\def\csname LTb\endcsname{\color{black}}%
      \expandafter\def\csname LTa\endcsname{\color{black}}%
      \expandafter\def\csname LT0\endcsname{\color[rgb]{1,0,0}}%
      \expandafter\def\csname LT1\endcsname{\color[rgb]{0,1,0}}%
      \expandafter\def\csname LT2\endcsname{\color[rgb]{0,0,1}}%
      \expandafter\def\csname LT3\endcsname{\color[rgb]{1,0,1}}%
      \expandafter\def\csname LT4\endcsname{\color[rgb]{0,1,1}}%
      \expandafter\def\csname LT5\endcsname{\color[rgb]{1,1,0}}%
      \expandafter\def\csname LT6\endcsname{\color[rgb]{0,0,0}}%
      \expandafter\def\csname LT7\endcsname{\color[rgb]{1,0.3,0}}%
      \expandafter\def\csname LT8\endcsname{\color[rgb]{0.5,0.5,0.5}}%
    \else
      \def\colorrgb#1{\color{black}}%
      \def\colorgray#1{\color[gray]{#1}}%
      \expandafter\def\csname LTw\endcsname{\color{white}}%
      \expandafter\def\csname LTb\endcsname{\color{black}}%
      \expandafter\def\csname LTa\endcsname{\color{black}}%
      \expandafter\def\csname LT0\endcsname{\color{black}}%
      \expandafter\def\csname LT1\endcsname{\color{black}}%
      \expandafter\def\csname LT2\endcsname{\color{black}}%
      \expandafter\def\csname LT3\endcsname{\color{black}}%
      \expandafter\def\csname LT4\endcsname{\color{black}}%
      \expandafter\def\csname LT5\endcsname{\color{black}}%
      \expandafter\def\csname LT6\endcsname{\color{black}}%
      \expandafter\def\csname LT7\endcsname{\color{black}}%
      \expandafter\def\csname LT8\endcsname{\color{black}}%
    \fi
  \fi
    \setlength{\unitlength}{0.0500bp}%
    \ifx\gptboxheight\undefined%
      \newlength{\gptboxheight}%
      \newlength{\gptboxwidth}%
      \newsavebox{\gptboxtext}%
    \fi%
    \setlength{\fboxrule}{0.5pt}%
    \setlength{\fboxsep}{1pt}%
\begin{picture}(7200.00,5040.00)%
    \gplgaddtomacro\gplbacktext{%
      \csname LTb\endcsname
      \put(982,704){\makebox(0,0)[r]{\strut{}$0$}}%
      \put(982,1232){\makebox(0,0)[r]{\strut{}$0.005$}}%
      \put(982,1760){\makebox(0,0)[r]{\strut{}$0.01$}}%
      \put(982,2288){\makebox(0,0)[r]{\strut{}$0.015$}}%
      \put(982,2817){\makebox(0,0)[r]{\strut{}$0.02$}}%
      \put(982,3345){\makebox(0,0)[r]{\strut{}$0.025$}}%
      \put(982,3873){\makebox(0,0)[r]{\strut{}$0.03$}}%
      \put(982,4401){\makebox(0,0)[r]{\strut{}$0.035$}}%
      \put(982,4929){\makebox(0,0)[r]{\strut{}$0.04$}}%
      \put(1114,484){\makebox(0,0){\strut{}$-3.5$}}%
      \put(1824,484){\makebox(0,0){\strut{}$-3$}}%
      \put(2534,484){\makebox(0,0){\strut{}$-2.5$}}%
      \put(3244,484){\makebox(0,0){\strut{}$-2$}}%
      \put(3955,484){\makebox(0,0){\strut{}$-1.5$}}%
      \put(4665,484){\makebox(0,0){\strut{}$-1$}}%
      \put(5375,484){\makebox(0,0){\strut{}$-0.5$}}%
      \put(6085,484){\makebox(0,0){\strut{}$0$}}%
    }%
    \gplgaddtomacro\gplfronttext{%
      \csname LTb\endcsname
      \put(102,2816){\rotatebox{-270}{\makebox(0,0){\strut{}\large $\Delta r_A/\text{mm}$}}}%
      \put(3599,30){\makebox(0,0){\strut{}\large $\mu^N_{void}/\alpha$}}%
    }%
    \gplgaddtomacro\gplbacktext{%
      \csname LTb\endcsname
      \put(1875,2772){\makebox(0,0)[r]{\strut{}\scriptsize$0.004$}}%
      \put(1875,3276){\makebox(0,0)[r]{\strut{}\scriptsize$0.008$}}%
      \put(1875,3779){\makebox(0,0)[r]{\strut{}\scriptsize$0.012$}}%
      \put(1875,4283){\makebox(0,0)[r]{\strut{}\scriptsize$0.016$}}%
      \put(2126,2300){\makebox(0,0){\strut{}\scriptsize$-0.6$}}%
      \put(2916,2300){\makebox(0,0){\strut{}\scriptsize$-0.4$}}%
      \put(3706,2300){\makebox(0,0){\strut{}\scriptsize$-0.2$}}%
      \put(4496,2300){\makebox(0,0){\strut{}$0$}}%
    }%
    \gplgaddtomacro\gplfronttext{%
      \csname LTb\endcsname
      \put(1334,3527){\rotatebox{-270}{\makebox(0,0){\strut{}}}}%
      \put(3311,2014){\makebox(0,0){\strut{}}}%
    }%
    \gplbacktext
    \put(0,0){\includegraphics{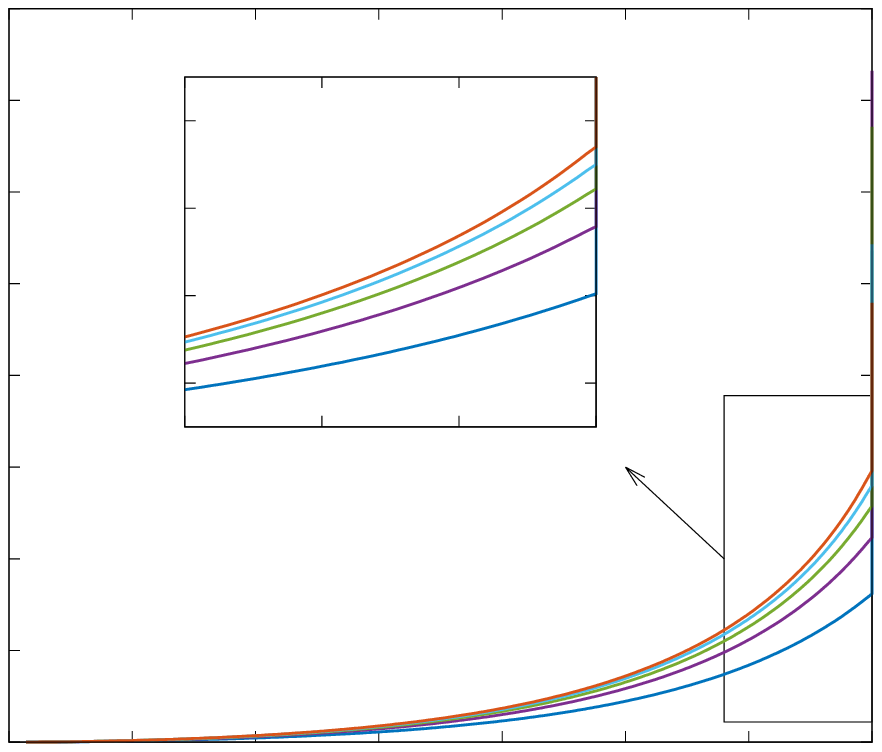}}%
    \gplfronttext
  \end{picture}%
\endgroup

%% file: mu_v-h_vn.tex
\begingroup
  \makeatletter
  \providecommand\color[2][]{%
    \GenericError{(gnuplot) \space\space\space\@spaces}{%
      Package color not loaded in conjunction with
      terminal option `colourtext'%
    }{See the gnuplot documentation for explanation.%
    }{Either use 'blacktext' in gnuplot or load the package
      color.sty in LaTeX.}%
    \renewcommand\color[2][]{}%
  }%
  \providecommand\includegraphics[2][]{%
    \GenericError{(gnuplot) \space\space\space\@spaces}{%
      Package graphicx or graphics not loaded%
    }{See the gnuplot documentation for explanation.%
    }{The gnuplot epslatex terminal needs graphicx.sty or graphics.sty.}%
    \renewcommand\includegraphics[2][]{}%
  }%
  \providecommand\rotatebox[2]{#2}%
  \@ifundefined{ifGPcolor}{%
    \newif\ifGPcolor
    \GPcolortrue
  }{}%
  \@ifundefined{ifGPblacktext}{%
    \newif\ifGPblacktext
    \GPblacktexttrue
  }{}%
  \let\gplgaddtomacro\g@addto@macro
  \gdef\gplbacktext{}%
  \gdef\gplfronttext{}%
  \makeatother
  \ifGPblacktext
    \def\colorrgb#1{}%
    \def\colorgray#1{}%
  \else
    \ifGPcolor
      \def\colorrgb#1{\color[rgb]{#1}}%
      \def\colorgray#1{\color[gray]{#1}}%
      \expandafter\def\csname LTw\endcsname{\color{white}}%
      \expandafter\def\csname LTb\endcsname{\color{black}}%
      \expandafter\def\csname LTa\endcsname{\color{black}}%
      \expandafter\def\csname LT0\endcsname{\color[rgb]{1,0,0}}%
      \expandafter\def\csname LT1\endcsname{\color[rgb]{0,1,0}}%
      \expandafter\def\csname LT2\endcsname{\color[rgb]{0,0,1}}%
      \expandafter\def\csname LT3\endcsname{\color[rgb]{1,0,1}}%
      \expandafter\def\csname LT4\endcsname{\color[rgb]{0,1,1}}%
      \expandafter\def\csname LT5\endcsname{\color[rgb]{1,1,0}}%
      \expandafter\def\csname LT6\endcsname{\color[rgb]{0,0,0}}%
      \expandafter\def\csname LT7\endcsname{\color[rgb]{1,0.3,0}}%
      \expandafter\def\csname LT8\endcsname{\color[rgb]{0.5,0.5,0.5}}%
    \else
      \def\colorrgb#1{\color{black}}%
      \def\colorgray#1{\color[gray]{#1}}%
      \expandafter\def\csname LTw\endcsname{\color{white}}%
      \expandafter\def\csname LTb\endcsname{\color{black}}%
      \expandafter\def\csname LTa\endcsname{\color{black}}%
      \expandafter\def\csname LT0\endcsname{\color{black}}%
      \expandafter\def\csname LT1\endcsname{\color{black}}%
      \expandafter\def\csname LT2\endcsname{\color{black}}%
      \expandafter\def\csname LT3\endcsname{\color{black}}%
      \expandafter\def\csname LT4\endcsname{\color{black}}%
      \expandafter\def\csname LT5\endcsname{\color{black}}%
      \expandafter\def\csname LT6\endcsname{\color{black}}%
      \expandafter\def\csname LT7\endcsname{\color{black}}%
      \expandafter\def\csname LT8\endcsname{\color{black}}%
    \fi
  \fi
    \setlength{\unitlength}{0.0500bp}%
    \ifx\gptboxheight\undefined%
      \newlength{\gptboxheight}%
      \newlength{\gptboxwidth}%
      \newsavebox{\gptboxtext}%
    \fi%
    \setlength{\fboxrule}{0.5pt}%
    \setlength{\fboxsep}{1pt}%
\begin{picture}(7200.00,5040.00)%
    \gplgaddtomacro\gplbacktext{%
      \csname LTb\endcsname
      \put(982,704){\makebox(0,0)[r]{\strut{}$0$}}%
      \put(982,1549){\makebox(0,0)[r]{\strut{}$0.05$}}%
      \put(982,2394){\makebox(0,0)[r]{\strut{}$0.1$}}%
      \put(982,3239){\makebox(0,0)[r]{\strut{}$0.15$}}%
      \put(982,4084){\makebox(0,0)[r]{\strut{}$0.2$}}%
      \put(982,4929){\makebox(0,0)[r]{\strut{}$0.25$}}%
      \put(1114,484){\makebox(0,0){\strut{}$-3.5$}}%
      \put(1824,484){\makebox(0,0){\strut{}$-3$}}%
      \put(2534,484){\makebox(0,0){\strut{}$-2.5$}}%
      \put(3244,484){\makebox(0,0){\strut{}$-2$}}%
      \put(3955,484){\makebox(0,0){\strut{}$-1.5$}}%
      \put(4665,484){\makebox(0,0){\strut{}$-1$}}%
      \put(5375,484){\makebox(0,0){\strut{}$-0.5$}}%
      \put(6085,484){\makebox(0,0){\strut{}$0$}}%
    }%
    \gplgaddtomacro\gplfronttext{%
      \csname LTb\endcsname
      \put(234,2816){\rotatebox{-270}{\makebox(0,0){\strut{}\large $h_{void}/\text{(mm\textsuperscript{3}/s)}$}}}%
      \put(3599,30){\makebox(0,0){\strut{}\large $\mu^N_{void}/\alpha$}}%
      \csname LTb\endcsname
      \put(2038,4701){\makebox(0,0)[r]{\strut{}  $10\%$}}%
      \csname LTb\endcsname
      \put(2038,4371){\makebox(0,0)[r]{\strut{}  $20\%$}}%
      \csname LTb\endcsname
      \put(2038,4041){\makebox(0,0)[r]{\strut{}  $30\%$}}%
      \csname LTb\endcsname
      \put(2038,3711){\makebox(0,0)[r]{\strut{}  $40\%$}}%
      \csname LTb\endcsname
      \put(2038,3381){\makebox(0,0)[r]{\strut{}  $50\%$}}%
    }%
    \gplbacktext
    \put(0,0){\includegraphics{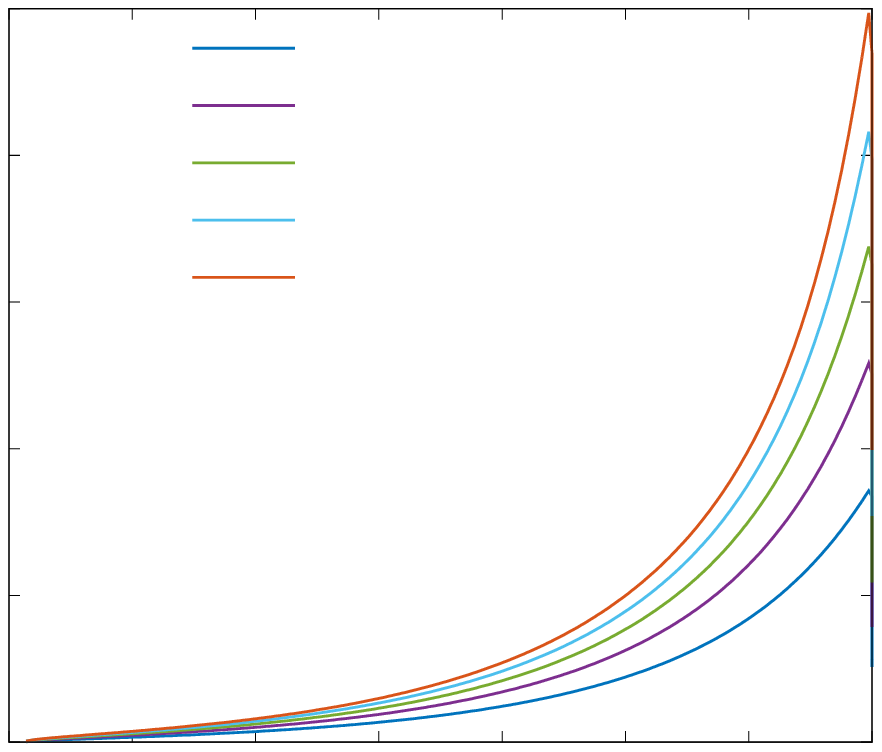}}%
    \gplfronttext
  \end{picture}%
\endgroup

%% file: mu_v-pbar11n.tex
\begingroup
  \makeatletter
  \providecommand\color[2][]{%
    \GenericError{(gnuplot) \space\space\space\@spaces}{%
      Package color not loaded in conjunction with
      terminal option `colourtext'%
    }{See the gnuplot documentation for explanation.%
    }{Either use 'blacktext' in gnuplot or load the package
      color.sty in LaTeX.}%
    \renewcommand\color[2][]{}%
  }%
  \providecommand\includegraphics[2][]{%
    \GenericError{(gnuplot) \space\space\space\@spaces}{%
      Package graphicx or graphics not loaded%
    }{See the gnuplot documentation for explanation.%
    }{The gnuplot epslatex terminal needs graphicx.sty or graphics.sty.}%
    \renewcommand\includegraphics[2][]{}%
  }%
  \providecommand\rotatebox[2]{#2}%
  \@ifundefined{ifGPcolor}{%
    \newif\ifGPcolor
    \GPcolortrue
  }{}%
  \@ifundefined{ifGPblacktext}{%
    \newif\ifGPblacktext
    \GPblacktexttrue
  }{}%
  \let\gplgaddtomacro\g@addto@macro
  \gdef\gplbacktext{}%
  \gdef\gplfronttext{}%
  \makeatother
  \ifGPblacktext
    \def\colorrgb#1{}%
    \def\colorgray#1{}%
  \else
    \ifGPcolor
      \def\colorrgb#1{\color[rgb]{#1}}%
      \def\colorgray#1{\color[gray]{#1}}%
      \expandafter\def\csname LTw\endcsname{\color{white}}%
      \expandafter\def\csname LTb\endcsname{\color{black}}%
      \expandafter\def\csname LTa\endcsname{\color{black}}%
      \expandafter\def\csname LT0\endcsname{\color[rgb]{1,0,0}}%
      \expandafter\def\csname LT1\endcsname{\color[rgb]{0,1,0}}%
      \expandafter\def\csname LT2\endcsname{\color[rgb]{0,0,1}}%
      \expandafter\def\csname LT3\endcsname{\color[rgb]{1,0,1}}%
      \expandafter\def\csname LT4\endcsname{\color[rgb]{0,1,1}}%
      \expandafter\def\csname LT5\endcsname{\color[rgb]{1,1,0}}%
      \expandafter\def\csname LT6\endcsname{\color[rgb]{0,0,0}}%
      \expandafter\def\csname LT7\endcsname{\color[rgb]{1,0.3,0}}%
      \expandafter\def\csname LT8\endcsname{\color[rgb]{0.5,0.5,0.5}}%
    \else
      \def\colorrgb#1{\color{black}}%
      \def\colorgray#1{\color[gray]{#1}}%
      \expandafter\def\csname LTw\endcsname{\color{white}}%
      \expandafter\def\csname LTb\endcsname{\color{black}}%
      \expandafter\def\csname LTa\endcsname{\color{black}}%
      \expandafter\def\csname LT0\endcsname{\color{black}}%
      \expandafter\def\csname LT1\endcsname{\color{black}}%
      \expandafter\def\csname LT2\endcsname{\color{black}}%
      \expandafter\def\csname LT3\endcsname{\color{black}}%
      \expandafter\def\csname LT4\endcsname{\color{black}}%
      \expandafter\def\csname LT5\endcsname{\color{black}}%
      \expandafter\def\csname LT6\endcsname{\color{black}}%
      \expandafter\def\csname LT7\endcsname{\color{black}}%
      \expandafter\def\csname LT8\endcsname{\color{black}}%
    \fi
  \fi
    \setlength{\unitlength}{0.0500bp}%
    \ifx\gptboxheight\undefined%
      \newlength{\gptboxheight}%
      \newlength{\gptboxwidth}%
      \newsavebox{\gptboxtext}%
    \fi%
    \setlength{\fboxrule}{0.5pt}%
    \setlength{\fboxsep}{1pt}%
\begin{picture}(7200.00,5040.00)%
    \gplgaddtomacro\gplbacktext{%
      \csname LTb\endcsname
      \put(982,1088){\makebox(0,0)[r]{\strut{}$-0.1$}}%
      \put(982,1856){\makebox(0,0)[r]{\strut{}$-0.08$}}%
      \put(982,2624){\makebox(0,0)[r]{\strut{}$-0.06$}}%
      \put(982,3393){\makebox(0,0)[r]{\strut{}$-0.04$}}%
      \put(982,4161){\makebox(0,0)[r]{\strut{}$-0.02$}}%
      \put(982,4929){\makebox(0,0)[r]{\strut{}$0$}}%
      \put(1114,484){\makebox(0,0){\strut{}$-3.5$}}%
      \put(1824,484){\makebox(0,0){\strut{}$-3$}}%
      \put(2534,484){\makebox(0,0){\strut{}$-2.5$}}%
      \put(3244,484){\makebox(0,0){\strut{}$-2$}}%
      \put(3955,484){\makebox(0,0){\strut{}$-1.5$}}%
      \put(4665,484){\makebox(0,0){\strut{}$-1$}}%
      \put(5375,484){\makebox(0,0){\strut{}$-0.5$}}%
      \put(6085,484){\makebox(0,0){\strut{}$0$}}%
    }%
    \gplgaddtomacro\gplfronttext{%
      \csname LTb\endcsname
      \put(102,2816){\rotatebox{-270}{\makebox(0,0){\strut{}\large $\lbar{P}{}_{11}/\gamma$}}}%
      \put(3599,30){\makebox(0,0){\strut{}\large $\mu^N_{void}/\alpha$}}%
      \csname LTb\endcsname
      \put(1774,2252){\makebox(0,0)[r]{\strut{}  $10\%$}}%
      \csname LTb\endcsname
      \put(1774,1922){\makebox(0,0)[r]{\strut{}  $20\%$}}%
      \csname LTb\endcsname
      \put(1774,1592){\makebox(0,0)[r]{\strut{}  $30\%$}}%
      \csname LTb\endcsname
      \put(1774,1262){\makebox(0,0)[r]{\strut{}  $40\%$}}%
      \csname LTb\endcsname
      \put(1774,932){\makebox(0,0)[r]{\strut{}  $50\%$}}%
    }%
    \gplbacktext
    \put(0,0){\includegraphics{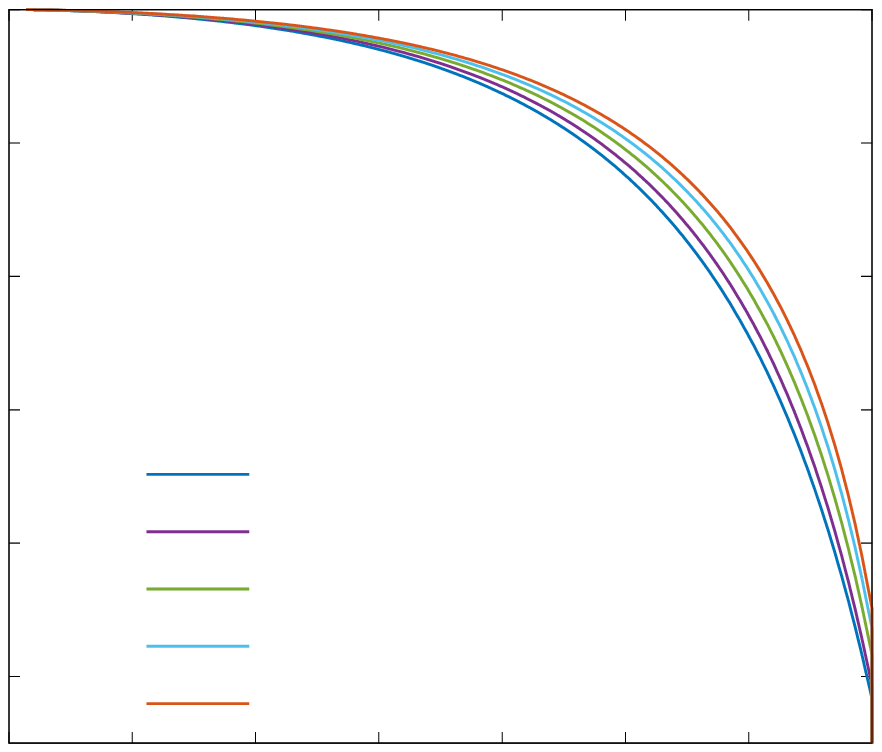}}%
    \gplfronttext
  \end{picture}%
\endgroup

%% file: f0-p11bar-vs-time.tex
\begingroup
  \makeatletter
  \providecommand\color[2][]{%
    \GenericError{(gnuplot) \space\space\space\@spaces}{%
      Package color not loaded in conjunction with
      terminal option `colourtext'%
    }{See the gnuplot documentation for explanation.%
    }{Either use 'blacktext' in gnuplot or load the package
      color.sty in LaTeX.}%
    \renewcommand\color[2][]{}%
  }%
  \providecommand\includegraphics[2][]{%
    \GenericError{(gnuplot) \space\space\space\@spaces}{%
      Package graphicx or graphics not loaded%
    }{See the gnuplot documentation for explanation.%
    }{The gnuplot epslatex terminal needs graphicx.sty or graphics.sty.}%
    \renewcommand\includegraphics[2][]{}%
  }%
  \providecommand\rotatebox[2]{#2}%
  \@ifundefined{ifGPcolor}{%
    \newif\ifGPcolor
    \GPcolortrue
  }{}%
  \@ifundefined{ifGPblacktext}{%
    \newif\ifGPblacktext
    \GPblacktexttrue
  }{}%
  \let\gplgaddtomacro\g@addto@macro
  \gdef\gplbacktext{}%
  \gdef\gplfronttext{}%
  \makeatother
  \ifGPblacktext
    \def\colorrgb#1{}%
    \def\colorgray#1{}%
  \else
    \ifGPcolor
      \def\colorrgb#1{\color[rgb]{#1}}%
      \def\colorgray#1{\color[gray]{#1}}%
      \expandafter\def\csname LTw\endcsname{\color{white}}%
      \expandafter\def\csname LTb\endcsname{\color{black}}%
      \expandafter\def\csname LTa\endcsname{\color{black}}%
      \expandafter\def\csname LT0\endcsname{\color[rgb]{1,0,0}}%
      \expandafter\def\csname LT1\endcsname{\color[rgb]{0,1,0}}%
      \expandafter\def\csname LT2\endcsname{\color[rgb]{0,0,1}}%
      \expandafter\def\csname LT3\endcsname{\color[rgb]{1,0,1}}%
      \expandafter\def\csname LT4\endcsname{\color[rgb]{0,1,1}}%
      \expandafter\def\csname LT5\endcsname{\color[rgb]{1,1,0}}%
      \expandafter\def\csname LT6\endcsname{\color[rgb]{0,0,0}}%
      \expandafter\def\csname LT7\endcsname{\color[rgb]{1,0.3,0}}%
      \expandafter\def\csname LT8\endcsname{\color[rgb]{0.5,0.5,0.5}}%
    \else
      \def\colorrgb#1{\color{black}}%
      \def\colorgray#1{\color[gray]{#1}}%
      \expandafter\def\csname LTw\endcsname{\color{white}}%
      \expandafter\def\csname LTb\endcsname{\color{black}}%
      \expandafter\def\csname LTa\endcsname{\color{black}}%
      \expandafter\def\csname LT0\endcsname{\color{black}}%
      \expandafter\def\csname LT1\endcsname{\color{black}}%
      \expandafter\def\csname LT2\endcsname{\color{black}}%
      \expandafter\def\csname LT3\endcsname{\color{black}}%
      \expandafter\def\csname LT4\endcsname{\color{black}}%
      \expandafter\def\csname LT5\endcsname{\color{black}}%
      \expandafter\def\csname LT6\endcsname{\color{black}}%
      \expandafter\def\csname LT7\endcsname{\color{black}}%
      \expandafter\def\csname LT8\endcsname{\color{black}}%
    \fi
  \fi
    \setlength{\unitlength}{0.0500bp}%
    \ifx\gptboxheight\undefined%
      \newlength{\gptboxheight}%
      \newlength{\gptboxwidth}%
      \newsavebox{\gptboxtext}%
    \fi%
    \setlength{\fboxrule}{0.5pt}%
    \setlength{\fboxsep}{1pt}%
\begin{picture}(7200.00,5040.00)%
    \gplgaddtomacro\gplbacktext{%
      \csname LTb\endcsname
      \put(982,704){\makebox(0,0)[r]{\strut{}$0.5$}}%
      \put(982,1173){\makebox(0,0)[r]{\strut{}$1.0$}}%
      \put(982,1643){\makebox(0,0)[r]{\strut{}$1.5$}}%
      \put(982,2112){\makebox(0,0)[r]{\strut{}$2.0$}}%
      \put(982,2582){\makebox(0,0)[r]{\strut{}$2.5$}}%
      \put(982,3051){\makebox(0,0)[r]{\strut{}$3.0$}}%
      \put(982,3521){\makebox(0,0)[r]{\strut{}$3.5$}}%
      \put(982,3990){\makebox(0,0)[r]{\strut{}$4.0$}}%
      \put(982,4460){\makebox(0,0)[r]{\strut{}$4.5$}}%
      \put(982,4929){\makebox(0,0)[r]{\strut{}$5.0$}}%
      \put(1114,484){\makebox(0,0){\strut{}$10$}}%
      \put(1735,484){\makebox(0,0){\strut{}$15$}}%
      \put(2357,484){\makebox(0,0){\strut{}$20$}}%
      \put(2978,484){\makebox(0,0){\strut{}$25$}}%
      \put(3600,484){\makebox(0,0){\strut{}$30$}}%
      \put(4221,484){\makebox(0,0){\strut{}$35$}}%
      \put(4842,484){\makebox(0,0){\strut{}$40$}}%
      \put(5464,484){\makebox(0,0){\strut{}$45$}}%
      \put(6085,484){\makebox(0,0){\strut{}$50$}}%
    }%
    \gplgaddtomacro\gplfronttext{%
      \csname LTb\endcsname
      \put(366,2816){\rotatebox{-270}{\makebox(0,0){\strut{} \large $t_{crit}/\text{sec}$}}}%
      \put(3599,154){\makebox(0,0){\strut{}\large $f_0/\%$}}%
      \csname LTb\endcsname
      \put(5098,4695){\makebox(0,0)[r]{\strut{}\footnotesize $M=1.0\cdot 10^{-4}\,\text{mm}^2/\text{(Ns)}$}}%
      \csname LTb\endcsname
      \put(5098,4354){\makebox(0,0)[r]{\strut{}\footnotesize$M=2.5\cdot 10^{-4}\,\text{mm}^2/\text{(Ns)}$}}%
      \csname LTb\endcsname
      \put(5098,4013){\makebox(0,0)[r]{\strut{}\footnotesize$M=5.0\cdot 10^{-4}\,\text{mm}^2/\text{(Ns)}$}}%
      \csname LTb\endcsname
      \put(5098,3672){\makebox(0,0)[r]{\strut{}\footnotesize$M=7.5\cdot 10^{-4}\,\text{mm}^2/\text{(Ns)}$}}%
      \csname LTb\endcsname
      \put(5098,3331){\makebox(0,0)[r]{\strut{}\footnotesize macro ($\Bk\to\Bzero$)}}%
    }%
    \gplbacktext
    \put(0,0){\includegraphics{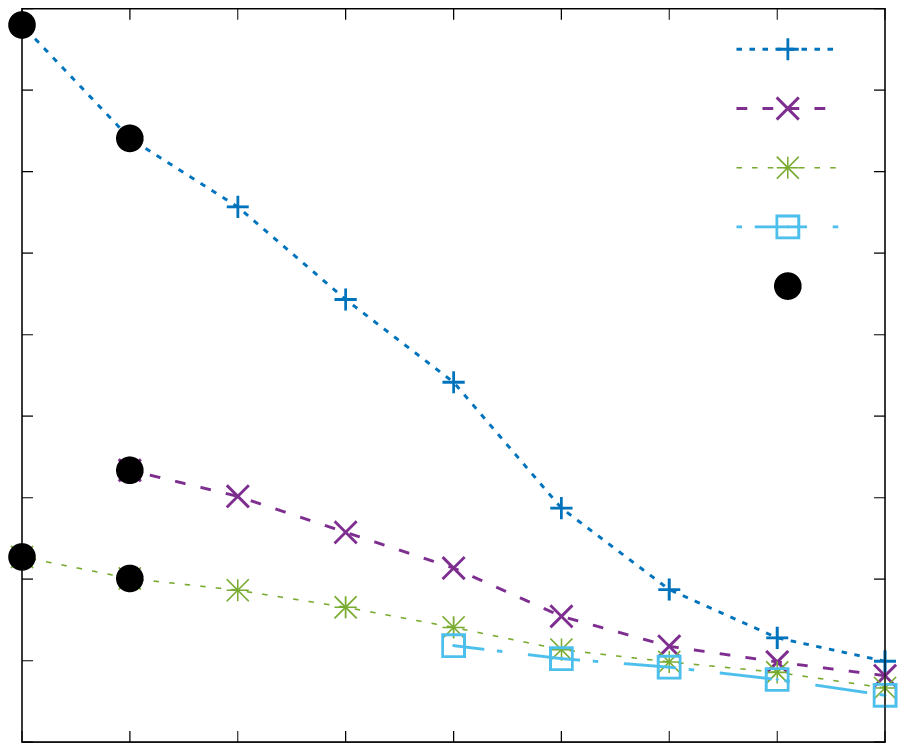}}%
    \gplfronttext
  \end{picture}%
\endgroup

%% file: f0-p11bar-vs-M.tex
\begingroup
  \makeatletter
  \providecommand\color[2][]{%
    \GenericError{(gnuplot) \space\space\space\@spaces}{%
      Package color not loaded in conjunction with
      terminal option `colourtext'%
    }{See the gnuplot documentation for explanation.%
    }{Either use 'blacktext' in gnuplot or load the package
      color.sty in LaTeX.}%
    \renewcommand\color[2][]{}%
  }%
  \providecommand\includegraphics[2][]{%
    \GenericError{(gnuplot) \space\space\space\@spaces}{%
      Package graphicx or graphics not loaded%
    }{See the gnuplot documentation for explanation.%
    }{The gnuplot epslatex terminal needs graphicx.sty or graphics.sty.}%
    \renewcommand\includegraphics[2][]{}%
  }%
  \providecommand\rotatebox[2]{#2}%
  \@ifundefined{ifGPcolor}{%
    \newif\ifGPcolor
    \GPcolortrue
  }{}%
  \@ifundefined{ifGPblacktext}{%
    \newif\ifGPblacktext
    \GPblacktexttrue
  }{}%
  \let\gplgaddtomacro\g@addto@macro
  \gdef\gplbacktext{}%
  \gdef\gplfronttext{}%
  \makeatother
  \ifGPblacktext
    \def\colorrgb#1{}%
    \def\colorgray#1{}%
  \else
    \ifGPcolor
      \def\colorrgb#1{\color[rgb]{#1}}%
      \def\colorgray#1{\color[gray]{#1}}%
      \expandafter\def\csname LTw\endcsname{\color{white}}%
      \expandafter\def\csname LTb\endcsname{\color{black}}%
      \expandafter\def\csname LTa\endcsname{\color{black}}%
      \expandafter\def\csname LT0\endcsname{\color[rgb]{1,0,0}}%
      \expandafter\def\csname LT1\endcsname{\color[rgb]{0,1,0}}%
      \expandafter\def\csname LT2\endcsname{\color[rgb]{0,0,1}}%
      \expandafter\def\csname LT3\endcsname{\color[rgb]{1,0,1}}%
      \expandafter\def\csname LT4\endcsname{\color[rgb]{0,1,1}}%
      \expandafter\def\csname LT5\endcsname{\color[rgb]{1,1,0}}%
      \expandafter\def\csname LT6\endcsname{\color[rgb]{0,0,0}}%
      \expandafter\def\csname LT7\endcsname{\color[rgb]{1,0.3,0}}%
      \expandafter\def\csname LT8\endcsname{\color[rgb]{0.5,0.5,0.5}}%
    \else
      \def\colorrgb#1{\color{black}}%
      \def\colorgray#1{\color[gray]{#1}}%
      \expandafter\def\csname LTw\endcsname{\color{white}}%
      \expandafter\def\csname LTb\endcsname{\color{black}}%
      \expandafter\def\csname LTa\endcsname{\color{black}}%
      \expandafter\def\csname LT0\endcsname{\color{black}}%
      \expandafter\def\csname LT1\endcsname{\color{black}}%
      \expandafter\def\csname LT2\endcsname{\color{black}}%
      \expandafter\def\csname LT3\endcsname{\color{black}}%
      \expandafter\def\csname LT4\endcsname{\color{black}}%
      \expandafter\def\csname LT5\endcsname{\color{black}}%
      \expandafter\def\csname LT6\endcsname{\color{black}}%
      \expandafter\def\csname LT7\endcsname{\color{black}}%
      \expandafter\def\csname LT8\endcsname{\color{black}}%
    \fi
  \fi
    \setlength{\unitlength}{0.0500bp}%
    \ifx\gptboxheight\undefined%
      \newlength{\gptboxheight}%
      \newlength{\gptboxwidth}%
      \newsavebox{\gptboxtext}%
    \fi%
    \setlength{\fboxrule}{0.5pt}%
    \setlength{\fboxsep}{1pt}%
\begin{picture}(7200.00,5040.00)%
    \gplgaddtomacro\gplbacktext{%
      \csname LTb\endcsname
      \put(982,704){\makebox(0,0)[r]{\strut{}$-0.50$}}%
      \put(982,1127){\makebox(0,0)[r]{\strut{}$-0.45$}}%
      \put(982,1549){\makebox(0,0)[r]{\strut{}$-0.40$}}%
      \put(982,1971){\makebox(0,0)[r]{\strut{}$-0.35$}}%
      \put(982,2394){\makebox(0,0)[r]{\strut{}$-0.30$}}%
      \put(982,2816){\makebox(0,0)[r]{\strut{}$-0.25$}}%
      \put(982,3239){\makebox(0,0)[r]{\strut{}$-0.20$}}%
      \put(982,3661){\makebox(0,0)[r]{\strut{}$-0.15$}}%
      \put(982,4084){\makebox(0,0)[r]{\strut{}$-0.10$}}%
      \put(982,4507){\makebox(0,0)[r]{\strut{}$-0.05$}}%
      \put(982,4929){\makebox(0,0)[r]{\strut{}$0$}}%
      \put(1114,484){\makebox(0,0){\strut{}$10$}}%
      \put(1735,484){\makebox(0,0){\strut{}$15$}}%
      \put(2357,484){\makebox(0,0){\strut{}$20$}}%
      \put(2978,484){\makebox(0,0){\strut{}$25$}}%
      \put(3600,484){\makebox(0,0){\strut{}$30$}}%
      \put(4221,484){\makebox(0,0){\strut{}$35$}}%
      \put(4842,484){\makebox(0,0){\strut{}$40$}}%
      \put(5464,484){\makebox(0,0){\strut{}$45$}}%
      \put(6085,484){\makebox(0,0){\strut{}$50$}}%
    }%
    \gplgaddtomacro\gplfronttext{%
      \csname LTb\endcsname
      \put(102,2816){\rotatebox{-270}{\makebox(0,0){\strut{}\large $\lbar{P}_{11}/\gamma$}}}%
      \put(3599,154){\makebox(0,0){\strut{}\large $f_0/\%$}}%
      \csname LTb\endcsname
      \put(5098,2301){\makebox(0,0)[r]{\strut{}\scriptsize$M=1.0\cdot 10^{-4}\,\text{mm}^2/\text{(Ns)}$}}%
      \csname LTb\endcsname
      \put(5098,1960){\makebox(0,0)[r]{\strut{}\scriptsize$M=2.5\cdot 10^{-4}\,\text{mm}^2/\text{(Ns)}$}}%
      \csname LTb\endcsname
      \put(5098,1619){\makebox(0,0)[r]{\strut{}\scriptsize$M=5.0\cdot 10^{-4}\,\text{mm}^2/\text{(Ns)}$}}%
      \csname LTb\endcsname
      \put(5098,1278){\makebox(0,0)[r]{\strut{}\scriptsize$M=7.5\cdot 10^{-4}\,\text{mm}^2/\text{(Ns)}$}}%
      \csname LTb\endcsname
      \put(5098,937){\makebox(0,0)[r]{\strut{}\footnotesize macro ($\Bk\to\Bzero$)}}%
    }%
    \gplgaddtomacro\gplbacktext{%
      \csname LTb\endcsname
      \put(1990,3498){\makebox(0,0)[r]{\strut{}\scriptsize$-0.46\,\,$}}%
      \put(1990,3780){\makebox(0,0)[r]{\strut{}\scriptsize$-0.44\,\,$}}%
      \put(1990,4062){\makebox(0,0)[r]{\strut{}\scriptsize$-0.42\,\,$}}%
      \put(1990,4344){\makebox(0,0)[r]{\strut{}\scriptsize$-0.40\,\,$}}%
      \put(1990,4626){\makebox(0,0)[r]{\strut{}\scriptsize$-0.38\,\,$}}%
      \put(2122,3207){\makebox(0,0){\strut{}\scriptsize $10$}}%
      \put(2952,3207){\makebox(0,0){\strut{}\scriptsize $15$}}%
      \put(3781,3207){\makebox(0,0){\strut{}\scriptsize $20$}}%
    }%
    \gplgaddtomacro\gplfronttext{%
      \csname LTb\endcsname
      \put(1330,4132){\rotatebox{-270}{\makebox(0,0){\strut{}}}}%
      \put(2951,2921){\makebox(0,0){\strut{}}}%
    }%
    \gplbacktext
    \put(0,0){\includegraphics{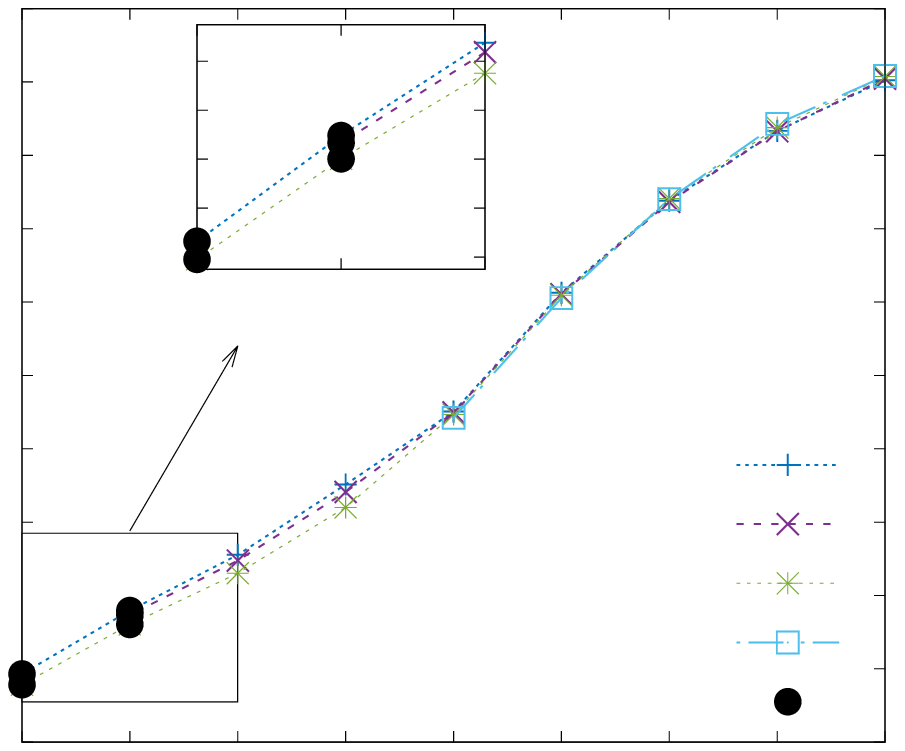}}%
    \gplfronttext
  \end{picture}%
\endgroup

%% file: f0-vs-time-alpha.tex
\begingroup
  \makeatletter
  \providecommand\color[2][]{%
    \GenericError{(gnuplot) \space\space\space\@spaces}{%
      Package color not loaded in conjunction with
      terminal option `colourtext'%
    }{See the gnuplot documentation for explanation.%
    }{Either use 'blacktext' in gnuplot or load the package
      color.sty in LaTeX.}%
    \renewcommand\color[2][]{}%
  }%
  \providecommand\includegraphics[2][]{%
    \GenericError{(gnuplot) \space\space\space\@spaces}{%
      Package graphicx or graphics not loaded%
    }{See the gnuplot documentation for explanation.%
    }{The gnuplot epslatex terminal needs graphicx.sty or graphics.sty.}%
    \renewcommand\includegraphics[2][]{}%
  }%
  \providecommand\rotatebox[2]{#2}%
  \@ifundefined{ifGPcolor}{%
    \newif\ifGPcolor
    \GPcolortrue
  }{}%
  \@ifundefined{ifGPblacktext}{%
    \newif\ifGPblacktext
    \GPblacktexttrue
  }{}%
  \let\gplgaddtomacro\g@addto@macro
  \gdef\gplbacktext{}%
  \gdef\gplfronttext{}%
  \makeatother
  \ifGPblacktext
    \def\colorrgb#1{}%
    \def\colorgray#1{}%
  \else
    \ifGPcolor
      \def\colorrgb#1{\color[rgb]{#1}}%
      \def\colorgray#1{\color[gray]{#1}}%
      \expandafter\def\csname LTw\endcsname{\color{white}}%
      \expandafter\def\csname LTb\endcsname{\color{black}}%
      \expandafter\def\csname LTa\endcsname{\color{black}}%
      \expandafter\def\csname LT0\endcsname{\color[rgb]{1,0,0}}%
      \expandafter\def\csname LT1\endcsname{\color[rgb]{0,1,0}}%
      \expandafter\def\csname LT2\endcsname{\color[rgb]{0,0,1}}%
      \expandafter\def\csname LT3\endcsname{\color[rgb]{1,0,1}}%
      \expandafter\def\csname LT4\endcsname{\color[rgb]{0,1,1}}%
      \expandafter\def\csname LT5\endcsname{\color[rgb]{1,1,0}}%
      \expandafter\def\csname LT6\endcsname{\color[rgb]{0,0,0}}%
      \expandafter\def\csname LT7\endcsname{\color[rgb]{1,0.3,0}}%
      \expandafter\def\csname LT8\endcsname{\color[rgb]{0.5,0.5,0.5}}%
    \else
      \def\colorrgb#1{\color{black}}%
      \def\colorgray#1{\color[gray]{#1}}%
      \expandafter\def\csname LTw\endcsname{\color{white}}%
      \expandafter\def\csname LTb\endcsname{\color{black}}%
      \expandafter\def\csname LTa\endcsname{\color{black}}%
      \expandafter\def\csname LT0\endcsname{\color{black}}%
      \expandafter\def\csname LT1\endcsname{\color{black}}%
      \expandafter\def\csname LT2\endcsname{\color{black}}%
      \expandafter\def\csname LT3\endcsname{\color{black}}%
      \expandafter\def\csname LT4\endcsname{\color{black}}%
      \expandafter\def\csname LT5\endcsname{\color{black}}%
      \expandafter\def\csname LT6\endcsname{\color{black}}%
      \expandafter\def\csname LT7\endcsname{\color{black}}%
      \expandafter\def\csname LT8\endcsname{\color{black}}%
    \fi
  \fi
    \setlength{\unitlength}{0.0500bp}%
    \ifx\gptboxheight\undefined%
      \newlength{\gptboxheight}%
      \newlength{\gptboxwidth}%
      \newsavebox{\gptboxtext}%
    \fi%
    \setlength{\fboxrule}{0.5pt}%
    \setlength{\fboxsep}{1pt}%
\begin{picture}(7200.00,5040.00)%
    \gplgaddtomacro\gplbacktext{%
      \csname LTb\endcsname
      \put(982,704){\makebox(0,0)[r]{\strut{}$0.8$}}%
      \put(982,1173){\makebox(0,0)[r]{\strut{}$1.0$}}%
      \put(982,1643){\makebox(0,0)[r]{\strut{}$1.2$}}%
      \put(982,2112){\makebox(0,0)[r]{\strut{}$1.4$}}%
      \put(982,2582){\makebox(0,0)[r]{\strut{}$1.6$}}%
      \put(982,3051){\makebox(0,0)[r]{\strut{}$1.8$}}%
      \put(982,3521){\makebox(0,0)[r]{\strut{}$2.0$}}%
      \put(982,3990){\makebox(0,0)[r]{\strut{}$2.2$}}%
      \put(982,4460){\makebox(0,0)[r]{\strut{}$2.4$}}%
      \put(982,4929){\makebox(0,0)[r]{\strut{}$2.6$}}%
      \put(1114,484){\makebox(0,0){\strut{}$10$}}%
      \put(1735,484){\makebox(0,0){\strut{}$15$}}%
      \put(2357,484){\makebox(0,0){\strut{}$20$}}%
      \put(2978,484){\makebox(0,0){\strut{}$25$}}%
      \put(3600,484){\makebox(0,0){\strut{}$30$}}%
      \put(4221,484){\makebox(0,0){\strut{}$35$}}%
      \put(4842,484){\makebox(0,0){\strut{}$40$}}%
      \put(5464,484){\makebox(0,0){\strut{}$45$}}%
      \put(6085,484){\makebox(0,0){\strut{}$50$}}%
    }%
    \gplgaddtomacro\gplfronttext{%
      \csname LTb\endcsname
      \put(366,2816){\rotatebox{-270}{\makebox(0,0){\strut{}\large $t_{crit}/\text{sec}$}}}%
      \put(3599,154){\makebox(0,0){\strut{}\large $f_0/\%$}}%
      \csname LTb\endcsname
      \put(5098,4695){\makebox(0,0)[r]{\strut{}\footnotesize$\alpha=20\,\text{N/mm\textsuperscript{2}}$}}%
      \csname LTb\endcsname
      \put(5098,4354){\makebox(0,0)[r]{\strut{}\footnotesize$\alpha=40\,\text{N/mm\textsuperscript{2}}$}}%
      \csname LTb\endcsname
      \put(5098,4013){\makebox(0,0)[r]{\strut{}\footnotesize$\alpha=60\,\text{N/mm\textsuperscript{2}}$}}%
      \csname LTb\endcsname
      \put(5098,3672){\makebox(0,0)[r]{\strut{}\footnotesize$\alpha=80\,\text{N/mm\textsuperscript{2}}$}}%
      \csname LTb\endcsname
      \put(5098,3331){\makebox(0,0)[r]{\strut{}\footnotesize macro ($\Bk\to\Bzero$)}}%
    }%
    \gplbacktext
    \put(0,0){\includegraphics{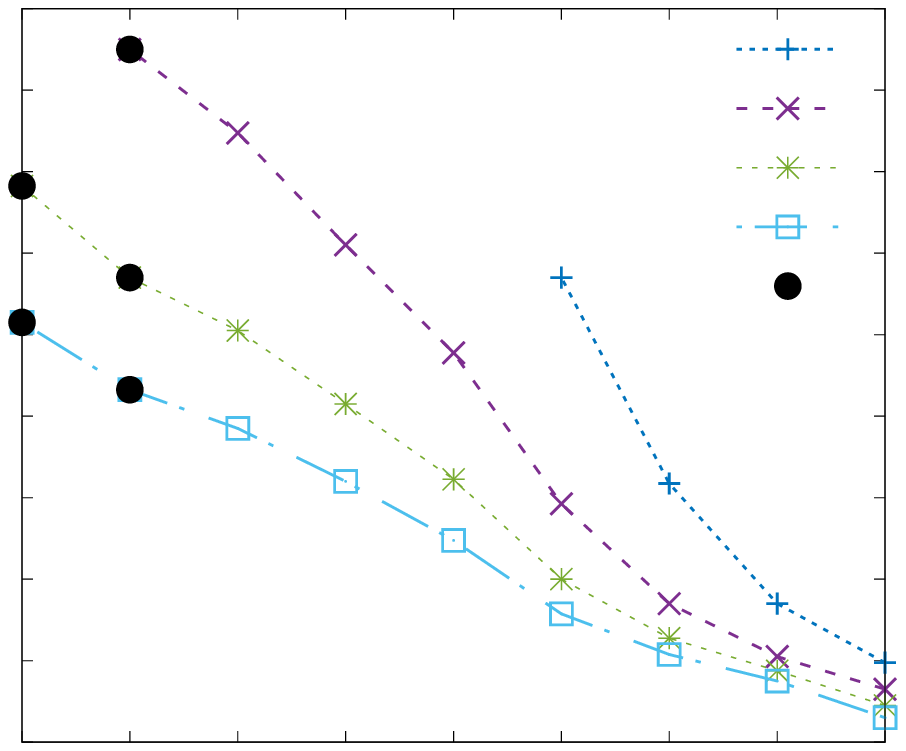}}%
    \gplfronttext
  \end{picture}%
\endgroup

%% file: f0-vs-p11bar-alpha.tex
\begingroup
  \makeatletter
  \providecommand\color[2][]{%
    \GenericError{(gnuplot) \space\space\space\@spaces}{%
      Package color not loaded in conjunction with
      terminal option `colourtext'%
    }{See the gnuplot documentation for explanation.%
    }{Either use 'blacktext' in gnuplot or load the package
      color.sty in LaTeX.}%
    \renewcommand\color[2][]{}%
  }%
  \providecommand\includegraphics[2][]{%
    \GenericError{(gnuplot) \space\space\space\@spaces}{%
      Package graphicx or graphics not loaded%
    }{See the gnuplot documentation for explanation.%
    }{The gnuplot epslatex terminal needs graphicx.sty or graphics.sty.}%
    \renewcommand\includegraphics[2][]{}%
  }%
  \providecommand\rotatebox[2]{#2}%
  \@ifundefined{ifGPcolor}{%
    \newif\ifGPcolor
    \GPcolortrue
  }{}%
  \@ifundefined{ifGPblacktext}{%
    \newif\ifGPblacktext
    \GPblacktexttrue
  }{}%
  \let\gplgaddtomacro\g@addto@macro
  \gdef\gplbacktext{}%
  \gdef\gplfronttext{}%
  \makeatother
  \ifGPblacktext
    \def\colorrgb#1{}%
    \def\colorgray#1{}%
  \else
    \ifGPcolor
      \def\colorrgb#1{\color[rgb]{#1}}%
      \def\colorgray#1{\color[gray]{#1}}%
      \expandafter\def\csname LTw\endcsname{\color{white}}%
      \expandafter\def\csname LTb\endcsname{\color{black}}%
      \expandafter\def\csname LTa\endcsname{\color{black}}%
      \expandafter\def\csname LT0\endcsname{\color[rgb]{1,0,0}}%
      \expandafter\def\csname LT1\endcsname{\color[rgb]{0,1,0}}%
      \expandafter\def\csname LT2\endcsname{\color[rgb]{0,0,1}}%
      \expandafter\def\csname LT3\endcsname{\color[rgb]{1,0,1}}%
      \expandafter\def\csname LT4\endcsname{\color[rgb]{0,1,1}}%
      \expandafter\def\csname LT5\endcsname{\color[rgb]{1,1,0}}%
      \expandafter\def\csname LT6\endcsname{\color[rgb]{0,0,0}}%
      \expandafter\def\csname LT7\endcsname{\color[rgb]{1,0.3,0}}%
      \expandafter\def\csname LT8\endcsname{\color[rgb]{0.5,0.5,0.5}}%
    \else
      \def\colorrgb#1{\color{black}}%
      \def\colorgray#1{\color[gray]{#1}}%
      \expandafter\def\csname LTw\endcsname{\color{white}}%
      \expandafter\def\csname LTb\endcsname{\color{black}}%
      \expandafter\def\csname LTa\endcsname{\color{black}}%
      \expandafter\def\csname LT0\endcsname{\color{black}}%
      \expandafter\def\csname LT1\endcsname{\color{black}}%
      \expandafter\def\csname LT2\endcsname{\color{black}}%
      \expandafter\def\csname LT3\endcsname{\color{black}}%
      \expandafter\def\csname LT4\endcsname{\color{black}}%
      \expandafter\def\csname LT5\endcsname{\color{black}}%
      \expandafter\def\csname LT6\endcsname{\color{black}}%
      \expandafter\def\csname LT7\endcsname{\color{black}}%
      \expandafter\def\csname LT8\endcsname{\color{black}}%
    \fi
  \fi
    \setlength{\unitlength}{0.0500bp}%
    \ifx\gptboxheight\undefined%
      \newlength{\gptboxheight}%
      \newlength{\gptboxwidth}%
      \newsavebox{\gptboxtext}%
    \fi%
    \setlength{\fboxrule}{0.5pt}%
    \setlength{\fboxsep}{1pt}%
\begin{picture}(7200.00,5040.00)%
    \gplgaddtomacro\gplbacktext{%
      \csname LTb\endcsname
      \put(982,704){\makebox(0,0)[r]{\strut{}$-0.50$}}%
      \put(982,1127){\makebox(0,0)[r]{\strut{}$-0.45$}}%
      \put(982,1549){\makebox(0,0)[r]{\strut{}$-0.40$}}%
      \put(982,1971){\makebox(0,0)[r]{\strut{}$-0.35$}}%
      \put(982,2394){\makebox(0,0)[r]{\strut{}$-0.30$}}%
      \put(982,2816){\makebox(0,0)[r]{\strut{}$-0.25$}}%
      \put(982,3239){\makebox(0,0)[r]{\strut{}$-0.20$}}%
      \put(982,3661){\makebox(0,0)[r]{\strut{}$-0.15$}}%
      \put(982,4084){\makebox(0,0)[r]{\strut{}$-0.10$}}%
      \put(982,4507){\makebox(0,0)[r]{\strut{}$-0.05$}}%
      \put(982,4929){\makebox(0,0)[r]{\strut{}$0$}}%
      \put(1114,484){\makebox(0,0){\strut{}$10$}}%
      \put(1735,484){\makebox(0,0){\strut{}$15$}}%
      \put(2357,484){\makebox(0,0){\strut{}$20$}}%
      \put(2978,484){\makebox(0,0){\strut{}$25$}}%
      \put(3600,484){\makebox(0,0){\strut{}$30$}}%
      \put(4221,484){\makebox(0,0){\strut{}$35$}}%
      \put(4842,484){\makebox(0,0){\strut{}$40$}}%
      \put(5464,484){\makebox(0,0){\strut{}$45$}}%
      \put(6085,484){\makebox(0,0){\strut{}$50$}}%
    }%
    \gplgaddtomacro\gplfronttext{%
      \csname LTb\endcsname
      \put(102,2816){\rotatebox{-270}{\makebox(0,0){\strut{}\large $\lbar{P}_{11}/\gamma$}}}%
      \put(3599,154){\makebox(0,0){\strut{}\large $f_0/\%$}}%
      \csname LTb\endcsname
      \put(5098,2301){\makebox(0,0)[r]{\strut{}\footnotesize$\alpha=20\,\text{N/mm\textsuperscript{2}}$}}%
      \csname LTb\endcsname
      \put(5098,1960){\makebox(0,0)[r]{\strut{}\footnotesize$\alpha=40\,\text{N/mm\textsuperscript{2}}$}}%
      \csname LTb\endcsname
      \put(5098,1619){\makebox(0,0)[r]{\strut{}\footnotesize$\alpha=60\,\text{N/mm\textsuperscript{2}}$}}%
      \csname LTb\endcsname
      \put(5098,1278){\makebox(0,0)[r]{\strut{}\footnotesize$\alpha=80\,\text{N/mm\textsuperscript{2}}$}}%
      \csname LTb\endcsname
      \put(5098,937){\makebox(0,0)[r]{\strut{}\footnotesize macro ($\Bk\to\Bzero$)}}%
    }%
    \gplgaddtomacro\gplbacktext{%
      \csname LTb\endcsname
      \put(1990,3498){\makebox(0,0)[r]{\strut{}\scriptsize$-0.46\,\,$}}%
      \put(1990,3780){\makebox(0,0)[r]{\strut{}\scriptsize$-0.44\,\,$}}%
      \put(1990,4062){\makebox(0,0)[r]{\strut{}\scriptsize$-0.42\,\,$}}%
      \put(1990,4344){\makebox(0,0)[r]{\strut{}\scriptsize$-0.40\,\,$}}%
      \put(1990,4626){\makebox(0,0)[r]{\strut{}\scriptsize$-0.38\,\,$}}%
      \put(2122,3207){\makebox(0,0){\strut{}\scriptsize$10$}}%
      \put(2952,3207){\makebox(0,0){\strut{}\scriptsize$15$}}%
      \put(3781,3207){\makebox(0,0){\strut{}\scriptsize$20$}}%
    }%
    \gplgaddtomacro\gplfronttext{%
      \csname LTb\endcsname
      \put(1330,4132){\rotatebox{-270}{\makebox(0,0){\strut{}}}}%
      \put(2951,2921){\makebox(0,0){\strut{}}}%
    }%
    \gplbacktext
    \put(0,0){\includegraphics{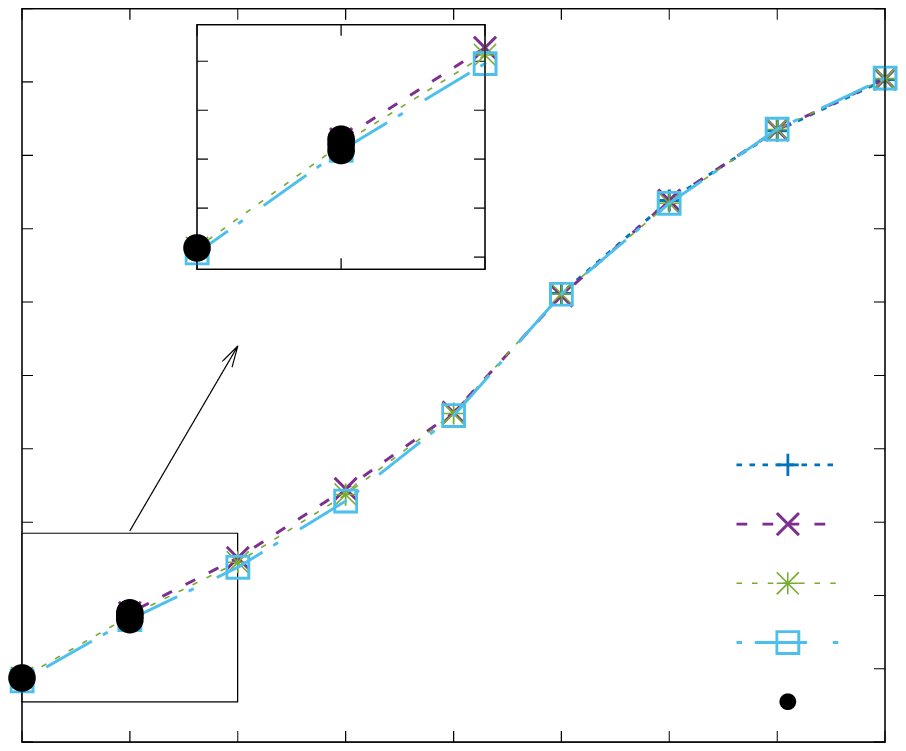}}%
    \gplfronttext
  \end{picture}%
\endgroup

%% file: 2lay_h-vs-time_ga.tex
\begingroup
  \makeatletter
  \providecommand\color[2][]{%
    \GenericError{(gnuplot) \space\space\space\@spaces}{%
      Package color not loaded in conjunction with
      terminal option `colourtext'%
    }{See the gnuplot documentation for explanation.%
    }{Either use 'blacktext' in gnuplot or load the package
      color.sty in LaTeX.}%
    \renewcommand\color[2][]{}%
  }%
  \providecommand\includegraphics[2][]{%
    \GenericError{(gnuplot) \space\space\space\@spaces}{%
      Package graphicx or graphics not loaded%
    }{See the gnuplot documentation for explanation.%
    }{The gnuplot epslatex terminal needs graphicx.sty or graphics.sty.}%
    \renewcommand\includegraphics[2][]{}%
  }%
  \providecommand\rotatebox[2]{#2}%
  \@ifundefined{ifGPcolor}{%
    \newif\ifGPcolor
    \GPcolortrue
  }{}%
  \@ifundefined{ifGPblacktext}{%
    \newif\ifGPblacktext
    \GPblacktexttrue
  }{}%
  \let\gplgaddtomacro\g@addto@macro
  \gdef\gplbacktext{}%
  \gdef\gplfronttext{}%
  \makeatother
  \ifGPblacktext
    \def\colorrgb#1{}%
    \def\colorgray#1{}%
  \else
    \ifGPcolor
      \def\colorrgb#1{\color[rgb]{#1}}%
      \def\colorgray#1{\color[gray]{#1}}%
      \expandafter\def\csname LTw\endcsname{\color{white}}%
      \expandafter\def\csname LTb\endcsname{\color{black}}%
      \expandafter\def\csname LTa\endcsname{\color{black}}%
      \expandafter\def\csname LT0\endcsname{\color[rgb]{1,0,0}}%
      \expandafter\def\csname LT1\endcsname{\color[rgb]{0,1,0}}%
      \expandafter\def\csname LT2\endcsname{\color[rgb]{0,0,1}}%
      \expandafter\def\csname LT3\endcsname{\color[rgb]{1,0,1}}%
      \expandafter\def\csname LT4\endcsname{\color[rgb]{0,1,1}}%
      \expandafter\def\csname LT5\endcsname{\color[rgb]{1,1,0}}%
      \expandafter\def\csname LT6\endcsname{\color[rgb]{0,0,0}}%
      \expandafter\def\csname LT7\endcsname{\color[rgb]{1,0.3,0}}%
      \expandafter\def\csname LT8\endcsname{\color[rgb]{0.5,0.5,0.5}}%
    \else
      \def\colorrgb#1{\color{black}}%
      \def\colorgray#1{\color[gray]{#1}}%
      \expandafter\def\csname LTw\endcsname{\color{white}}%
      \expandafter\def\csname LTb\endcsname{\color{black}}%
      \expandafter\def\csname LTa\endcsname{\color{black}}%
      \expandafter\def\csname LT0\endcsname{\color{black}}%
      \expandafter\def\csname LT1\endcsname{\color{black}}%
      \expandafter\def\csname LT2\endcsname{\color{black}}%
      \expandafter\def\csname LT3\endcsname{\color{black}}%
      \expandafter\def\csname LT4\endcsname{\color{black}}%
      \expandafter\def\csname LT5\endcsname{\color{black}}%
      \expandafter\def\csname LT6\endcsname{\color{black}}%
      \expandafter\def\csname LT7\endcsname{\color{black}}%
      \expandafter\def\csname LT8\endcsname{\color{black}}%
    \fi
  \fi
    \setlength{\unitlength}{0.0500bp}%
    \ifx\gptboxheight\undefined%
      \newlength{\gptboxheight}%
      \newlength{\gptboxwidth}%
      \newsavebox{\gptboxtext}%
    \fi%
    \setlength{\fboxrule}{0.5pt}%
    \setlength{\fboxsep}{1pt}%
\begin{picture}(7200.00,5040.00)%
    \gplgaddtomacro\gplbacktext{%
      \csname LTb\endcsname
      \put(982,704){\makebox(0,0)[r]{\strut{}$0.7$}}%
      \put(982,1549){\makebox(0,0)[r]{\strut{}$0.8$}}%
      \put(982,2394){\makebox(0,0)[r]{\strut{}$0.9$}}%
      \put(982,3239){\makebox(0,0)[r]{\strut{}$1.0$}}%
      \put(982,4084){\makebox(0,0)[r]{\strut{}$1.1$}}%
      \put(982,4929){\makebox(0,0)[r]{\strut{}$1.2$}}%
      \put(1114,484){\makebox(0,0){\strut{}$0.02$}}%
      \put(2771,484){\makebox(0,0){\strut{}$0.03$}}%
      \put(4428,484){\makebox(0,0){\strut{}$0.04$}}%
      \put(6085,484){\makebox(0,0){\strut{}$0.05$}}%
    }%
    \gplgaddtomacro\gplfronttext{%
      \csname LTb\endcsname
      \put(366,2816){\rotatebox{-270}{\makebox(0,0){\strut{}\large$t_{crit}/\text{sec}$}}}%
      \put(3599,154){\makebox(0,0){\strut{}\large$\omega/\text{mm}$}}%
      \csname LTb\endcsname
      \put(5098,1960){\makebox(0,0)[r]{\strut{}$\gamma_f/\gamma_m=5.0$}}%
      \csname LTb\endcsname
      \put(5098,1619){\makebox(0,0)[r]{\strut{}$\gamma_f/\gamma_m=20.0$}}%
      \csname LTb\endcsname
      \put(5098,1278){\makebox(0,0)[r]{\strut{}$\gamma_f/\gamma_m=100.0$}}%
      \csname LTb\endcsname
      \put(5098,937){\makebox(0,0)[r]{\strut{}macro $(\Bk\to\Bzero)$}}%
    }%
    \gplbacktext
    \put(0,0){\includegraphics{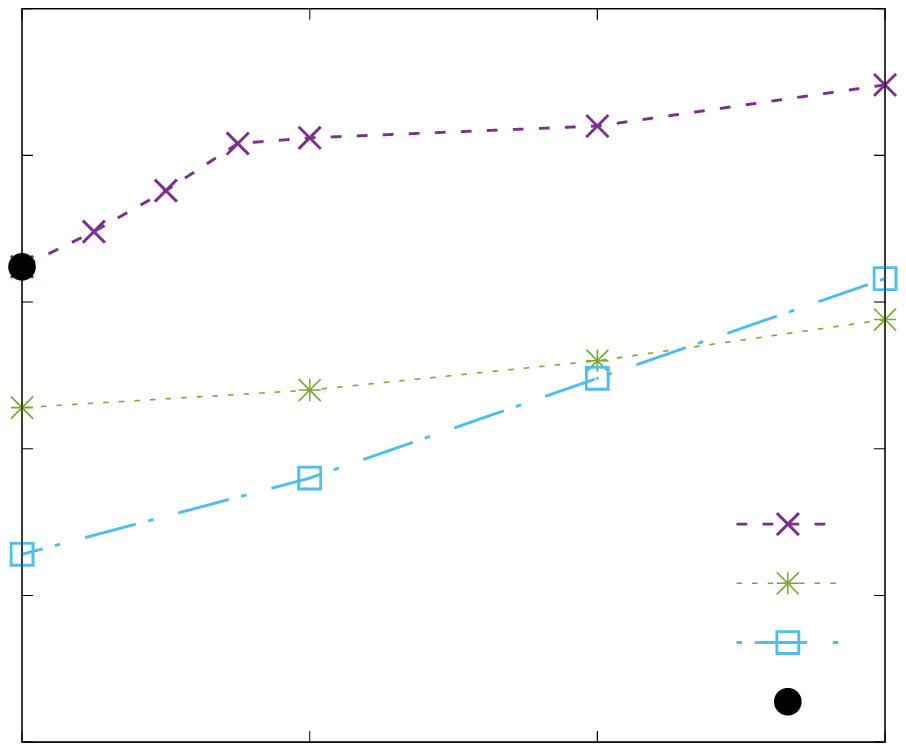}}%
    \gplfronttext
  \end{picture}%
\endgroup

%% file: 2lay_h-vs-pbar11_ga.tex
\begingroup
  \makeatletter
  \providecommand\color[2][]{%
    \GenericError{(gnuplot) \space\space\space\@spaces}{%
      Package color not loaded in conjunction with
      terminal option `colourtext'%
    }{See the gnuplot documentation for explanation.%
    }{Either use 'blacktext' in gnuplot or load the package
      color.sty in LaTeX.}%
    \renewcommand\color[2][]{}%
  }%
  \providecommand\includegraphics[2][]{%
    \GenericError{(gnuplot) \space\space\space\@spaces}{%
      Package graphicx or graphics not loaded%
    }{See the gnuplot documentation for explanation.%
    }{The gnuplot epslatex terminal needs graphicx.sty or graphics.sty.}%
    \renewcommand\includegraphics[2][]{}%
  }%
  \providecommand\rotatebox[2]{#2}%
  \@ifundefined{ifGPcolor}{%
    \newif\ifGPcolor
    \GPcolortrue
  }{}%
  \@ifundefined{ifGPblacktext}{%
    \newif\ifGPblacktext
    \GPblacktexttrue
  }{}%
  \let\gplgaddtomacro\g@addto@macro
  \gdef\gplbacktext{}%
  \gdef\gplfronttext{}%
  \makeatother
  \ifGPblacktext
    \def\colorrgb#1{}%
    \def\colorgray#1{}%
  \else
    \ifGPcolor
      \def\colorrgb#1{\color[rgb]{#1}}%
      \def\colorgray#1{\color[gray]{#1}}%
      \expandafter\def\csname LTw\endcsname{\color{white}}%
      \expandafter\def\csname LTb\endcsname{\color{black}}%
      \expandafter\def\csname LTa\endcsname{\color{black}}%
      \expandafter\def\csname LT0\endcsname{\color[rgb]{1,0,0}}%
      \expandafter\def\csname LT1\endcsname{\color[rgb]{0,1,0}}%
      \expandafter\def\csname LT2\endcsname{\color[rgb]{0,0,1}}%
      \expandafter\def\csname LT3\endcsname{\color[rgb]{1,0,1}}%
      \expandafter\def\csname LT4\endcsname{\color[rgb]{0,1,1}}%
      \expandafter\def\csname LT5\endcsname{\color[rgb]{1,1,0}}%
      \expandafter\def\csname LT6\endcsname{\color[rgb]{0,0,0}}%
      \expandafter\def\csname LT7\endcsname{\color[rgb]{1,0.3,0}}%
      \expandafter\def\csname LT8\endcsname{\color[rgb]{0.5,0.5,0.5}}%
    \else
      \def\colorrgb#1{\color{black}}%
      \def\colorgray#1{\color[gray]{#1}}%
      \expandafter\def\csname LTw\endcsname{\color{white}}%
      \expandafter\def\csname LTb\endcsname{\color{black}}%
      \expandafter\def\csname LTa\endcsname{\color{black}}%
      \expandafter\def\csname LT0\endcsname{\color{black}}%
      \expandafter\def\csname LT1\endcsname{\color{black}}%
      \expandafter\def\csname LT2\endcsname{\color{black}}%
      \expandafter\def\csname LT3\endcsname{\color{black}}%
      \expandafter\def\csname LT4\endcsname{\color{black}}%
      \expandafter\def\csname LT5\endcsname{\color{black}}%
      \expandafter\def\csname LT6\endcsname{\color{black}}%
      \expandafter\def\csname LT7\endcsname{\color{black}}%
      \expandafter\def\csname LT8\endcsname{\color{black}}%
    \fi
  \fi
    \setlength{\unitlength}{0.0500bp}%
    \ifx\gptboxheight\undefined%
      \newlength{\gptboxheight}%
      \newlength{\gptboxwidth}%
      \newsavebox{\gptboxtext}%
    \fi%
    \setlength{\fboxrule}{0.5pt}%
    \setlength{\fboxsep}{1pt}%
\begin{picture}(7200.00,5040.00)%
    \gplgaddtomacro\gplbacktext{%
      \csname LTb\endcsname
      \put(982,704){\makebox(0,0)[r]{\strut{}$-1.0$}}%
      \put(982,1232){\makebox(0,0)[r]{\strut{}$-0.9$}}%
      \put(982,1760){\makebox(0,0)[r]{\strut{}$-0.8$}}%
      \put(982,2288){\makebox(0,0)[r]{\strut{}$-0.7$}}%
      \put(982,2816){\makebox(0,0)[r]{\strut{}$-0.6$}}%
      \put(982,3345){\makebox(0,0)[r]{\strut{}$-0.5$}}%
      \put(982,3873){\makebox(0,0)[r]{\strut{}$-0.4$}}%
      \put(982,4401){\makebox(0,0)[r]{\strut{}$-0.3$}}%
      \put(982,4929){\makebox(0,0)[r]{\strut{}$-0.2$}}%
      \put(1114,484){\makebox(0,0){\strut{}$0.02$}}%
      \put(2771,484){\makebox(0,0){\strut{}$0.03$}}%
      \put(4428,484){\makebox(0,0){\strut{}$0.04$}}%
      \put(6085,484){\makebox(0,0){\strut{}$0.05$}}%
    }%
    \gplgaddtomacro\gplfronttext{%
      \csname LTb\endcsname
      \put(234,2816){\rotatebox{-270}{\makebox(0,0){\strut{}\large$\lbar{P}_{11}/\gamma_m$}}}%
      \put(3599,154){\makebox(0,0){\strut{}\large $\omega/\text{mm}$}}%
      \csname LTb\endcsname
      \put(2962,2301){\makebox(0,0)[r]{\strut{}$\gamma_f/\gamma_m=1.0$}}%
      \csname LTb\endcsname
      \put(2962,1960){\makebox(0,0)[r]{\strut{}$\gamma_f/\gamma_m=5.0$}}%
      \csname LTb\endcsname
      \put(2962,1619){\makebox(0,0)[r]{\strut{}$\gamma_f/\gamma_m=20.0$}}%
      \csname LTb\endcsname
      \put(2962,1278){\makebox(0,0)[r]{\strut{}$\gamma_f/\gamma_m=100.0$}}%
      \csname LTb\endcsname
      \put(2962,937){\makebox(0,0)[r]{\strut{}macro $(\Bk\to\Bzero)$}}%
    }%
    \gplbacktext
    \put(0,0){\includegraphics{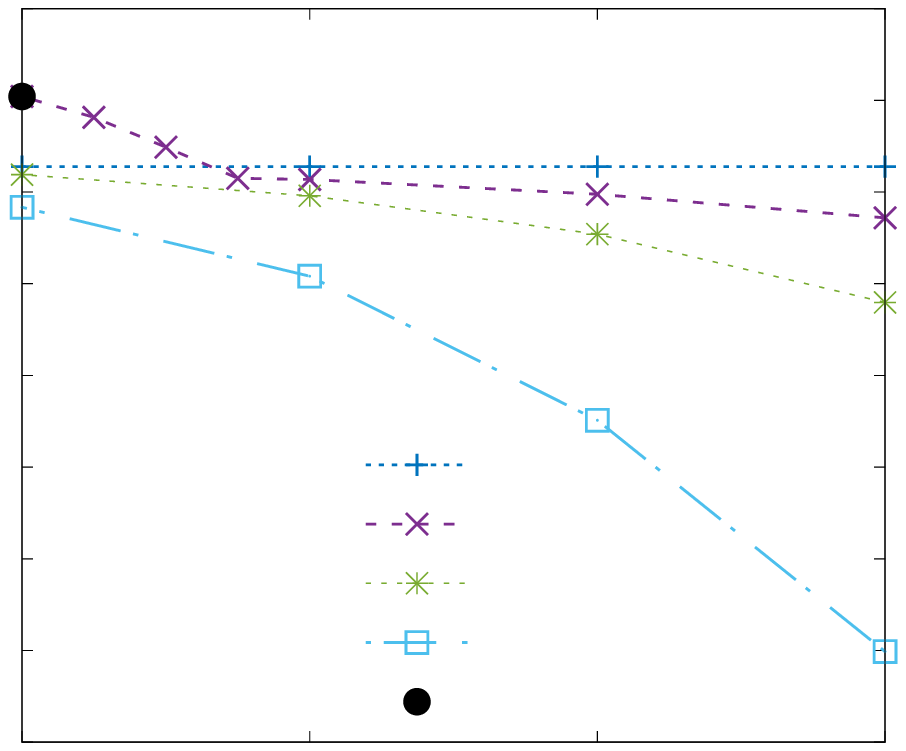}}%
    \gplfronttext
  \end{picture}%
\endgroup

%% file: 2lay_gamma-vs-time_h.tex
\begingroup
  \makeatletter
  \providecommand\color[2][]{%
    \GenericError{(gnuplot) \space\space\space\@spaces}{%
      Package color not loaded in conjunction with
      terminal option `colourtext'%
    }{See the gnuplot documentation for explanation.%
    }{Either use 'blacktext' in gnuplot or load the package
      color.sty in LaTeX.}%
    \renewcommand\color[2][]{}%
  }%
  \providecommand\includegraphics[2][]{%
    \GenericError{(gnuplot) \space\space\space\@spaces}{%
      Package graphicx or graphics not loaded%
    }{See the gnuplot documentation for explanation.%
    }{The gnuplot epslatex terminal needs graphicx.sty or graphics.sty.}%
    \renewcommand\includegraphics[2][]{}%
  }%
  \providecommand\rotatebox[2]{#2}%
  \@ifundefined{ifGPcolor}{%
    \newif\ifGPcolor
    \GPcolortrue
  }{}%
  \@ifundefined{ifGPblacktext}{%
    \newif\ifGPblacktext
    \GPblacktexttrue
  }{}%
  \let\gplgaddtomacro\g@addto@macro
  \gdef\gplbacktext{}%
  \gdef\gplfronttext{}%
  \makeatother
  \ifGPblacktext
    \def\colorrgb#1{}%
    \def\colorgray#1{}%
  \else
    \ifGPcolor
      \def\colorrgb#1{\color[rgb]{#1}}%
      \def\colorgray#1{\color[gray]{#1}}%
      \expandafter\def\csname LTw\endcsname{\color{white}}%
      \expandafter\def\csname LTb\endcsname{\color{black}}%
      \expandafter\def\csname LTa\endcsname{\color{black}}%
      \expandafter\def\csname LT0\endcsname{\color[rgb]{1,0,0}}%
      \expandafter\def\csname LT1\endcsname{\color[rgb]{0,1,0}}%
      \expandafter\def\csname LT2\endcsname{\color[rgb]{0,0,1}}%
      \expandafter\def\csname LT3\endcsname{\color[rgb]{1,0,1}}%
      \expandafter\def\csname LT4\endcsname{\color[rgb]{0,1,1}}%
      \expandafter\def\csname LT5\endcsname{\color[rgb]{1,1,0}}%
      \expandafter\def\csname LT6\endcsname{\color[rgb]{0,0,0}}%
      \expandafter\def\csname LT7\endcsname{\color[rgb]{1,0.3,0}}%
      \expandafter\def\csname LT8\endcsname{\color[rgb]{0.5,0.5,0.5}}%
    \else
      \def\colorrgb#1{\color{black}}%
      \def\colorgray#1{\color[gray]{#1}}%
      \expandafter\def\csname LTw\endcsname{\color{white}}%
      \expandafter\def\csname LTb\endcsname{\color{black}}%
      \expandafter\def\csname LTa\endcsname{\color{black}}%
      \expandafter\def\csname LT0\endcsname{\color{black}}%
      \expandafter\def\csname LT1\endcsname{\color{black}}%
      \expandafter\def\csname LT2\endcsname{\color{black}}%
      \expandafter\def\csname LT3\endcsname{\color{black}}%
      \expandafter\def\csname LT4\endcsname{\color{black}}%
      \expandafter\def\csname LT5\endcsname{\color{black}}%
      \expandafter\def\csname LT6\endcsname{\color{black}}%
      \expandafter\def\csname LT7\endcsname{\color{black}}%
      \expandafter\def\csname LT8\endcsname{\color{black}}%
    \fi
  \fi
    \setlength{\unitlength}{0.0500bp}%
    \ifx\gptboxheight\undefined%
      \newlength{\gptboxheight}%
      \newlength{\gptboxwidth}%
      \newsavebox{\gptboxtext}%
    \fi%
    \setlength{\fboxrule}{0.5pt}%
    \setlength{\fboxsep}{1pt}%
\begin{picture}(7200.00,5040.00)%
    \gplgaddtomacro\gplbacktext{%
      \csname LTb\endcsname
      \put(982,704){\makebox(0,0)[r]{\strut{}$0.80$}}%
      \put(982,1173){\makebox(0,0)[r]{\strut{}$0.85$}}%
      \put(982,1643){\makebox(0,0)[r]{\strut{}$0.90$}}%
      \put(982,2112){\makebox(0,0)[r]{\strut{}$0.95$}}%
      \put(982,2582){\makebox(0,0)[r]{\strut{}$1.00$}}%
      \put(982,3051){\makebox(0,0)[r]{\strut{}$1.05$}}%
      \put(982,3521){\makebox(0,0)[r]{\strut{}$1.10$}}%
      \put(982,3990){\makebox(0,0)[r]{\strut{}$1.15$}}%
      \put(982,4460){\makebox(0,0)[r]{\strut{}$1.20$}}%
      \put(982,4929){\makebox(0,0)[r]{\strut{}$1.25$}}%
      \put(1114,484){\makebox(0,0){\strut{}$0$}}%
      \put(1611,484){\makebox(0,0){\strut{}$10$}}%
      \put(2108,484){\makebox(0,0){\strut{}$20$}}%
      \put(2605,484){\makebox(0,0){\strut{}$30$}}%
      \put(3102,484){\makebox(0,0){\strut{}$40$}}%
      \put(3600,484){\makebox(0,0){\strut{}$50$}}%
      \put(4097,484){\makebox(0,0){\strut{}$60$}}%
      \put(4594,484){\makebox(0,0){\strut{}$70$}}%
      \put(5091,484){\makebox(0,0){\strut{}$80$}}%
      \put(5588,484){\makebox(0,0){\strut{}$90$}}%
      \put(6085,484){\makebox(0,0){\strut{}$100$}}%
    }%
    \gplgaddtomacro\gplfronttext{%
      \csname LTb\endcsname
      \put(234,2816){\rotatebox{-270}{\makebox(0,0){\strut{}\large $t_{crit}/sec$}}}%
      \put(3599,154){\makebox(0,0){\strut{}\large $\gamma_f/\gamma_m$}}%
      \csname LTb\endcsname
      \put(5098,4695){\makebox(0,0)[r]{\strut{}$\omega=0.02\,\text{mm}$}}%
      \csname LTb\endcsname
      \put(5098,4354){\makebox(0,0)[r]{\strut{}$\omega=0.03\,\text{mm}$}}%
      \csname LTb\endcsname
      \put(5098,4013){\makebox(0,0)[r]{\strut{}$\omega=0.04\,\text{mm}$}}%
      \csname LTb\endcsname
      \put(5098,3672){\makebox(0,0)[r]{\strut{}$\omega=0.05\,\text{mm}$}}%
      \csname LTb\endcsname
      \put(5098,3331){\makebox(0,0)[r]{\strut{}macro ($\Bk\to\Bzero$)}}%
    }%
    \gplbacktext
    \put(0,0){\includegraphics{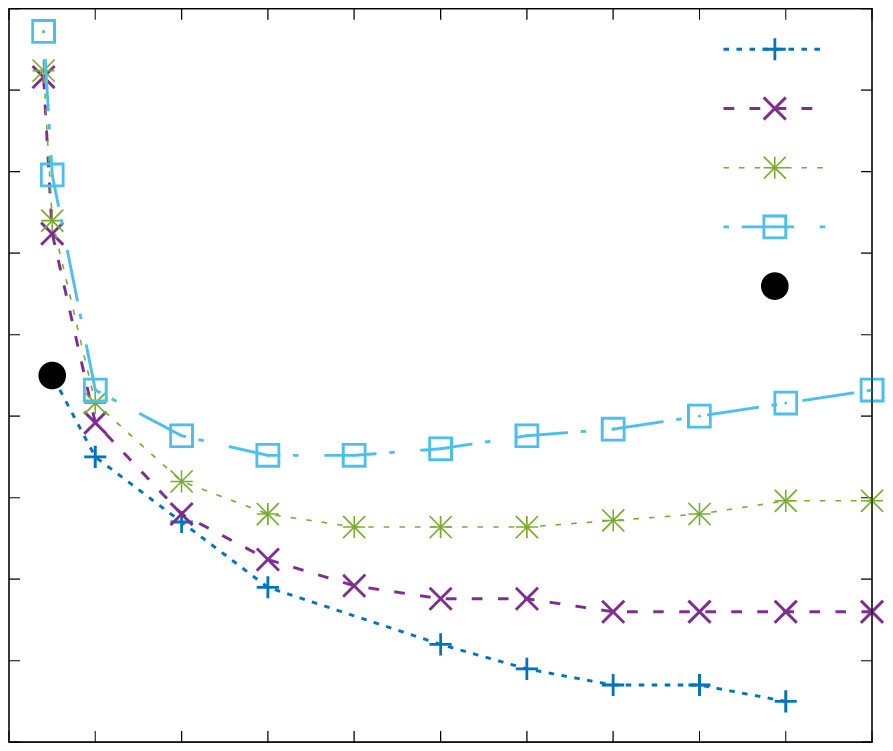}}%
    \gplfronttext
  \end{picture}%
\endgroup

%% file: 2lay_gamma-vs-pbar11_h.tex
\begingroup
  \makeatletter
  \providecommand\color[2][]{%
    \GenericError{(gnuplot) \space\space\space\@spaces}{%
      Package color not loaded in conjunction with
      terminal option `colourtext'%
    }{See the gnuplot documentation for explanation.%
    }{Either use 'blacktext' in gnuplot or load the package
      color.sty in LaTeX.}%
    \renewcommand\color[2][]{}%
  }%
  \providecommand\includegraphics[2][]{%
    \GenericError{(gnuplot) \space\space\space\@spaces}{%
      Package graphicx or graphics not loaded%
    }{See the gnuplot documentation for explanation.%
    }{The gnuplot epslatex terminal needs graphicx.sty or graphics.sty.}%
    \renewcommand\includegraphics[2][]{}%
  }%
  \providecommand\rotatebox[2]{#2}%
  \@ifundefined{ifGPcolor}{%
    \newif\ifGPcolor
    \GPcolortrue
  }{}%
  \@ifundefined{ifGPblacktext}{%
    \newif\ifGPblacktext
    \GPblacktexttrue
  }{}%
  \let\gplgaddtomacro\g@addto@macro
  \gdef\gplbacktext{}%
  \gdef\gplfronttext{}%
  \makeatother
  \ifGPblacktext
    \def\colorrgb#1{}%
    \def\colorgray#1{}%
  \else
    \ifGPcolor
      \def\colorrgb#1{\color[rgb]{#1}}%
      \def\colorgray#1{\color[gray]{#1}}%
      \expandafter\def\csname LTw\endcsname{\color{white}}%
      \expandafter\def\csname LTb\endcsname{\color{black}}%
      \expandafter\def\csname LTa\endcsname{\color{black}}%
      \expandafter\def\csname LT0\endcsname{\color[rgb]{1,0,0}}%
      \expandafter\def\csname LT1\endcsname{\color[rgb]{0,1,0}}%
      \expandafter\def\csname LT2\endcsname{\color[rgb]{0,0,1}}%
      \expandafter\def\csname LT3\endcsname{\color[rgb]{1,0,1}}%
      \expandafter\def\csname LT4\endcsname{\color[rgb]{0,1,1}}%
      \expandafter\def\csname LT5\endcsname{\color[rgb]{1,1,0}}%
      \expandafter\def\csname LT6\endcsname{\color[rgb]{0,0,0}}%
      \expandafter\def\csname LT7\endcsname{\color[rgb]{1,0.3,0}}%
      \expandafter\def\csname LT8\endcsname{\color[rgb]{0.5,0.5,0.5}}%
    \else
      \def\colorrgb#1{\color{black}}%
      \def\colorgray#1{\color[gray]{#1}}%
      \expandafter\def\csname LTw\endcsname{\color{white}}%
      \expandafter\def\csname LTb\endcsname{\color{black}}%
      \expandafter\def\csname LTa\endcsname{\color{black}}%
      \expandafter\def\csname LT0\endcsname{\color{black}}%
      \expandafter\def\csname LT1\endcsname{\color{black}}%
      \expandafter\def\csname LT2\endcsname{\color{black}}%
      \expandafter\def\csname LT3\endcsname{\color{black}}%
      \expandafter\def\csname LT4\endcsname{\color{black}}%
      \expandafter\def\csname LT5\endcsname{\color{black}}%
      \expandafter\def\csname LT6\endcsname{\color{black}}%
      \expandafter\def\csname LT7\endcsname{\color{black}}%
      \expandafter\def\csname LT8\endcsname{\color{black}}%
    \fi
  \fi
    \setlength{\unitlength}{0.0500bp}%
    \ifx\gptboxheight\undefined%
      \newlength{\gptboxheight}%
      \newlength{\gptboxwidth}%
      \newsavebox{\gptboxtext}%
    \fi%
    \setlength{\fboxrule}{0.5pt}%
    \setlength{\fboxsep}{1pt}%
\begin{picture}(7200.00,5040.00)%
    \gplgaddtomacro\gplbacktext{%
      \csname LTb\endcsname
      \put(982,704){\makebox(0,0)[r]{\strut{}$-1.0$}}%
      \put(982,1232){\makebox(0,0)[r]{\strut{}$-0.9$}}%
      \put(982,1760){\makebox(0,0)[r]{\strut{}$-0.8$}}%
      \put(982,2288){\makebox(0,0)[r]{\strut{}$-0.7$}}%
      \put(982,2816){\makebox(0,0)[r]{\strut{}$-0.6$}}%
      \put(982,3345){\makebox(0,0)[r]{\strut{}$-0.5$}}%
      \put(982,3873){\makebox(0,0)[r]{\strut{}$-0.4$}}%
      \put(982,4401){\makebox(0,0)[r]{\strut{}$-0.3$}}%
      \put(982,4929){\makebox(0,0)[r]{\strut{}$-0.2$}}%
      \put(1114,484){\makebox(0,0){\strut{}$0$}}%
      \put(1611,484){\makebox(0,0){\strut{}$10$}}%
      \put(2108,484){\makebox(0,0){\strut{}$20$}}%
      \put(2605,484){\makebox(0,0){\strut{}$30$}}%
      \put(3102,484){\makebox(0,0){\strut{}$40$}}%
      \put(3600,484){\makebox(0,0){\strut{}$50$}}%
      \put(4097,484){\makebox(0,0){\strut{}$60$}}%
      \put(4594,484){\makebox(0,0){\strut{}$70$}}%
      \put(5091,484){\makebox(0,0){\strut{}$80$}}%
      \put(5588,484){\makebox(0,0){\strut{}$90$}}%
      \put(6085,484){\makebox(0,0){\strut{}$100$}}%
    }%
    \gplgaddtomacro\gplfronttext{%
      \csname LTb\endcsname
      \put(234,2816){\rotatebox{-270}{\makebox(0,0){\strut{}\large $\lbar{P}_{11}/\gamma_m$}}}%
      \put(3599,154){\makebox(0,0){\strut{}\large $\gamma_f/\gamma_m$}}%
      \csname LTb\endcsname
      \put(2698,2301){\makebox(0,0)[r]{\strut{}$\omega=0.02\,\text{mm}$}}%
      \csname LTb\endcsname
      \put(2698,1960){\makebox(0,0)[r]{\strut{}$\omega=0.03\,\text{mm}$}}%
      \csname LTb\endcsname
      \put(2698,1619){\makebox(0,0)[r]{\strut{}$\omega=0.04\,\text{mm}$}}%
      \csname LTb\endcsname
      \put(2698,1278){\makebox(0,0)[r]{\strut{}$\omega=0.05\,\text{mm}$}}%
      \csname LTb\endcsname
      \put(2698,937){\makebox(0,0)[r]{\strut{}macro ($\Bk\to\Bzero$)}}%
    }%
    \gplbacktext
    \put(0,0){\includegraphics{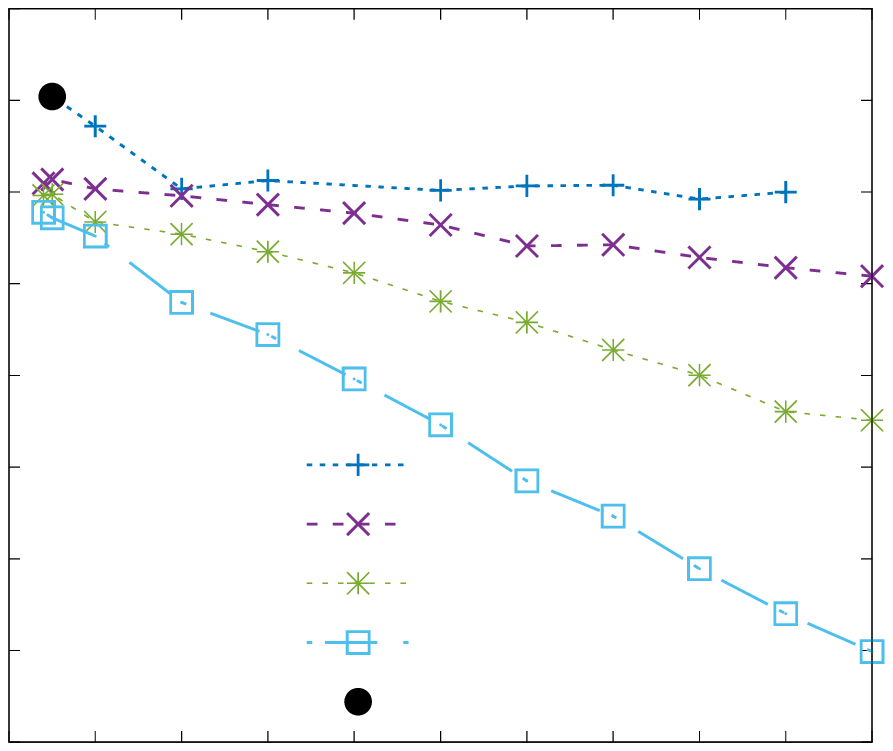}}%
    \gplfronttext
  \end{picture}%
\endgroup

%% file: 2lay_time-vs-h_M.tex
\begingroup
  \makeatletter
  \providecommand\color[2][]{%
    \GenericError{(gnuplot) \space\space\space\@spaces}{%
      Package color not loaded in conjunction with
      terminal option `colourtext'%
    }{See the gnuplot documentation for explanation.%
    }{Either use 'blacktext' in gnuplot or load the package
      color.sty in LaTeX.}%
    \renewcommand\color[2][]{}%
  }%
  \providecommand\includegraphics[2][]{%
    \GenericError{(gnuplot) \space\space\space\@spaces}{%
      Package graphicx or graphics not loaded%
    }{See the gnuplot documentation for explanation.%
    }{The gnuplot epslatex terminal needs graphicx.sty or graphics.sty.}%
    \renewcommand\includegraphics[2][]{}%
  }%
  \providecommand\rotatebox[2]{#2}%
  \@ifundefined{ifGPcolor}{%
    \newif\ifGPcolor
    \GPcolortrue
  }{}%
  \@ifundefined{ifGPblacktext}{%
    \newif\ifGPblacktext
    \GPblacktexttrue
  }{}%
  \let\gplgaddtomacro\g@addto@macro
  \gdef\gplbacktext{}%
  \gdef\gplfronttext{}%
  \makeatother
  \ifGPblacktext
    \def\colorrgb#1{}%
    \def\colorgray#1{}%
  \else
    \ifGPcolor
      \def\colorrgb#1{\color[rgb]{#1}}%
      \def\colorgray#1{\color[gray]{#1}}%
      \expandafter\def\csname LTw\endcsname{\color{white}}%
      \expandafter\def\csname LTb\endcsname{\color{black}}%
      \expandafter\def\csname LTa\endcsname{\color{black}}%
      \expandafter\def\csname LT0\endcsname{\color[rgb]{1,0,0}}%
      \expandafter\def\csname LT1\endcsname{\color[rgb]{0,1,0}}%
      \expandafter\def\csname LT2\endcsname{\color[rgb]{0,0,1}}%
      \expandafter\def\csname LT3\endcsname{\color[rgb]{1,0,1}}%
      \expandafter\def\csname LT4\endcsname{\color[rgb]{0,1,1}}%
      \expandafter\def\csname LT5\endcsname{\color[rgb]{1,1,0}}%
      \expandafter\def\csname LT6\endcsname{\color[rgb]{0,0,0}}%
      \expandafter\def\csname LT7\endcsname{\color[rgb]{1,0.3,0}}%
      \expandafter\def\csname LT8\endcsname{\color[rgb]{0.5,0.5,0.5}}%
    \else
      \def\colorrgb#1{\color{black}}%
      \def\colorgray#1{\color[gray]{#1}}%
      \expandafter\def\csname LTw\endcsname{\color{white}}%
      \expandafter\def\csname LTb\endcsname{\color{black}}%
      \expandafter\def\csname LTa\endcsname{\color{black}}%
      \expandafter\def\csname LT0\endcsname{\color{black}}%
      \expandafter\def\csname LT1\endcsname{\color{black}}%
      \expandafter\def\csname LT2\endcsname{\color{black}}%
      \expandafter\def\csname LT3\endcsname{\color{black}}%
      \expandafter\def\csname LT4\endcsname{\color{black}}%
      \expandafter\def\csname LT5\endcsname{\color{black}}%
      \expandafter\def\csname LT6\endcsname{\color{black}}%
      \expandafter\def\csname LT7\endcsname{\color{black}}%
      \expandafter\def\csname LT8\endcsname{\color{black}}%
    \fi
  \fi
    \setlength{\unitlength}{0.0500bp}%
    \ifx\gptboxheight\undefined%
      \newlength{\gptboxheight}%
      \newlength{\gptboxwidth}%
      \newsavebox{\gptboxtext}%
    \fi%
    \setlength{\fboxrule}{0.5pt}%
    \setlength{\fboxsep}{1pt}%
\begin{picture}(7200.00,5040.00)%
    \gplgaddtomacro\gplbacktext{%
      \csname LTb\endcsname
      \put(982,704){\makebox(0,0)[r]{\strut{}$1.00$}}%
      \put(982,1173){\makebox(0,0)[r]{\strut{}$1.02$}}%
      \put(982,1643){\makebox(0,0)[r]{\strut{}$1.04$}}%
      \put(982,2112){\makebox(0,0)[r]{\strut{}$1.06$}}%
      \put(982,2582){\makebox(0,0)[r]{\strut{}$1.08$}}%
      \put(982,3051){\makebox(0,0)[r]{\strut{}$1.10$}}%
      \put(982,3521){\makebox(0,0)[r]{\strut{}$1.12$}}%
      \put(982,3990){\makebox(0,0)[r]{\strut{}$1.14$}}%
      \put(982,4460){\makebox(0,0)[r]{\strut{}$1.16$}}%
      \put(982,4929){\makebox(0,0)[r]{\strut{}$1.18$}}%
      \put(1114,484){\makebox(0,0){\strut{}$0.02$}}%
      \put(2771,484){\makebox(0,0){\strut{}$0.03$}}%
      \put(4428,484){\makebox(0,0){\strut{}$0.04$}}%
      \put(6085,484){\makebox(0,0){\strut{}$0.05$}}%
    }%
    \gplgaddtomacro\gplfronttext{%
      \csname LTb\endcsname
      \put(234,2816){\rotatebox{-270}{\makebox(0,0){\strut{}\large $t_{crit}/\text{sec}$}}}%
      \put(3599,154){\makebox(0,0){\strut{}\large $\omega/\text{mm}$}}%
      \csname LTb\endcsname
      \put(5098,2642){\makebox(0,0)[r]{\strut{}\footnotesize$M_f/M_m=0.01$}}%
      \csname LTb\endcsname
      \put(5098,2301){\makebox(0,0)[r]{\strut{}\footnotesize$M_f/M_m=0.1$}}%
      \csname LTb\endcsname
      \put(5098,1960){\makebox(0,0)[r]{\strut{}\footnotesize$M_f/M_m=1.0$}}%
      \csname LTb\endcsname
      \put(5098,1619){\makebox(0,0)[r]{\strut{}\footnotesize$M_f/M_m=10.0$}}%
      \csname LTb\endcsname
      \put(5098,1278){\makebox(0,0)[r]{\strut{}\footnotesize$M_f/M_m=100.0$}}%
      \csname LTb\endcsname
      \put(5098,937){\makebox(0,0)[r]{\strut{}\footnotesize macro ($\Bk\to\Bzero$)}}%
    }%
    \gplbacktext
    \put(0,0){\includegraphics{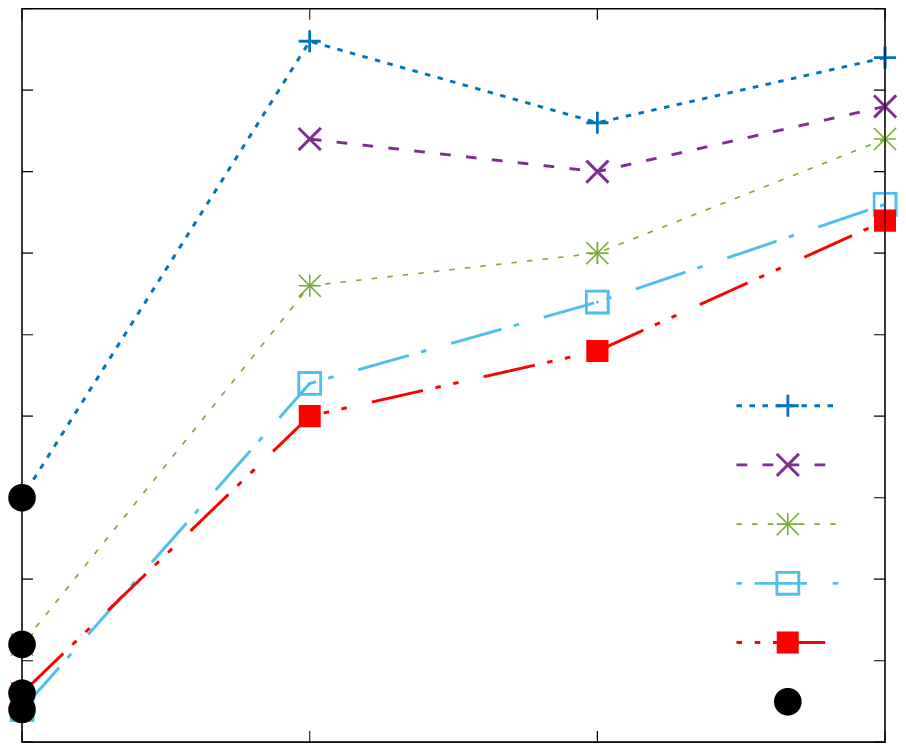}}%
    \gplfronttext
  \end{picture}%
\endgroup

%% file: 2lay_p11bar-vs-h_M.tex
\begingroup
  \makeatletter
  \providecommand\color[2][]{%
    \GenericError{(gnuplot) \space\space\space\@spaces}{%
      Package color not loaded in conjunction with
      terminal option `colourtext'%
    }{See the gnuplot documentation for explanation.%
    }{Either use 'blacktext' in gnuplot or load the package
      color.sty in LaTeX.}%
    \renewcommand\color[2][]{}%
  }%
  \providecommand\includegraphics[2][]{%
    \GenericError{(gnuplot) \space\space\space\@spaces}{%
      Package graphicx or graphics not loaded%
    }{See the gnuplot documentation for explanation.%
    }{The gnuplot epslatex terminal needs graphicx.sty or graphics.sty.}%
    \renewcommand\includegraphics[2][]{}%
  }%
  \providecommand\rotatebox[2]{#2}%
  \@ifundefined{ifGPcolor}{%
    \newif\ifGPcolor
    \GPcolortrue
  }{}%
  \@ifundefined{ifGPblacktext}{%
    \newif\ifGPblacktext
    \GPblacktexttrue
  }{}%
  \let\gplgaddtomacro\g@addto@macro
  \gdef\gplbacktext{}%
  \gdef\gplfronttext{}%
  \makeatother
  \ifGPblacktext
    \def\colorrgb#1{}%
    \def\colorgray#1{}%
  \else
    \ifGPcolor
      \def\colorrgb#1{\color[rgb]{#1}}%
      \def\colorgray#1{\color[gray]{#1}}%
      \expandafter\def\csname LTw\endcsname{\color{white}}%
      \expandafter\def\csname LTb\endcsname{\color{black}}%
      \expandafter\def\csname LTa\endcsname{\color{black}}%
      \expandafter\def\csname LT0\endcsname{\color[rgb]{1,0,0}}%
      \expandafter\def\csname LT1\endcsname{\color[rgb]{0,1,0}}%
      \expandafter\def\csname LT2\endcsname{\color[rgb]{0,0,1}}%
      \expandafter\def\csname LT3\endcsname{\color[rgb]{1,0,1}}%
      \expandafter\def\csname LT4\endcsname{\color[rgb]{0,1,1}}%
      \expandafter\def\csname LT5\endcsname{\color[rgb]{1,1,0}}%
      \expandafter\def\csname LT6\endcsname{\color[rgb]{0,0,0}}%
      \expandafter\def\csname LT7\endcsname{\color[rgb]{1,0.3,0}}%
      \expandafter\def\csname LT8\endcsname{\color[rgb]{0.5,0.5,0.5}}%
    \else
      \def\colorrgb#1{\color{black}}%
      \def\colorgray#1{\color[gray]{#1}}%
      \expandafter\def\csname LTw\endcsname{\color{white}}%
      \expandafter\def\csname LTb\endcsname{\color{black}}%
      \expandafter\def\csname LTa\endcsname{\color{black}}%
      \expandafter\def\csname LT0\endcsname{\color{black}}%
      \expandafter\def\csname LT1\endcsname{\color{black}}%
      \expandafter\def\csname LT2\endcsname{\color{black}}%
      \expandafter\def\csname LT3\endcsname{\color{black}}%
      \expandafter\def\csname LT4\endcsname{\color{black}}%
      \expandafter\def\csname LT5\endcsname{\color{black}}%
      \expandafter\def\csname LT6\endcsname{\color{black}}%
      \expandafter\def\csname LT7\endcsname{\color{black}}%
      \expandafter\def\csname LT8\endcsname{\color{black}}%
    \fi
  \fi
    \setlength{\unitlength}{0.0500bp}%
    \ifx\gptboxheight\undefined%
      \newlength{\gptboxheight}%
      \newlength{\gptboxwidth}%
      \newsavebox{\gptboxtext}%
    \fi%
    \setlength{\fboxrule}{0.5pt}%
    \setlength{\fboxsep}{1pt}%
\begin{picture}(7200.00,5040.00)%
    \gplgaddtomacro\gplbacktext{%
      \csname LTb\endcsname
      \put(982,704){\makebox(0,0)[r]{\strut{}$-0.44$}}%
      \put(982,1232){\makebox(0,0)[r]{\strut{}$-0.42$}}%
      \put(982,1760){\makebox(0,0)[r]{\strut{}$-0.40$}}%
      \put(982,2288){\makebox(0,0)[r]{\strut{}$-0.38$}}%
      \put(982,2817){\makebox(0,0)[r]{\strut{}$-0.36$}}%
      \put(982,3345){\makebox(0,0)[r]{\strut{}$-0.34$}}%
      \put(982,3873){\makebox(0,0)[r]{\strut{}$-0.32$}}%
      \put(982,4401){\makebox(0,0)[r]{\strut{}$-0.30$}}%
      \put(982,4929){\makebox(0,0)[r]{\strut{}$-0.28$}}%
      \put(1114,484){\makebox(0,0){\strut{}$0.02$}}%
      \put(2771,484){\makebox(0,0){\strut{}$0.03$}}%
      \put(4428,484){\makebox(0,0){\strut{}$0.04$}}%
      \put(6085,484){\makebox(0,0){\strut{}$0.05$}}%
    }%
    \gplgaddtomacro\gplfronttext{%
      \csname LTb\endcsname
      \put(102,2816){\rotatebox{-270}{\makebox(0,0){\strut{}\large $\lbar{P}_{11}/\gamma_m$}}}%
      \put(3599,154){\makebox(0,0){\strut{}\large $\omega/\text{mm}$}}%
      \csname LTb\endcsname
      \put(5098,4695){\makebox(0,0)[r]{\strut{}\footnotesize$M_f/M_m=0.01$}}%
      \csname LTb\endcsname
      \put(5098,4354){\makebox(0,0)[r]{\strut{}\footnotesize$M_f/M_m=0.1$}}%
      \csname LTb\endcsname
      \put(5098,4013){\makebox(0,0)[r]{\strut{}\footnotesize$M_f/M_m=1.0$}}%
      \csname LTb\endcsname
      \put(5098,3672){\makebox(0,0)[r]{\strut{}\footnotesize$M_f/M_m=10.0$}}%
      \csname LTb\endcsname
      \put(5098,3331){\makebox(0,0)[r]{\strut{}\footnotesize$M_f/M_m=100.0$}}%
      \csname LTb\endcsname
      \put(5098,2990){\makebox(0,0)[r]{\strut{}\footnotesize macro ($\Bk\to\Bzero$)}}%
    }%
    \gplbacktext
    \put(0,0){\includegraphics{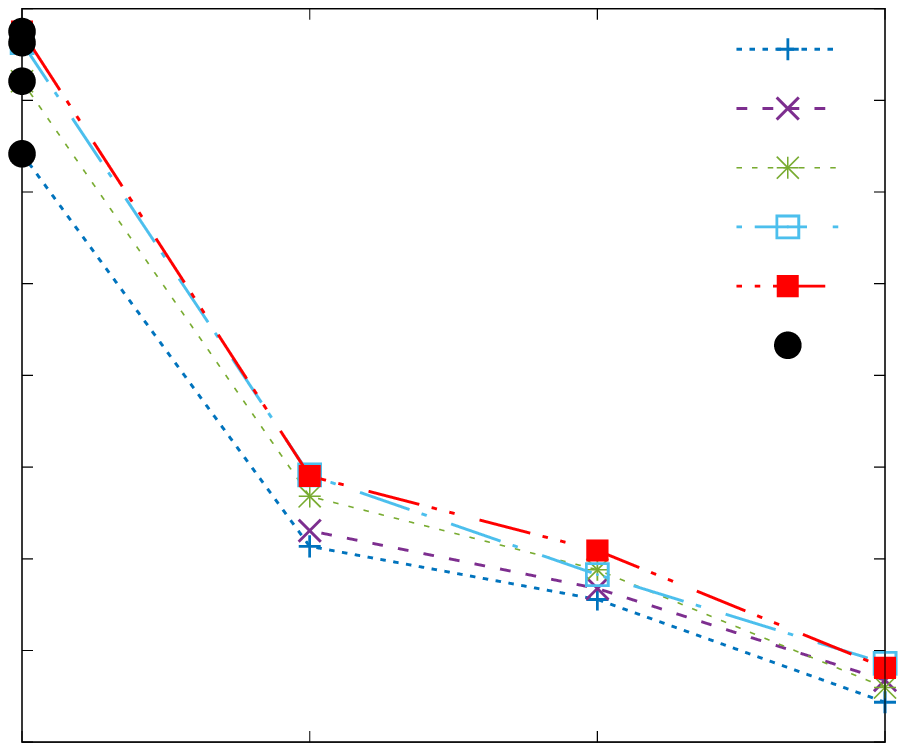}}%
    \gplfronttext
  \end{picture}%
\endgroup

%% file: b-vs-time-alpha.tex
\begingroup
  \makeatletter
  \providecommand\color[2][]{%
    \GenericError{(gnuplot) \space\space\space\@spaces}{%
      Package color not loaded in conjunction with
      terminal option `colourtext'%
    }{See the gnuplot documentation for explanation.%
    }{Either use 'blacktext' in gnuplot or load the package
      color.sty in LaTeX.}%
    \renewcommand\color[2][]{}%
  }%
  \providecommand\includegraphics[2][]{%
    \GenericError{(gnuplot) \space\space\space\@spaces}{%
      Package graphicx or graphics not loaded%
    }{See the gnuplot documentation for explanation.%
    }{The gnuplot epslatex terminal needs graphicx.sty or graphics.sty.}%
    \renewcommand\includegraphics[2][]{}%
  }%
  \providecommand\rotatebox[2]{#2}%
  \@ifundefined{ifGPcolor}{%
    \newif\ifGPcolor
    \GPcolortrue
  }{}%
  \@ifundefined{ifGPblacktext}{%
    \newif\ifGPblacktext
    \GPblacktexttrue
  }{}%
  \let\gplgaddtomacro\g@addto@macro
  \gdef\gplbacktext{}%
  \gdef\gplfronttext{}%
  \makeatother
  \ifGPblacktext
    \def\colorrgb#1{}%
    \def\colorgray#1{}%
  \else
    \ifGPcolor
      \def\colorrgb#1{\color[rgb]{#1}}%
      \def\colorgray#1{\color[gray]{#1}}%
      \expandafter\def\csname LTw\endcsname{\color{white}}%
      \expandafter\def\csname LTb\endcsname{\color{black}}%
      \expandafter\def\csname LTa\endcsname{\color{black}}%
      \expandafter\def\csname LT0\endcsname{\color[rgb]{1,0,0}}%
      \expandafter\def\csname LT1\endcsname{\color[rgb]{0,1,0}}%
      \expandafter\def\csname LT2\endcsname{\color[rgb]{0,0,1}}%
      \expandafter\def\csname LT3\endcsname{\color[rgb]{1,0,1}}%
      \expandafter\def\csname LT4\endcsname{\color[rgb]{0,1,1}}%
      \expandafter\def\csname LT5\endcsname{\color[rgb]{1,1,0}}%
      \expandafter\def\csname LT6\endcsname{\color[rgb]{0,0,0}}%
      \expandafter\def\csname LT7\endcsname{\color[rgb]{1,0.3,0}}%
      \expandafter\def\csname LT8\endcsname{\color[rgb]{0.5,0.5,0.5}}%
    \else
      \def\colorrgb#1{\color{black}}%
      \def\colorgray#1{\color[gray]{#1}}%
      \expandafter\def\csname LTw\endcsname{\color{white}}%
      \expandafter\def\csname LTb\endcsname{\color{black}}%
      \expandafter\def\csname LTa\endcsname{\color{black}}%
      \expandafter\def\csname LT0\endcsname{\color{black}}%
      \expandafter\def\csname LT1\endcsname{\color{black}}%
      \expandafter\def\csname LT2\endcsname{\color{black}}%
      \expandafter\def\csname LT3\endcsname{\color{black}}%
      \expandafter\def\csname LT4\endcsname{\color{black}}%
      \expandafter\def\csname LT5\endcsname{\color{black}}%
      \expandafter\def\csname LT6\endcsname{\color{black}}%
      \expandafter\def\csname LT7\endcsname{\color{black}}%
      \expandafter\def\csname LT8\endcsname{\color{black}}%
    \fi
  \fi
    \setlength{\unitlength}{0.0500bp}%
    \ifx\gptboxheight\undefined%
      \newlength{\gptboxheight}%
      \newlength{\gptboxwidth}%
      \newsavebox{\gptboxtext}%
    \fi%
    \setlength{\fboxrule}{0.5pt}%
    \setlength{\fboxsep}{1pt}%
\begin{picture}(7200.00,5040.00)%
    \gplgaddtomacro\gplbacktext{%
      \csname LTb\endcsname
      \put(982,704){\makebox(0,0)[r]{\strut{}$0.055$}}%
      \put(982,1126){\makebox(0,0)[r]{\strut{}$0.060$}}%
      \put(982,1549){\makebox(0,0)[r]{\strut{}$0.065$}}%
      \put(982,1972){\makebox(0,0)[r]{\strut{}$0.070$}}%
      \put(982,2394){\makebox(0,0)[r]{\strut{}$0.075$}}%
      \put(982,2817){\makebox(0,0)[r]{\strut{}$0.080$}}%
      \put(982,3239){\makebox(0,0)[r]{\strut{}$0.085$}}%
      \put(982,3662){\makebox(0,0)[r]{\strut{}$0.090$}}%
      \put(982,4084){\makebox(0,0)[r]{\strut{}$0.095$}}%
      \put(982,4507){\makebox(0,0)[r]{\strut{}$0.100$}}%
      \put(982,4929){\makebox(0,0)[r]{\strut{}$0.105$}}%
      \put(1114,484){\makebox(0,0){\strut{}$6$}}%
      \put(1735,484){\makebox(0,0){\strut{}}}%
      \put(2357,484){\makebox(0,0){\strut{}$7$}}%
      \put(2978,484){\makebox(0,0){\strut{}}}%
      \put(3600,484){\makebox(0,0){\strut{}$8$}}%
      \put(4221,484){\makebox(0,0){\strut{}}}%
      \put(4842,484){\makebox(0,0){\strut{}$9$}}%
      \put(5464,484){\makebox(0,0){\strut{}}}%
      \put(6085,484){\makebox(0,0){\strut{}$10$}}%
    }%
    \gplgaddtomacro\gplfronttext{%
      \csname LTb\endcsname
      \put(102,2816){\rotatebox{-270}{\makebox(0,0){\strut{}\large $t_{crit}/sec$}}}%
      \put(3599,154){\makebox(0,0){\strut{}\large $b/\text{mm}$}}%
      \csname LTb\endcsname
      \put(5098,4695){\makebox(0,0)[r]{\strut{}$\beta=45^\circ$}}%
      \csname LTb\endcsname
      \put(5098,4354){\makebox(0,0)[r]{\strut{}$\beta=60^\circ$}}%
      \csname LTb\endcsname
      \put(5098,4013){\makebox(0,0)[r]{\strut{}$\beta=75^\circ$}}%
      \csname LTb\endcsname
      \put(5098,3672){\makebox(0,0)[r]{\strut{}$\Bk\to\Bzero$}}%
    }%
    \gplbacktext
    \put(0,0){\includegraphics{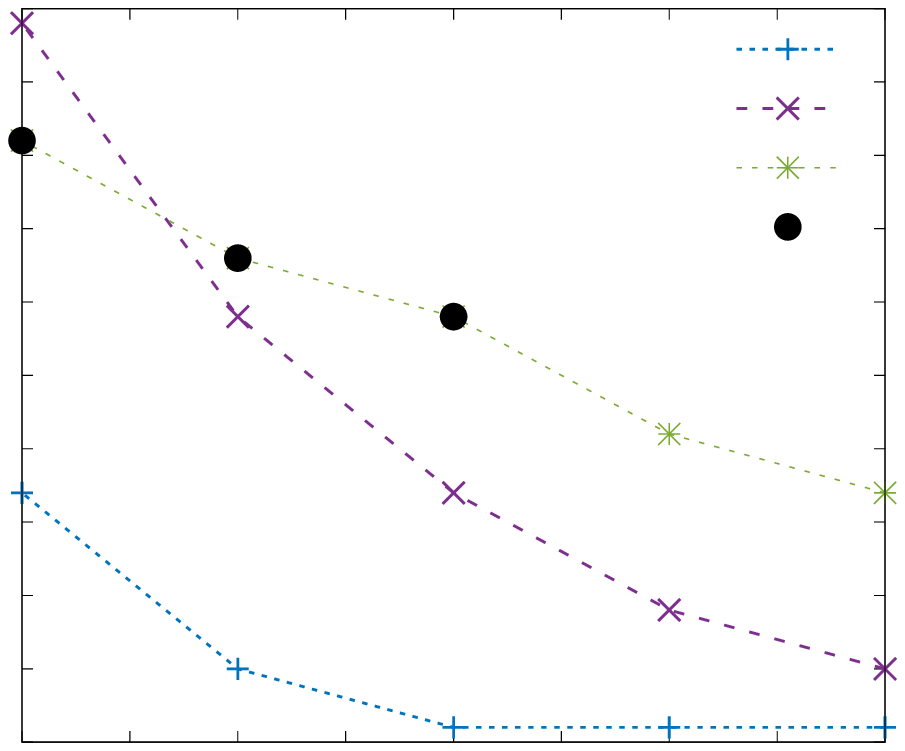}}%
    \gplfronttext
  \end{picture}%
\endgroup

%% file: b-vs-pbar-alpha.tex
\begingroup
  \makeatletter
  \providecommand\color[2][]{%
    \GenericError{(gnuplot) \space\space\space\@spaces}{%
      Package color not loaded in conjunction with
      terminal option `colourtext'%
    }{See the gnuplot documentation for explanation.%
    }{Either use 'blacktext' in gnuplot or load the package
      color.sty in LaTeX.}%
    \renewcommand\color[2][]{}%
  }%
  \providecommand\includegraphics[2][]{%
    \GenericError{(gnuplot) \space\space\space\@spaces}{%
      Package graphicx or graphics not loaded%
    }{See the gnuplot documentation for explanation.%
    }{The gnuplot epslatex terminal needs graphicx.sty or graphics.sty.}%
    \renewcommand\includegraphics[2][]{}%
  }%
  \providecommand\rotatebox[2]{#2}%
  \@ifundefined{ifGPcolor}{%
    \newif\ifGPcolor
    \GPcolortrue
  }{}%
  \@ifundefined{ifGPblacktext}{%
    \newif\ifGPblacktext
    \GPblacktexttrue
  }{}%
  \let\gplgaddtomacro\g@addto@macro
  \gdef\gplbacktext{}%
  \gdef\gplfronttext{}%
  \makeatother
  \ifGPblacktext
    \def\colorrgb#1{}%
    \def\colorgray#1{}%
  \else
    \ifGPcolor
      \def\colorrgb#1{\color[rgb]{#1}}%
      \def\colorgray#1{\color[gray]{#1}}%
      \expandafter\def\csname LTw\endcsname{\color{white}}%
      \expandafter\def\csname LTb\endcsname{\color{black}}%
      \expandafter\def\csname LTa\endcsname{\color{black}}%
      \expandafter\def\csname LT0\endcsname{\color[rgb]{1,0,0}}%
      \expandafter\def\csname LT1\endcsname{\color[rgb]{0,1,0}}%
      \expandafter\def\csname LT2\endcsname{\color[rgb]{0,0,1}}%
      \expandafter\def\csname LT3\endcsname{\color[rgb]{1,0,1}}%
      \expandafter\def\csname LT4\endcsname{\color[rgb]{0,1,1}}%
      \expandafter\def\csname LT5\endcsname{\color[rgb]{1,1,0}}%
      \expandafter\def\csname LT6\endcsname{\color[rgb]{0,0,0}}%
      \expandafter\def\csname LT7\endcsname{\color[rgb]{1,0.3,0}}%
      \expandafter\def\csname LT8\endcsname{\color[rgb]{0.5,0.5,0.5}}%
    \else
      \def\colorrgb#1{\color{black}}%
      \def\colorgray#1{\color[gray]{#1}}%
      \expandafter\def\csname LTw\endcsname{\color{white}}%
      \expandafter\def\csname LTb\endcsname{\color{black}}%
      \expandafter\def\csname LTa\endcsname{\color{black}}%
      \expandafter\def\csname LT0\endcsname{\color{black}}%
      \expandafter\def\csname LT1\endcsname{\color{black}}%
      \expandafter\def\csname LT2\endcsname{\color{black}}%
      \expandafter\def\csname LT3\endcsname{\color{black}}%
      \expandafter\def\csname LT4\endcsname{\color{black}}%
      \expandafter\def\csname LT5\endcsname{\color{black}}%
      \expandafter\def\csname LT6\endcsname{\color{black}}%
      \expandafter\def\csname LT7\endcsname{\color{black}}%
      \expandafter\def\csname LT8\endcsname{\color{black}}%
    \fi
  \fi
    \setlength{\unitlength}{0.0500bp}%
    \ifx\gptboxheight\undefined%
      \newlength{\gptboxheight}%
      \newlength{\gptboxwidth}%
      \newsavebox{\gptboxtext}%
    \fi%
    \setlength{\fboxrule}{0.5pt}%
    \setlength{\fboxsep}{1pt}%
\begin{picture}(7200.00,5040.00)%
    \gplgaddtomacro\gplbacktext{%
      \csname LTb\endcsname
      \put(982,704){\makebox(0,0)[r]{\strut{}$-0.0065$}}%
      \put(982,1232){\makebox(0,0)[r]{\strut{}$-0.0060$}}%
      \put(982,1760){\makebox(0,0)[r]{\strut{}$-0.0055$}}%
      \put(982,2288){\makebox(0,0)[r]{\strut{}$-0.0050$}}%
      \put(982,2817){\makebox(0,0)[r]{\strut{}$-0.0045$}}%
      \put(982,3345){\makebox(0,0)[r]{\strut{}$-0.0040$}}%
      \put(982,3873){\makebox(0,0)[r]{\strut{}$-0.0035$}}%
      \put(982,4401){\makebox(0,0)[r]{\strut{}$-0.0030$}}%
      \put(982,4929){\makebox(0,0)[r]{\strut{}$-0.0025$}}%
      \put(1114,484){\makebox(0,0){\strut{}$6$}}%
      \put(1735,484){\makebox(0,0){\strut{}}}%
      \put(2357,484){\makebox(0,0){\strut{}$7$}}%
      \put(2978,484){\makebox(0,0){\strut{}}}%
      \put(3600,484){\makebox(0,0){\strut{}$8$}}%
      \put(4221,484){\makebox(0,0){\strut{}}}%
      \put(4842,484){\makebox(0,0){\strut{}$9$}}%
      \put(5464,484){\makebox(0,0){\strut{}}}%
      \put(6085,484){\makebox(0,0){\strut{}$10$}}%
    }%
    \gplgaddtomacro\gplfronttext{%
      \csname LTb\endcsname
      \put(-162,2816){\rotatebox{-270}{\makebox(0,0){\strut{}\large $\lbar{P}_{11}/\gamma_m$}}}%
      \put(3599,154){\makebox(0,0){\strut{}\large $b/\text{mm}$}}%
      \csname LTb\endcsname
      \put(5098,1960){\makebox(0,0)[r]{\strut{}$\beta=45^\circ$}}%
      \csname LTb\endcsname
      \put(5098,1619){\makebox(0,0)[r]{\strut{}$\beta=60^\circ$}}%
      \csname LTb\endcsname
      \put(5098,1278){\makebox(0,0)[r]{\strut{}$\beta=75^\circ$}}%
      \csname LTb\endcsname
      \put(5098,937){\makebox(0,0)[r]{\strut{}$\Bk\to\Bzero$}}%
    }%
    \gplbacktext
    \psfrag{b}      [cb][cb]{$b$}
    \psfrag{alpha}  [lb][lb]{$\beta$}    
    \put(0,0){\includegraphics{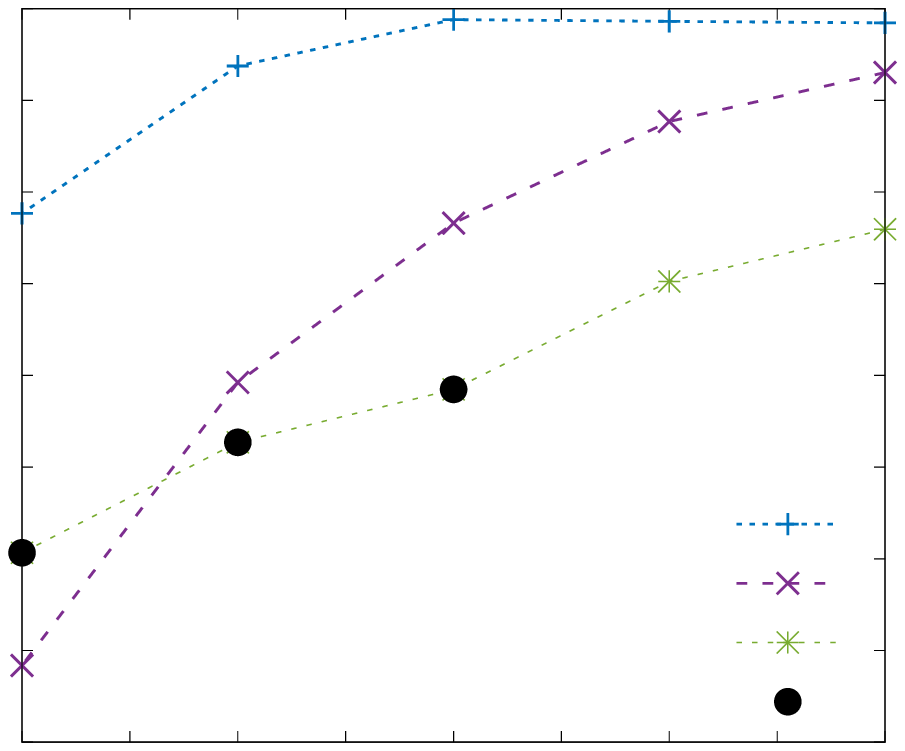}
    \put(-11100,3900){\makebox{
      \includegraphics[height=1.5cm]{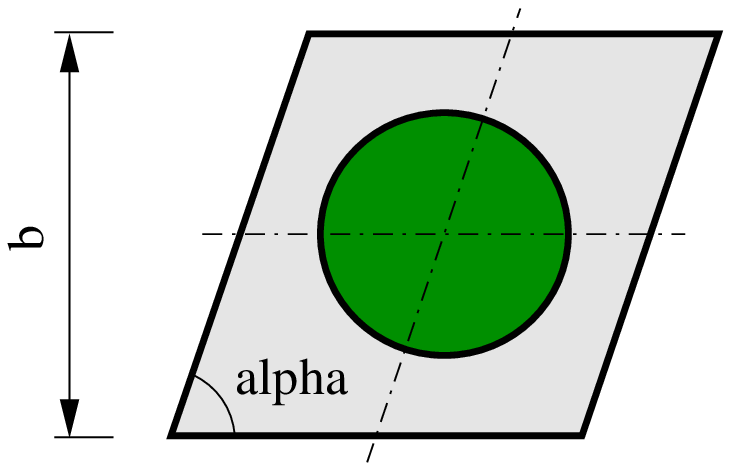}%
    }}
    \put(-4500,1200){\makebox{
      \includegraphics[height=1.5cm]{./fluid-bc3d_inset}%
    }}
}%
    \gplfronttext
  \end{picture}%
\endgroup